\documentclass{elsarticle}
\usepackage{latexsym,amssymb,amsmath}
\usepackage{eucal}
\usepackage{mathrsfs}
\usepackage[all]{xy}
\usepackage{color}
\usepackage{tikz}

\newtheorem{thm}{Theorem}[section]
\newtheorem{lem}[thm]{Lemma}
\newtheorem{prop}[thm]{Proposition}
\newtheorem{cor}[thm]{Corollary}

\newdefinition{defn}[thm]{Definition}
\newdefinition{rmk}[thm]{Remark}
\newdefinition{ex}[thm]{Example}
\newdefinition{exs}[thm]{Examples}
\newdefinition{open}[thm]{Open Question}
\newdefinition{question}[thm]{Question}
\newdefinition{notn}[thm]{Notation}
\newdefinition{ackn}[thm]{Acknowledgments}
\newdefinition{org}[thm]{Organization of the article}
\newdefinition{rel}[thm]{Related work}

\newproof{proof}{Proof}

\DeclareMathAlphabet{\mathbfit}{OT1}{cmr}{bx}{it}

\newcommand{\cat}[1]{\mathbf{{#1}}}

\newcommand{\op}[1]{\mathscr{#1}}

\newcommand{\colim}{\operatorname{colim}}
\newcommand{\N}{\mathbf{N}}

\newcommand{\Zmod}[1]{\mathbf{F}_{#1}}
\newcommand{\Hom}{\operatorname{Hom}}

\newcommand{\End}{\operatorname{End}}

\newcommand{\alg}[1]{{\op {#1}}\text{-}\cat{Alg}}

\newcommand{\coalg}[1]{{\op{#1}}\text{-}\cat{Coalg}}

\newcommand{\Ssigma}{\op{A}^\perp}
\newcommand{\diffract}{\Phi}

\newcommand{\bimod}[2]{{}_{\op{#1}}\cat{Mod}_{\op{#2}}}
\newcommand{\bicomod}[2]{{}_{{#1}}\cat{Comod}_{{#2}}}

\newcommand{\acirc}{\underset {\op A}{\circ}}
\newcommand{\ind}{\operatorname{Ind}}
\newcommand{\lin}{\operatorname{Lin}}
\newcommand\Om{\Omega}
\newcommand{\im}{\operatorname{Im}}
\newcommand{\ob}{\operatorname{Ob}}

\newcommand{\del}{\partial}
\newcommand{\susp}{\mathfrak S}

\newcommand\egal[2]{\overset {#1}{\underset {#2}\rightrightarrows }}
\newcommand\adjunct[2]{\overset {#1}{\underset {#2}\rightleftarrows }}
\newcommand{\opbar}{\mathbf{B}}
\newcommand{\ve}{\varepsilon}

\newcommand{\dual}{\sharp}
\newcommand{\sgn}{\mathrm{sgn}}
\newcommand{\monad}[1]{\mathsf{#1}}
\newcommand{\coker}{\operatorname{coker}}
\newcommand{\Mod}[2]{{}_{\op{#1}}\cat{Mod}_{\op{#2}}}
\newcommand{\Comod}[2]{{}_{\op{#1}}\cat{Comod}_{\op{#2}}}

\newcommand{\In}{\operatorname{In}}
\newcommand{\Out}{\operatorname{Out}}

\usepackage{hyperref}

\title{Twisting structures and strongly homotopy morphisms\tnoteref{pep}}

\author{Kathryn Hess\fnref{kh}}
\author{Jonathan Scott}

\address{Institut de g\'eom\'etrie, alg\`ebre et topologie (IGAT) \\
    \'Ecole Polytechnique F\'ed\'erale de Lausanne \\
    CH-1015 Lausanne \\
    Switzerland}
    \ead{kathryn.hess@epfl.ch}

\address{Department of Mathematics \\ Cleveland State University \\ 2121 Euclid Ave., RT 1515 \\ Cleveland OH 44115-2214 USA }
    \ead{j.a.scott3@csuohio.edu}

\address{Department of Mathematics and Statistics \\ University of Ottawa \\ 585 King Edward Ave. \\ Ottawa ON K1N 6N5 Canada}
    \ead{pparent@uottawa.ca}
    
\tnotetext[pep]{With an appendix by P.-E. Parent (University of Ottawa)}

\fntext[kh]{K.H. thanks the Institut Mittag-Leffler (Djursholm,
Sweden) for its hospitality during a crucial phase of this research and the Midwest Topology Network for its financial support.}

    \date{\today}

\begin{document}

\begin{abstract}
{In an application of} the notion of twisting structures {introduced by Hess and Lack \cite{hess-lack}, we} define twisted composition products {of symmetric sequences of chain complexes that are degreewise projective and finitely generated.  Let $\op Q$ be a cooperad, and let $\opbar \op P$ denote the bar construction on an operad $\op P$.} To each {morphism of cooperads} $g:\op{Q} \rightarrow \opbar\op{P}$ is associated a $\op{P}$-co-ring, $K(g)$, which generalizes the two-sided Koszul and bar constructions.  When the counit $K(g) \rightarrow \op{P}$ is a quasi-isomorphism, we show that the Kleisli category for $K(g)$ is isomorphic to the category of $\op{P}$-algebras and {of their} morphisms up to  strong homotopy{, and} we give the classifying morphisms for both strict and homotopy $\op{P}$-algebras.  Parametrized morphisms {of (co)associative chain (co)algebras} up to strong homotopy are also introduced and studied, {and a general existence theorem is proved}.  {In the appendix, we study the co-ring associated to the canonical morphism of cooperads $\op A^\perp\to \opbar \op A$, which is exactly the two-sided Koszul resolution of the associative operad $\op A$, also known as the Alexander-Whitney co-ring.}
\end{abstract}


\maketitle

 \section{Introduction}

When working with a given type of algebra, it is often possible to form a ``standard construction'', a co-free coalgebra of some sort in which the algebraic structure is encoded in a quadratic differential.  For example, for an associative or Lie algebra we have respectively the bar construction and the Chevalley-Eilenberg complex~\cite{chevalley-eilenberg:48}.  We can furthermore reconstruct the original algebra from its standard construction, at least if $2$ is invertible in the ground ring.  Thus we often identify the algebra and its standard construction.

The family of higher homotopies that constitute a strongly homotopy-multiplicative map $f : A \Rightarrow A'$ can be ``rolled up'' into a single morphism of coassociative coalgebras, $F : BA \rightarrow BA'$~\cite{gugenheim-munkholm:74}.  Similarly, if $C_{*}(-)$ denotes the Chevalley-Eilenberg complex of a DG Lie algebra, which is a commutative associative coalgebra, then a morphism of coalgebras, $F: C_{*}(L) \rightarrow C_{*}(L')$ `is' a strong homotopy Lie morphism.

More generally, let $\op{P}$ be a quadratic operad with quadratic dual cooperad $\op{P}^{\bot}$~\cite{getzler-jones:94}.  A $\op{P}$-algebra $A$ has a standard `bar' construction $B_{\op{P}}A$ that is a co-free connected $\op{P}^{\bot}$-coalgebra, provided with a quadratic differential that encodes the algebraic structure of $A$. A strong homotopy (henceforth abbreviated SH) $\op{P}$-morphism from $A$ to another $\op{P}$-algebra $A'$ is then a morphism of $\op{P}^{\bot}$-coalgebras, $B_{\op{P}}A \rightarrow B_{\op{P}}A'$.  We remark that since we are using $\op{P}^{\bot}$ instead of the quadratic dual \emph{operad} of Ginzburg and Kapranov~\cite{ginzburg-kapranov:94}, our bar constructions are shifted by one degree from the classical constructions.

Dually, there is a standard `cobar' construction for $\op{P}$-coalgebras, that takes values in $\op{P}^{!}$-algebras, and we obtain the notion of SH morphisms of $\op{P}$-coalgebras.

The goal of this article is to describe the various categories with SH morphisms ``operadically'', when working in the category $\cat {dgProj}$  of degreewise finitely generated and projective chain complexes over a commutative ring $R$. Our first result reads as follows (Proposition \ref{prop:quad-co-ring} and Theorem \ref{thm:kleisli}).

\begin{thm}
Let $\op{P}$ be a Koszul operad in $\cat {dgProj}$.  The two-sided Koszul resolution $K(\op{P})$ of $\op{P}$ is naturally a $\op{P}$-co-ring, and therefore induces a comonad on the category of $\op P$-algebras in $\cat {dgProj}$.  The associated Kleisli category, with $\op{P}$-algebras as objects and morphisms parametrized by $K(\op{P})$, is isomorphic to the category of $\op{P}$-algebras in $\cat {dgProj}$ and their SH morphisms.
\end{thm}

A similar result for $\op{P}^{\bot}$-coalgebras holds as well.

Recall that a $\op{P}_{\infty}$-algebra is an algebra over a cofibrant replacement of $\op{P}$, with respect to some reasonable model structure on the category of operads~\cite{berger-moerdijk:03}.  If $\op{P}$ is Koszul, then we take as our replacement the operad $\mathbf{\Omega} \op{P}^{\bot}$~\cite{ginzburg-kapranov:94}, where $\mathbf{\Omega}$ denotes the operadic cobar construction.  We construct a Koszul resolution $K'(\op{P})$ for $\mathbf{\Omega} \op{P}^{\bot}$, and obtain the following result (Theorem \ref{thm:p-infty}).

\begin{thm}
With the notation and hypotheses above, the resolution $K'(\op{P})$ is naturally a $\mathbf{\Omega}\op{P}^{\bot}$-co-ring {and therefore induces a comonad on the category of $\mathbf{\Omega}\op{P}^{\bot}$-algebras} in $\cat {dgProj}$.  The associated Kleisli category, with $\mathbf{\Omega}\op{P}^{\bot}$-algebras as objects, and morphisms {parametrized by $K'(\op{P})$}, is isomorphic to the category of $\op{P}_{\infty}$-algebras in $\cat {dgProj}$ and their SH morphisms.
\end{thm}

We consider moreover parametrization of SH morphisms of associative algebras (respectively, coassociative coalgebras)  over cooperads of chain coalgebras (respectively, over operads of chain algebras) (Definition \ref{defn:parameter}).  We prove that such parametrized SH categories also admit an operadic, (co)Kleisli description (Theorem \ref{thm:kleisli-x}), which then enables us to establish useful existence results for natural, parametrized SH structures (Theorems \ref{thm:x-exist-bar} and \ref{thm:x-exist-cobar}).

In Section~\ref{sec:twisting}, we present our main tool in the study of morphism sets: the notion of \emph{twisting structures}.  While twisting cochains in the category of differential graded symmetric sequences go back as early~\cite{getzler-jones:94}, twisting structures on a category  allow for the definition of twisted products via classifying morphisms rather than twisting cochains, and can therefore be applied more generally.  We define twisting structures and twisted products, prove the expected adjunction relations, and define the standard (dual) constructions associated to a classifying morphism.  One result of the adjunction relations turns out to be critical to our development: a standard construction associated to any classifying map admits a natural comonoidal structure.

In Section~\ref{sec:operad}, we review differential graded symmetric sequences and (co)operads.  We use the bar construction with coefficients to define a twisting structure on the category of symmetric sequences of chain complexes.

In Section~\ref{sec:shp}, we recall the dual cooperad $\op{P}^{\bot}$ of a weight-graded operad $\op{P}$; the inclusion $\kappa_{\op P}:\op{P}^{\bot} \rightarrow \opbar\op{P}$ is a classifying morphism.  We show that if $\op{P}$ is Koszul (i.e., the  inclusion $\kappa_{\op P}$ is a quasi-isomorphism), then the Kleisli category for the standard construction $K(\kappa_{\op P})$ is isomorphic to the category of $\op{P}$-algebras and SH morphisms.  Considering instead the unit of the adjunction, $\eta: \op{P}^{\bot} \rightarrow \opbar\mathbf{\Omega}\op{P}^\bot$, we obtain the $\op{P}_{\infty}$ category.  We then briefly discuss the general, non-Koszul, case, where we use the identity on $\opbar\op{P}$ and the unit $\opbar\op{P} \rightarrow \opbar\mathbf{\Omega} \opbar\op{P}$, respectively, for our classifying morphisms.

In Section~\ref{sec:delooping}, we introduce and discuss the concept of parametrized SH morphisms of (co)associative coalgebras.

In the first appendix, P.-E. Parent discusses the special case of the Koszul resolution of the associative operad, $\op{A}$, and shows that an algebra with an SH-comultiplicative diagonal is not the same thing as an algebra with cup-$i$ products, as one might have expected.

The second appendix is devoted to the somewhat technical proof of the existence of parametrized SH morphisms of (co)associative (co)algebras.

\subsection{Comparison to other approaches}

The original work in this direction, done in the category of spaces, was by Iwase and Mimura~\cite{iwase-mimura:89}, who considered $A_{\infty}$-maps of $A_{\infty}$-spaces.

Markl~\cite{markl:04} considered the coloured operad, or multicategory, $\op{P}(\cdot \rightarrow \cdot)$, whose algebras are morphisms of $\op{P}$-algebras.  An SH $\op{P}$-algebra morphism is then an algebra for a cofibrant replacement of $\op{P}(\cdot \rightarrow \cdot)$; Markl shows, among other things, that SH morphisms are homotopy invariant.  We can relate our co-rings to the multicategory approach, and answer Markl's Problem 26 of~\cite{markl:02}: as described below, we can construct a map $\op{P}'(a \rightarrow c) \rightarrow \op{P}'(b \rightarrow c) * \op{P}'(a \rightarrow b)$ (where $*$ is the free product) such that the resulting composition of SH morphisms is associative.

Let $\op{P}$ be an operad, and let $K$ be the $\mathbf{\Omega} \opbar\op{P}$-co-ring associated to the classifying morphism $\opbar\op{P} \rightarrow \opbar\mathbf{\Omega}\opbar\op{P}$.     A category enriched in symmetric sequences is a special case of a multicategory.  Consider two such categories: the first, $\cat{M}$, has two objects $x$ and $y$, with hom objects $\cat{M}(x,x) = \cat{M}(y,y) = \op{P}'$, $\cat{M}(x,y) = K(\op{P}')$.   The second, $\cat{N}$, has three objects, $a$, $b$, and $c$, with hom objects $\cat{N}(a,a) = \cat{N}(b,b) = \cat{N}(c,c) = \op{P}'$, and $\cat{N}(a,b) = \cat{N}(b,c) = K(\op{P}')$, $\cat{N}(a,c) = K(\op{P}') \circ_{\op{P}'} K(\op{P}')$.   An enriched functor $F: \cat{M} \rightarrow \cat{N}$ is defined by $F(x) = a$, $F(y) = c$, and $F:\cat{M}(x,y) \rightarrow \cat{N}(a,c)$ is the diagonal in $K(\op{P}')$.  Then pulling back along $F$ provides the desired associative composition in the SH category, from the multicategory point of view.  Note that in our formulation of homotopy morphisms, it is obvious how to define their composition and indeed the whole category structure.

Coloured operads have the advantage that they can be applied to many more situations than can homotopy morphisms.  Berger and Moerdijk~\cite{berger-moerdijk:07} consider coloured operads in general, and study conditions under which, in particular, SH morphisms are rectifiable.

Leinster~\cite{leinster:00} defines the category of homotopy algebras over an operad $\op{P}$ to be the category of colax representations, $X : \widehat{\op{P}} \rightarrow \cat{M}$, where $\widehat{\op{P}}$ is the strict monoidal category with the nonnegative integers as objects, $+$ and $0$ for the monoidal structure, and $\widehat{\op{P}}(m,n) = \op{P}^{\odot m}(n)$ (see Section~\ref{subsec:sym-seq} for the graded tensor product $\odot$ of symmetric sequences).  The morphisms are monoidal transformations.  In this formulation, morphisms do indeed commute with algebraic structure weakly, but homotopy invariance is not immediately obvious to us.

\subsection{Acknowledgments}  The authors would like to express their heartfelt appreciation to Paul-Eug\`ene Parent for his helpful participation in the early stages of this project and in particular for his contribution to our understanding of the Alexander-Whitney co-ring, as detailed in the first appendix to this article.

This has proved to be a very long-term project, which has evolved significantly over the past six years.   Earlier versions of our approach to describing strongly homotopy morphisms via co-rings can be found on the arXiv~\cite{hps:05}.  Results from these earlier manuscripts, in particular concerning the Alexander-Whitney co-ring, which have already been applied in various articles and theses (e.g., \cite{hpst:canonicalAH}, \cite{hps:cohoch}, \cite{hess-levi:07}, \cite{hess-rognes}, \cite {naito} and \cite{boyle}), are also stated and proved here.

\subsection{Notation}\label{notn:intro}

 If $\cat C$ is a category, and $A$ and $B$ are objects in $\cat C$, then  $\cat C(A,B)$ denotes the set of morphisms from $A$ to $B$.

If $\mathsf T=(T,\mu, \eta)$ is a monad on a category $\cat C$, then $\cat C_{\mathsf T}$ denotes the Kleisli category determined by $\mathsf T$, with $\ob \cat C_{\mathsf T}=\ob \cat C$ and $\cat C_{\mathsf T}(A,B)= \cat C(A,TB)$. If $f\in \cat C_{\mathsf T}(A,B)$ and $g\in \cat C_{\mathsf T}(B,C)$, then their composite in $\cat C_{\mathsf T}$ is defined to be the composite of
 $$A\xrightarrow f TB\xrightarrow {Tg} T^2C \xrightarrow {\mu _{C}} TC$$
 in $\cat C$.

 Dually, if $\mathsf K=(K, \Delta, \ve)$ is a comonad on $\cat C$, then ${}_{\mathsf K}\cat C$ denotes the coKleisli category determined by $\mathsf K$, with $\ob {}_{\mathsf K}\cat C=\ob \cat C$ and ${}_{\mathsf K}\cat C(A,B)= \cat C(KA,B)$. If $f\in {}_{\mathsf K}\cat C(A,B)$ and $g\in {}_{\mathsf K}\cat C(B,C)$, then their composite in ${}_{\mathsf K}\cat C$ is defined to be the composite of
 $$KA\xrightarrow {\Delta_{K}} K^2A\xrightarrow {Kf} KB \xrightarrow {g} C$$
 in $\cat C$.

 \section{Twisting structures}\label{sec:twisting}

The goal of this section is to introduce a categorical structure that conveniently generalizes both  twisting cochains from differential graded coalgebras to differential graded algebras and twisting functions from simplicial sets to simplicial groups.  Hess and Lack first formulated such a definition in \cite{hess-lack}, though in an even more highly categorical manner.

Throughout this section $(\cat M, \otimes, I)$ denotes a monoidal category that admits equalizers and coequalizers.  Let $\cat {Mon}$ and $\cat {Comon}$ denote the categories of monoids and of comonoids in $\cat M$.  If $A$ is a monoid in $\cat M$, then ${}_{A}\cat {Mod}$ and $\cat {Mod}_{A}$ are the categories of left $A$-modules and of right $A$-modules.  Similarly, ${}_{C}\cat {Comod}$ and $\cat {Comod}_{C}$ denote the categories of left and right comodules over a comonoid $C$.

The following definitions are classical.

\begin{defn}  Let $A$ be a monoid in $\cat M$.  Let $(M,\rho)$ be a right $A$-module, and let $(N,\lambda)$ be a left $A$-module.  The \emph{tensor product of $M$ and $N$ over $A$} is the coequalizer
 $$M\underset A\otimes N:=\operatorname{coequal}(M\otimes A\otimes N \egal{\rho\otimes N}{M\otimes \lambda} M\otimes N).$$

 Let $C$ be a comonoid in $\cat M$.  Let $(M,\rho)$ be a right $C$-comodule, and let $(N,\lambda)$ be a left $C$-comodule.  The \emph{cotensor product of $M$ and $N$ over $C$} is the equalizer
 $$M\underset C\square N:=\operatorname{equal}(M\otimes N\egal{\rho\otimes N}{M\otimes \lambda} M\otimes C\otimes N).$$
 \end{defn}

We often consider objects of $\cat M$ that are endowed with either two actions or two coactions or an action and a coaction, for which we introduce the following notation.

\begin{notn} Let $A$ and $A'$ be monoids in $\cat M$, and let $C$ and $C'$ be comonoids in $\cat M$.   We consider the following classes of objects in $\cat M$ that are endowed with two structures.
\medskip
\begin{align*}
{}_{A}\cat {Mix}_{A'}=&{}_{A}\cat {Mod}_{A'}\\
=&\{ (M, \rho, \lambda)\mid (M,\rho)\in \cat {Mod}_{A'}, (M,\lambda)\in {}_{A}\cat {Mod}, \lambda (A\otimes \rho)=\rho(\lambda \otimes A')\}
\end{align*}
\medskip
$${}_{A}\cat {Mix}_{C}=\{ (M, \rho, \lambda)\mid (M,\rho)\in \cat {Comod}_{C}, (M,\lambda)\in {}_{A}\cat {Mod}, (\lambda\otimes C) (A\otimes \rho)=\rho\lambda\}$$
\medskip
$${}_{C}\cat {Mix}_{A}=\{ (M, \rho, \lambda)\mid (M,\rho)\in \cat {Mod}_{A}, (M,\lambda)\in {}_{C}\cat {Comod}, (C\otimes \rho) (\lambda \otimes A)=\lambda\rho\}$$
\medskip
\begin{align*}
{}_{C}\cat {Mix}_{C'}=&{}_{C}\cat {Comod}_{C'}\\
&=\{ (M, \rho, \lambda)\mid (M,\rho)\in \cat {Comod}_{C'}, (M,\lambda)\in {}_{C}\cat {Comod}, (\lambda \otimes C')\rho=(C\otimes \rho)\lambda\}
\end{align*}
\end{notn}

Tensor and cotensor products must commute, in the sense of the following definition, if we wish to define twisting structures.

 \begin{defn} The monoidal category $\cat M$ is \emph{twistable} if the tensor and cotensor products defined above restrict and corestrict to {bifunctors}
 $$-\underset {A}\otimes -:{}_{X}\cat {Mix}_{A}\times {}_{A}\cat {Mix}_{Y}\to {}_{X}\cat {Mix}_{Y}$$
 and
 $$-\underset C\square -: {}_{X}\cat {Mix}_{C}\times {}_{C}\cat {Mix}_{Y}\to {}_{X}\cat {Mix}_{Y},$$
for all (co)monoids $X$ and $Y$.  Furthermore, the tensoring and cotensoring must be associative up to isomorphism in the obvious sense.
\end{defn}

\begin{ex}  The categories of sets and of simplicial sets are clearly twistable, as is the category $\cat{dgProj}$ of degree-wise finitely generated, projective chain complexes over a commutative ring $R$.  We prove in section \ref{sssec:twistability} that the category of symmetric sequences in $\cat{dgProj}$ is also twistable (Theorem \ref{thm:twistable}).
\end{ex}

\begin{rmk}  If $\cat M$ is twistable, then we can define a category $\cat {Mix}$ with
$$\ob \cat {Mix}=\ob \cat {Mon} \cup \ob \cat {Comon}$$
and
$$\cat {Mix}(X,Y)={}_{X}\cat {Mix}_{Y}/\cong,$$
where $\cong$ denotes the isomorphism relation.  Composition of morphisms in $\cat {Mix}$ is given by tensoring over monoids and cotensoring over comonoids.

There are functors
$$J: \cat {Mon} \to \cat {Mix}\quad \text{ and }\quad \widetilde J: \cat {Mon}^{op} \to \cat {Mix}$$
specified on objects by $J(A)=A=\widetilde J(A)$.  On morphisms we set
$$J(f:A\to A')= {}_{f}A'\in \cat {Mix}(A,A'),$$ where ${}_{f}A'$ denotes the isomorphism class of $A'$ seen as a right  $A'$-module via its own multiplication and as left $A$-module via $f$, while
$$\widetilde J(f:A\to A')= A'_{f}\in \cat {Mix}(A',A),$$
where $A'_{f}$ denotes the isomorphism class of $A'$ seen as a left $A'$-module via its own multiplication and as a right $A$-module via $f$.

Similarly, there are functors
$$coJ: \cat {Comon} \to \cat {Mix}\quad \text {  and  } \quad co\widetilde J:\cat{Comon}^{op}\to \cat {Mix}$$
specified on objects by $coJ(C)=C=co\widetilde J(C)$.  On morphisms we set
$$coJ(g:C\to C')= C_{g}\in \cat{Mix}(C,C'),$$
where $C_{g}$ denotes the isomorphism class of $C$ seen as a left  $C$-comodule via its own comultiplication and as right $C$-comodule via $g$, while
$$co\widetilde J(g:C\to C')= {}_{g}C\in \cat {Mix}(C',C),$$ where ${}_{g}C$ denotes the isomorphism class of $C$ seen as a right $C$-comodule via its own multiplication and as a left $C'$-comodule via $f$.
\end {rmk}

The following definition is a slight variant of that formulated in \cite{hess-lack}.

 \begin{defn}\label{defn:twisting}  A \emph{right twisting structure} on a twistable monoidal category $\cat M$ consists of a functor
 $$B : \cat {Mon}\to \cat {Comon}$$
 together with natural transformations
 $$E:coJ\circ B\Rightarrow J$$
 and
$$\widetilde E: \widetilde J \Rightarrow co\widetilde J\circ B$$
 of functors from $\cat {Mon}$ (respectively, $\cat{Mon}^{op}$) to $\cat {Mix}$ and natural morphisms for all monoids $A$
 $$\delta_{A}: B A\to EA\underset A\otimes \widetilde EA$$
 of $B A$-bicomodules
 and
 $$\mu_{A}: \widetilde EA \underset {B A} \square EA\to A$$
 of $A$-bimodules such that the following diagrams commute.
 $$\xymatrix{\widetilde EA \cong  \widetilde E A\underset {B A}\square B A\ar [rr]^{\widetilde EA \underset {B A}\square \delta_{A}}\ar@{=}[drr]&&\widetilde EA \underset {B A}\square EA\underset A\otimes \widetilde EA \ar[d]^{\mu _{A}\underset A\otimes \widetilde EA}\\
 &&\widetilde EA}$$
  $$\xymatrix{ EA \cong  B A\underset {B A}\square E A\ar [rr]^{ \delta_{A} \underset {B A}\square E A}\ar@{=}[drr]&&EA\underset A\otimes \widetilde EA \underset {B A}\square EA\ar[d]^{EA\underset A\otimes \mu _{A}}\\
 && EA}$$
  \end{defn}

The choice of terminology above is motivated by the existence of right twisting structures when $\cat M$ is the category of either connected chain complexes or reduced simplicial sets, where the right twisting structure can be defined by either twisting cochains or twisting functions, both of which {are} classical notions.  The details of both of these cases can be found in \cite{hess-lack}.  In this paper we generalize the chain complex case in section \ref{subsec:twist-chcx}, showing that the category of symmetric sequences of chain complexes admits a right twisting structure.

\begin{rmk}  It may be helpful to unfold the definition above.  If $(B, E, \widetilde E, \delta, \mu)$ is a right twisting structure on $\cat M$, then for all monoids $A$,
$$EA\in {}_{B A}\cat {Mix}_{A} \quad\text{ and }\quad \widetilde EA\in {}_{A}\cat {Mix}_{B A}.$$
Moreover, for all monoid morphisms $f:A\to A'$, the naturality of $E$ and of $\widetilde E$ implies that
$$EA\underset A\otimes {}_{f}A'\cong B A_{B f}\underset {B A'}\square EA'$$
and
$$A'_{f}\underset A\otimes \widetilde EA\cong \widetilde EA'\underset {B A'}\square {}_{B f}B A.$$
 \end{rmk}

For any monoid $A$, we think of $EA$ and $\widetilde EA$ as  the total spaces of the ``universal  right $A$-bundle'' and the ``universal  left $A$-bundle''  in $\cat M$, respectively.  This vision of their role motivates the following definition.

\begin{defn}\label{defn:classifying-twisting}  Let $C$ be a comonoid, and let $A$ be a monoid in a twistable monoidal category $\cat M$ endowed with a right twisting structure.  A \emph{classifying morphism} between $C$ and $A$ is a morphism of comonoids $g:C\to B A$.

Given a right $C$-comodule $V$ and a left $A$-module $W$, we define the \emph{twisted product} of $V$ and $W$ {with respect to the classifying morphism $g$} as
\[
    V \otimes_{g} W = V_{g} \underset{BA}{\square} EA \underset{A}{\otimes} W,
\]
where $V_{g}$ is $V$ considered as a right $BA$-comodule via $g$.
If $X$ is a right $A$-module and $Y$ is a left $C$-comodule, then
\[
    X \otimes_{g} Y = X \underset{A}{\otimes} \tilde{E}A \underset{BA}{\square} ({}_{g}Y),
\]
where ${}_{g}Y$ is $Y$ considered as a left $BA$-comodule via $g$.

The \emph{right $A$-bundle induced by $g$} is
$C_{g} \otimes_{g} A\in  {}_{C}\cat {Mix}_{A}$,
while the \emph{left $A$-bundle induced by $g$} is
$A \otimes_{g}C\in {}_{A}\cat {Mix}_{C}$.
\end{defn}

The next theorem and its corollary are essential to establishing our operadic characterization of strongly homotopy maps.

 \begin{thm}\label{thm:key-adjunction}
 Let $C$ and $C'$ be comonoids, and let $A$ and $A'$ be monoids  in a twistable monoidal category $\cat M$ endowed with a right twisting structure.  Let $g: C\to B A$ and $g': C'\to B A'$ be classifying morphisms. If
 \[
    (g,g')_{*}:{}_{C}\cat {Comod}_{C'}\to {}_{A}\cat {Mod}_{A'}
 \]
 denotes the functor specified by
 \[
    (g,g')_{*}(N)= A \otimes_{g} N \otimes_{g'} A'
 \]
 for all $(C,C')$-bicomodules $N$, and
 \[
    (g,g')^*:{}_{A}\cat {Mod}_{A'}\to{}_{C}\cat {Comod}_{C'}
 \]
denotes the functor specified by
 \[
    (g,g')^*(M)=C \otimes_{g} M \otimes_{g'} C',
 \]
then the functors
\[
    (g,g')_{*}:{}_{C}\cat {Comod}_{C'}\adjunct{}{}{}_{A}\cat {Mod}_{A'}:(g,g')^*
\]
form an adjunction.
 \end{thm}

 \begin{proof}  It clearly suffices to consider the universal case, i.e., $C=B A$, $C'=B A'$ and $g=Id_{B A}$, $g'=Id_{B A'}$.

 We define a unit natural transformation $\eta: Id \Rightarrow (g,g')^*(g,g')_{*}$ by
 $$\eta_{N}=\delta_{A}\underset{B A}\square N\underset {B A'}\square {\delta_{A'}}: N\to  EA\underset A\otimes \widetilde EA\underset{B A}\square N\underset {B A'}\square EA'\underset {A'}\otimes \widetilde EA'$$
 and a counit natural transformation $\ve: (g,g')_{*}(g,g')^*\Rightarrow Id$ by
 $$\ve_{M}=\mu_{A}\underset A\otimes M\underset {A'} \otimes\mu_{A'}:\widetilde EA \underset {B A} \square EA \underset A\otimes M\underset {A'} \otimes \widetilde EA' \underset {B A'} \square EA'\to M.$$
Since $(\mu _{A}\underset A\otimes \widetilde EA)(\widetilde EA \underset {B A}\square \delta_{A})=Id_{\widetilde EA}$ and $(EA\underset A\otimes \mu _{A})(\delta_{A} \underset {B A}\square E A)=Id_{EA}$, we can show easily that
 $$\ve_{(g,g)_{*}}\circ (g,g')_{*}\eta$$
 is the identity natural transformation on $(g,g')_{*}$ and that
 $$(g,g')^*\ve\circ \eta_{(g,g')^*}$$
 is the identity natural transformation on $(g,g')^*$.  In other words, $\eta$ and $\ve$ are the unit and counit of an adjunction.
 \end{proof}

Specializing to $g=Id_{I}$ or $g'=Id_{I}$, we obtain the next result.

 \begin{cor}\label{cor:adjunct-oneside}
 Let $C$ and $A$ be comonoid and a monoid in a twistable monoidal category $\cat M$ endowed with a right twisting structure.  If $g: C\to B A$ is a classifying morphism, then there are adjunctions
 \[
    (g_{\ell})_{*}:{}_{C}\cat {Comod}\adjunct{}{} {}_{A}\cat {Mod}:(g_{\ell})^*
\]
and
\[
    (g_{r})_{*}:\cat {Comod}_{C}\adjunct{}{}\cat {Mod}_{A}:(g_{r}),^*
\]
where
\[
    (g_{\ell})_{*}(N)= A \otimes_{g}N\quad\text{and}\quad (g_{r})_{*}(N')= N' \otimes_{g} A
\]
for all left $C$-comodules $N$ and right $C$-comodules $N'$, and
\[
    (g_{\ell})^*(M)=C \otimes_{g} M \quad\text{and}\quad (g_{r})^*(M')=M'\otimes_{g}C),
\]
 for all left $A$-modules $M$ and right $A$-comodules $M'$.
 \end{cor}

The following proposition plays a critical role in the rest of this paper.  Recall that if $A$ is a monoid in a monoidal category $\cat M$, then an \emph{$A$-co-ring} is a comonoid in the bimodule category $(\bimod AA, \underset A\otimes)$.  Dually, if $C$ is a comonoid in $\cat M$, then a \emph{$C$-ring} is monoid in the bicomodule category $(\bicomod CC, \underset C\square)$.

  \begin{prop}\label{prop:param-ring}
 Let $(B , E, \widetilde E, \delta, \mu)$ be a right twisting structure on a twistable monoidal category $\cat M$.  Let $g:C \rightarrow BA$ be a classifying morphism.
 \begin{enumerate}
 \item If $N$ is a comonoid in $({}_{C}\cat {Comod}_{C},\underset{C}\square, C)$, then  $A \otimes_{g} N \otimes_{g} A$ admits a natural $A$-co-ring structure.
 \medskip

\item If $M$ is a monoid in $({}_{A}\cat {Mod}_{A},\underset{A}\otimes,  A)$, then  $C \otimes_{g} M \otimes_{g} C$ admits a natural $C$-ring structure.
\end{enumerate}
 \end{prop}

 \begin{proof} (1) Since $N$ is a comonoid, $N \underset{C}{\square}(-)$ is the underlying functor of a comonad on ${}_{C}\cat{Comod}$.  Since $(A \otimes_{g} (-),C \otimes_{g}(-))$ form an adjoint pair, $A \otimes_{g} N \underset{C}{\square} C \otimes_{g}(-) = A \otimes_{g} N \otimes_{g}(-)$ is {also} the underlying functor of a comonad.  It follows that $A \otimes_{g} N \otimes_{g} A$ is an $A$-co-ring, with structure maps coming from the comultiplication and counit of the comonad applied to $A$.

 The proof of (2) is dual to that of (1).
\end{proof}

Considering $C$ as a bicomodule over itself and $A$ as a bimodule over itself, in the obvious way, we obtain important special cases of the constructions considered in the proposition above.

\begin{defn}\label{defn:std-constr}
Let $g:C\rightarrow BA$ be a classifying morphism.
Treating $A$ as a bimodule over itself, and $C$ as a bicomodule over itself, we obtain the \emph{standard construction {on $g$}}
\[
    K(g) = A \otimes_{g} C \otimes_{g} A.
\]
and the \emph{dual standard construction {on $g$}}
\[
    T(g) = C \otimes_{g} A \otimes_{g} C.
\]
\end{defn}

Observe that $C$ is a comonoid in $({}_{C}\cat {Comod}_{C},\underset{C}\square, C)$, with the unique isomorphism $C\xrightarrow\cong C\underset C\square C$ as comultiplication.  Similarly, $A$ is a monoid in $({}_{A}\cat {Mod}_{A},\underset{A}\otimes,  A)$, where the multiplication is $A\underset A\otimes A \xrightarrow \cong A$.  The next proposition is therefore an immediate consequence of  Proposition \ref{prop:param-ring}.

\begin{prop}\label{prop:std-coring}
{For any classifying morphism $g: C\to BA$, the standard construction $K(g)$ is an $A$-co-ring, and the dual standard construction $T(g)$ is a $C$-ring; these structures are natural in $A$ and in $C$.}
\end{prop}

For future reference, we note that any twisted product with respect to a classifying morphism $g$ can be computed in terms of the standard construction $K(g)$ and the dual standard construction $T(g)$. The proof consists of straightforward calculation.

\begin{prop}
Let $g : C \rightarrow BA$ be a classifying morphism.
\begin{enumerate}
\item If $M$ is a right $A$-module and $N$ is a left $A$-module, then $$M \underset A\otimes K(g) \underset A\otimes N \cong M \otimes_{g} C \otimes_{g} N.$$
\item If $U$ is a right $C$-comodule and $V$ is a left $C$-comodule, then $$U \underset{C}{\square} T(g) \underset{C}{\square} V \cong U \otimes_{g} T(g) \otimes_{g} V.$$
\end{enumerate}
The  isomorphisms {above} are natural in all variables.
\end{prop}

 \begin{rmk}  There is a strictly dual notion of left twisting structures $(\Om, P, \widetilde P, \mu, \delta)$ on twistable model categories.  If $\cat M$ admits both left and right twisting structures, and $(\Om , B)$ is an adjoint pair of functors, then the natural transformations $P$ and $\widetilde P$ can be deduced from $E$ and $\widetilde E$, as the right and left $\Om C$-bundles induced by the unit $\eta: C\to B \Om C$ of the $(\Om, B)$-adjunction, i.e.,
$$PC:=C_{\eta}\underset {B \Om C}\square E\Om C\in {}_{C}\cat {Mix}_{\Om C}$$
and
$$\widetilde PC:=\widetilde E\Om C \underset {B \Om C}\square {}_{\eta}C\in {}_{\Om C}\cat {Mix}_{C}.$$
The natural transformations $E$ and $\widetilde E$ can be similarly deduced from $P$ and $\widetilde P$, using the counit of the  $(\Om, B )$-adjunction.  Moreover, given classifying morphisms $g:C\to B A$ and $g':C'\to B A'$, an adjunction
 $${}_{C}\cat {Comod}_{C'}\adjunct{}{}{}_{A}\cat {Mod}_{A'}$$
 can be constructed using the transpose $g^\flat : \Om C \to A$ of $g$, as well as $PC$ and $\widetilde PC$.
\end{rmk}

\section{Symmetric sequences of chain complexes}\label{sec:operad}

In this section, we recall the definition of the monoidal category of symmetric sequences with the composition product.  Operads are monoids in this category; with some hypotheses, we can consider cooperads to be comonoids.  The functors $\mathbf E$ and $\widetilde{\mathbf E}$ are the two acyclic bar constructions as defined in~\cite{ginzburg-kapranov:94,fresse:04}; we recall the definitions, define the morphisms necessary for our twisting structures, and verify that the various identities hold.

\subsection{Symmetric sequences}\label{subsec:sym-seq}

We work in the closed symmetric monoidal category $\cat{dgM}$ of differential graded (DG) modules over an arbitrary commutative ring,
$R$, furnished with the (graded) tensor product $\otimes = \otimes_{R}$. The graded hom functor $\Hom(B,-)$ is right adjoint to $- \otimes B$ for all DG modules $B$. The \emph{linear dual} of $B$ is the DG module $B^{\dual} = \Hom(B,R)$. We denote by $\Sigma_{n}$ the symmetric group of permutations of $[n]=\{ 1 , \ldots, n \}$.

The \emph{suspension} $sX$ of a DG module $X$ is defined by $(sX)_{n} \cong X_{n-1}$, $\partial(sx) = -s\partial(x)$.  For the sake of the Koszul convention, we treat $s$ as a symbol of degree $1$ and as a natural isomorphism of degree $1$. Thus $s^{n} = s \circ \cdots \circ s$ ($n$ times) has degree $n$.  The inverse $s^{-1}$ is called the \emph{desuspension}.  If $f : X \rightarrow Y$ is a linear map of degree $k$, then $sf : sX \rightarrow sY$ is the linear map of degree $k$ defined by $(sf)(x) = (-1)^{k} f(sx)$.

If $X$ and $Y$ are DG modules, then an isomorphism $s^{2}(X \otimes Y) \xrightarrow{\cong} sX \otimes sY$ is defined by $s^{2}(x \otimes y) \mapsto (-1)^{\deg x} sx \otimes sy$.

A \emph{symmetric sequence of DG modules} is a sequence $\op{X} = ( \op{X}(n) )_{n \geq 0}$, where $\op{X}(n)$ is a
right module over the symmetric group $\Sigma_{n}$, for all $n \geq 0$. The parameter $n$ is referred to as the \emph{arity}. A morphism of symmetric sequences $\varphi:\op{X} \rightarrow \op{Y}$ is a sequence of morphisms, $(\varphi_{n} : \op{X}(n) \rightarrow \op{Y}(n) )_{n \geq 0}$, where each $\varphi_{n}$ is $\Sigma_{n}$-equivariant.  The category of symmetric sequences and their morphisms is denoted $\cat{dgM}^{\Sigma}$.  We say that $\op{X}$ is \emph{connected} if $\op{X}(0) = 0$, and \emph{projective} if each $\op{X}(n)$ is $R[\Sigma_{n}]$-projective.

When convenient, we may index our symmetric sequences by finite sets in the usual way~\cite{fresse:00}.  Let $I$ be a finite set of cardinality $n$.  If $\op{X}$ is a symmetric sequence, we set $\op{X}(I) = \left(\bigoplus_{\alpha} \op{X}(n)\right)/\sim$, where the direct sum is over all bijections $\alpha:[n] \rightarrow I$, and $(x\cdot\sigma)_{\alpha} \sim x_{\alpha\sigma^{-1}}$.  Note that this is the same as the colimit of the diagram with one map $X_{\beta} \rightarrow X_{\beta\sigma}$, $x_{\beta} \mapsto (x\sigma)_{\beta\sigma}$, for each $\beta:[n]\xrightarrow{\cong} I$ and for each $\sigma \in \Sigma_{n}$.

The \emph{graded tensor product} of symmetric sequences $\op{X}$ and $\op{Y}$ is the symmetric sequence $\op{X} \odot \op{Y}$ defined by
\[
    \op{X} \odot \op{Y} (I) = \bigoplus_{I_{1} \amalg I_{2} = I} \op{X}(I_{1}) \otimes \op{Y}(I_{2}).
\]
%
%
%
%
The \emph{composition product} of symmetric sequences is then given by
\[
    \op{X} \circ \op{Y} = \bigoplus_{m \geq 0} \op{X}(m)
    \underset{\Sigma_{m}}{\otimes} \op{Y}^{\odot m}.
\]
The composition product is associative~\cite{joyal:81,smirnov:81}.  It has a unit: the
symmetric sequence $\op{J}$ defined by
\[
    \op{J}(n) = \left\{
        \begin{array}{cc}
            R   & \text{if }n=1, \\
            0   & \text{otherwise.}
        \end{array}
        \right.
\]
By~\cite{rezk:96} for example, $(\cat{dgM}^{\Sigma},\circ,\op{J})$
is a right-closed monoidal category.

\subsection{Operads}\label{ssec:operad-basics}

An \emph{operad} is a monoid in $(\cat{dgM}^{\Sigma}, \circ,
\op{J})$.  Thus, an operad is a symmetric sequence $\op{P}$ equipped
with an associative multiplication $\gamma : \op{P} \circ \op{P} \rightarrow
\op{P}$ that is unital with respect to the unit $\eta : \op{J} \rightarrow \op{P}$.  A morphism of operads is then a morphism of monoids.  We denote by $\cat{Op}$ the category of operads and operad morphisms.

An \emph{augmented operad} is an operad $\op{P}$ along with a morphism of operads, $\epsilon : \op{P} \rightarrow \op{J}$.  The \emph{augmentation ideal} of $\op{P}$ is $\tilde{\op{P}} = \ker\epsilon$.


\begin{ex}\label{ex:assoc}
The \emph{associative operad} $\op{A}$ is defined as the symmetric
sequence $\op{A}(n) = R[\Sigma_{n}]$ with $\Sigma_{n}$ acting on
the right by multiplication.  The composition product is defined by
feeding permutations within blocks into a block permutation.

The \emph{commutative operad} $\op{C}$ is the symmetric sequence
$\op{C}(n) = R$ of trivial representations of $\Sigma_{n}$.  The
composition product is defined by multiplication.
\end{ex}

\begin{ex}\label{ex:susp-op}
The \emph{suspension operad} $\op{S}$ plays an important role in the
theory of quadratic operads.  Let $\op{S}(n)$ be the free graded
$R$-module generated by an element $s_{n-1}$ of degree $n-1$,
equipped with the sign representation of $\Sigma_{n}$. Suppose
$\vec{m} = (m_{1}, \ldots, m_{n}) \in I_{n,m}$. There exists $\sigma
\in \Sigma_{n}$ such that $\sigma\vec{m} = ( m_{\sigma^{-1}(1)},
\ldots, m_{\sigma^{-1}(n)} )$ is non-increasing. Denote by
$\sigma_{\vec{m}} \in \Sigma_{m}$ the block permutation determined
by $\sigma$ and $\vec{m}$.  Let $\kappa(\sigma,\vec{m})$ be the sign
introduced by the Koszul rule in the map,
\[
    s_{m_{1}-1} \otimes \cdots \otimes s_{m_{n}-1}
    \mapsto \kappa(\sigma,\vec{m})
    s_{m_{\sigma^{-1}(1)}-1} \otimes \cdots \otimes s_{m_{\sigma^{-1}(n)}-1}.
\]
Clearly, $\kappa(\rho\sigma, \vec{m}) = \kappa(\rho,\sigma\vec{m})
\kappa(\sigma,\vec{m})$.  In particular, $\kappa(\sigma,\vec{m}) =
\kappa(\sigma^{-1}, \sigma\vec{m})$. Define
\begin{equation}\label{eq:susp}
    \gamma_{\vec{m}}(s_{n-1} \otimes s_{m_{1} - 1} \otimes \cdots
    \otimes s_{m_{n}-1}) = \alpha(\sigma,\vec{m})
    s_{m-1},
\end{equation}
where
\[
    \alpha(\sigma,\vec{m}) = (\sgn \sigma)(\sgn
    \sigma_{\vec{m}})\kappa(\sigma,\vec{m}).
\]
We need to verify that $\alpha(\sigma,\vec{m})$ is independent of
$\sigma$. If $\rho \in \Sigma_{n}$ is another permutation such that
$\rho \vec{m}$ is non-increasing, then $\rho\vec{m}$ and
$\sigma\vec{m}$ are equal as lists.  Thus $\rho = \pi\sigma$, where
$\pi$ fixes $\sigma\vec{m}$.  Since $\sgn$ is a group homomorphism,
and since $(\pi\sigma)_{\vec{m}} =
\pi_{\sigma\vec{m}}\sigma_{\vec{m}}$, it follows that
\[
    \alpha(\rho,\vec{m}) =
    \alpha(\pi,\sigma\vec{m})\alpha(\sigma,\vec{m}).
\]
Thus it suffices to show that $\alpha(\pi,\sigma\vec{m}) = 1$. Let
$k_{i} = m_{\sigma^{-1}(i)}$, for $i = 1, \ldots, n$. Since
$\sigma\vec{m}$ is non-increasing, $\pi$ is the composition of
transpositions $\tau_{j} = (j,j+1)$, where $k_{j} = k_{j+1}$.  Thus
$\sgn \tau_{j} = -1$ and $\sgn (\tau_{j})_{\vec{k}} = (-1)^{k_{j}}$.
Now, $\kappa(\tau_{j},\vec{k})$ is the sign introduced by
interchanging two copies of $s_{k_{j}-1}$, so
$\kappa(\tau_{j},\vec{k}) = (-1)^{k_{j}-1}$.  Thus
\[
    \alpha(\tau_{j},\vec{k}) = (-1)(-1)^{k_{j}}(-1)^{k_{j}-1} = 1.
\]
It follows that $\alpha(\pi,\sigma\vec{m}) =1$, and so
\eqref{eq:susp} is independent of the choice of $\sigma$.
Furthermore, the $\gamma_{\vec{m}}$ have the equivariance properties
to define a product $\gamma : \op{S} \circ \op{S} \rightarrow
\op{S}$, making $\op{S}$ into an operad.  The unit is defined by the
isomorphism $\eta : \op{J}(1) \cong \op{S}(1)$.
\end{ex}

\subsubsection{Algebras}\label{sssec:alg}

The \emph{Schur functor} $T : \cat{dgM}^{\Sigma} \rightarrow \End(\cat{dgM})$ is
defined on objects by
\[
    T_{\op{X}}(V) = \bigoplus_{n \geq 0}
                        \op{X}(n) \underset{\Sigma_{n}}{\otimes} V^{\otimes
                        n}
\]
for $\op{X} \in \cat{dgM}^{\Sigma}$.  A morphism of symmetric
sequences $\varphi : \op{X} \rightarrow \op{Y}$ induces a natural
transformation of Schur functors, $T_{\op{X}} \Rightarrow T_{\op{Y}}$.

The functor $T$ is monoidal : $T_{\op{X}\circ \op{Y}} \cong
T_{\op{X}} \circ T_{\op{Y}}$. If $\op{P}$ is an operad, then
$T_{\op{P}}$ is the underlying functor of a monad
$\monad{T}_{\op{P}}$, with unit defined by the obvious inclusion on
the $n=1$ summand, and multiplication defined by the product
$\op{P}\circ \op{P} \rightarrow \op{P}$.  A $\op{P}$-\emph{algebra}
is an algebra for the monad $\monad{T}_{\op{P}}$. Thus a
$\op{P}$-algebra is a differential graded module $A$, equipped with
structure morphisms for all $n \geq 0$,
\[
    \lambda_{n} : \op{P}(n) \otimes A^{\otimes n} \rightarrow A
\]
that are equivariant, associative and unital.  We will use the
notation
\[
    p(a_{1}, \ldots, a_{n}) := \lambda_{n} ( p \otimes a_{1} \otimes
    \cdots \otimes a_{n})
\]
for all $p \in \op{P}(n)$, $a_{1}, \ldots, a_{n} \in A$.

A morphism of $\op{P}$-algebras is a chain map that preserves the
structure morphisms.  The category of $\op{P}$-algebras and
morphisms is denoted $\alg{P}$.

Let $V \in \cat{dgM}$.  The \emph{free $\op{P}$-algebra on $V$} is the DG module $T_{\op{P}}(V)$, with structure morphism and unit coming from those of the monad $\monad{T}_{\op{P}}$.

\subsubsection{Derivations}\label{sec:derivation}

Let $A$ be a $\op{P}$-algebra.  Let $g : A \rightarrow A$ be a
linear map of degree $k$.  We say that $g$ is a
$\op{P}$-\emph{derivation} if
\[
    g(p(a_{1}, \ldots, a_{n})) = \sum_{i = 1}^{n} (-1)^{k
    \ell_{i}} p(a_{1}, \ldots, g(a_{i}), \ldots, a_{n}),
\]
(where $\ell_{i} = \deg a_{1} + \cdots + \deg a_{i-1}$) for all $n
\geq 0$,  $p \in \op{P}(n)$, and $a_{1}, \ldots , a_{n} \in A$.

\subsubsection{Normalization}\label{sssec:norm}

(See~\cite[1.1.12]{fresse:00}.) We may suppose that our operad $\op{P}$ satisfies $\op{P}(0)=0$, by taking the augmentation ideals of our algebras.  More fully, we note that $P = \op{P}(0)$ is the initial $\op{P}$-algebra. Indeed, in arity zero, the composition product is
the sum of morphisms
\[
    \gamma : \op{P}(n) \otimes \op{P}(0)^{\otimes n} \rightarrow
    \op{P}(0)
\]
that are precisely the structure maps of a $\op{P}$-algebra. Let $A$
be a $\op{P}$-algebra.  Then the structure morphism $\op{P}(0)
\otimes A^{\otimes 0} \rightarrow A$ provides a map, $\eta_{A} : P
\rightarrow A$.  Since the structure morphism in $A$ is associative,
$\eta_{A}$ is a morphism of $\op{P}$-algebras.  Since this is the
only possible $\op{P}$-algebra  morphism from $P$ to $A$, it follows
that $P$ is initial.

Consider the category of \emph{augmented $\op{P}$-algebras}: the
objects are $\op{P}$-algebra morphisms $\varepsilon_{A} : A
\rightarrow P$; morphisms are commutative triangles
\[
    \xymatrix{
        A \ar[rr]^{\varphi} \ar[dr]_{\varepsilon_{A}}
            & & B \ar[dl]^{\varepsilon_{B}} \\
        & P.
    }
\]
We set $IA = \ker\varepsilon_{A}$;  this is the \emph{augmentation
ideal} of $A$.  As in the classical case, $\varepsilon_{A} \circ
\eta_{A} = 1_{A}$, so $A \cong P \oplus IA$. Let $\tilde{\op{P}}(n) = \op{P}(n)$ if $n > 0$, while
$\tilde{\op{P}}(0) = 0$.  Then $\tilde{\op{P}}$ is a sub-operad of
$\op{P}$, and $IA$ is a $\tilde{\op{P}}$-algebra.

\subsubsection{Algebras as left modules}\label{sssec:alg-leftmod}

Let $z : \cat{dgM} \rightarrow \cat{dgM}^{\Sigma}$ be the functor
defined by
\[
    z(V)(n) = \left\{
        \begin{array}{cc}
            V   & \text{if }n=0 \\
            0   & \text{otherwise.}
        \end{array}
        \right.
\]
Then $z$ restricts to define a functor $\alg{P} \rightarrow
\Mod{P}{}$.  Indeed, it is an easy exercise to show that
$z(T_{\op{P}}(A)) \cong \op{P} \circ z(A)$, so applying $z$ to the
structure morphism $ T_{\op{P}}(A) \rightarrow A$, we obtain the
structure morphism
\[
    \op{A} \circ z(A) \cong z(T_{\op{P}}(A)) \rightarrow z(A),
\]
so $z(A)$ is a left $\op{P}$-module.  It is immediate from the definition that $z$ is full and faithful.

Since $\alg{P}$ embeds in ${}_{\op{P}}\cat{Mod}$, it is often easier to state and prove results about left modules rather than algebras.

\subsection{Cooperads}

A \emph{cooperad} is a symmetric sequence $\op{Q}$ along with a counit $\varepsilon : \op{Q}(1) \rightarrow R$ and comultiplications
\begin{equation}\label{eqn:coop-struct}
    \psi_{\vec{n}} : \op{Q}(n) \rightarrow \op{Q}(m) \otimes
    \op{Q}(n_{1}) \otimes \cdots \otimes \op{Q}(n_{m})
\end{equation}
for all $\vec{n} = (n_{1}, \ldots, n_{m}) \in I_{m,n}$, that are coassociative, counital, and equivariant.

If $\op{Q}$ is connected, then the above morphisms yield a sequence of $\Sigma_{n}$-equivariant morphisms,
\[
    \tilde{\psi}_{n} : \op{Q}(n) \rightarrow \bigoplus_{m = 1}^{n} \left(
        \op{Q}(m) \otimes \op{Q}[m,n] \right)^{\Sigma_{m}}
\]
that we may compose with the natural map from fixed points to orbits, to obtain a morphism of symmetric sequences, $\psi : \op{Q} \rightarrow \op{Q} \circ \op{Q}$.  In fact, by~\cite[Proposition 1.1.15]{fresse:00}, the $\psi_{\vec{n}}$ may be recovered from $\psi$ and the two notions are equivalent. \emph{Henceforth, we assume that our cooperads are connected, and so a cooperad is simply a comonoid with respect to the composition product.}

A morphism of cooperads is a morphism of symmetric sequences, $\varphi : \op{Q} \rightarrow \op{K}$, that commutes with the
structure morphisms \eqref{eqn:coop-struct}.  We denote by $\cat{CoOp}$ the category of cooperads and cooperad morphisms.

\begin{ex}
The structure morphisms in the suspension operad $\op{S}$ are
isomorphisms, and so $\op{S}$ is also a cooperad.
\end{ex}

\subsubsection{Coalgebras}\label{sssec:coalg}

Let $\op{X} \in \cat{dgM}^{\Sigma}$.  Let $\cat{dgM}_{+}$ be the full subcategory of $\cat{dgM}$ consisting of all chain complexes concentrated in strictly positive degrees. Let $V \in \cat{dgM}_{+}$. We define
\[
    \Gamma_{\op{X}}(V) = \bigoplus_{n \geq 1} (\op{X}(n) \otimes
    V^{\otimes n})^{\Sigma_{n}}.
\]
By~\cite[Propositions 1.1.9 and 1.1.15]{fresse:00},
\[
    \Gamma : \cat{dgM}^{\Sigma} \rightarrow \End(\cat{dgM}_{+}), \quad \op X \mapsto \Gamma_{\op X}
\]
is monoidal when restricted to connected projective symmetric sequences of finite type.  Thus, if $\op Q$ is a projective cooperad of finite type, then the composition diagonal in $\op{Q}$ makes $\Gamma_{\op{Q}}$ the underlying functor of a comonad, $\monad{\Gamma}_{\op{Q}}$.  A $\op{Q}$-\emph{coalgebra} is a coalgebra over $\monad{\Gamma}_{\op{Q}}$.  Thus, a $\op{Q}$-coalgebra is a DG module $C$ along with structure morphisms $\rho_{n}:C \rightarrow \op{Q} \otimes C^{\otimes n}$ for all $n$, that are appropriately counital, coassociative, and equivariant.

A $\op{Q}$-coalgebra $C$ is \emph{co-nilpotent} if for all $c \in C$, there exists $N \geq 1$
such that $\rho_{n}(c) = 0$ whenever $n \geq N$.

The functor $\Gamma_{\op{Q}}$ is the right adjoint to the forgetful functor from co-nilpotent {$\op{Q}$}-coalgebras to $\cat{dgM}_{+}$. That is, if $C$ is a co-nilpotent $\op{Q}$-coalgebra and $f : C \rightarrow V$ is a chain map, then $f$ lifts uniquely through the projection $\Gamma_{\op{Q}}(V) \rightarrow V$ to define a morphism of coalgebras, $\phi : A \rightarrow \Gamma_{\op{P}}(V)$.  To obtain $\phi$, one simply sums the composites
\[
    C \rightarrow \op{Q} \otimes C^{\otimes n} \rightarrow \op{Q} \otimes V^{\otimes n}.
\]
Therefore, $\Gamma_{\op{Q}}(V)$ is the \emph{co-free co-nilpotent $\op{Q}$-coalgebra co-generated by $V$}.

\subsubsection{Coderivations}

Let $C$ be a $\op{Q}$-coalgebra.  A map of degree $k$, $g : C \rightarrow C$, is called a
$\op{Q}$-\emph{coderivation} if the following diagram commutes for all $n \geq 0$.
\[
    \xymatrix{
        C  \ar[r]^{\rho} \ar[d]_{g}
            & \op{Q}(n) \otimes C^{\otimes n} \ar[d]^{1 \otimes \sum_{j=0}^{n-1} 1^{\otimes j} \otimes g
            \otimes 1^{\otimes n-j-1}} \\
        C \ar[r]^{\rho}
            & \op{Q}(n) \otimes C^{\otimes n}
    }
\]
Let $V \in \cat{dgM}$.  A map $t : \Gamma_{\op{Q}}(V) \rightarrow V$ of degree $k$ lifts uniquely
to determine a $\op{Q}$-coderivation $g$ on $\Gamma_{\op{Q}}(V)$.

\subsubsection{Coalgebras as left comodules}

Let $\op{Q}$ be a projective cooperad of finite type. Recall the functor $z : \cat{dgM} \rightarrow \cat{dgM}^{\Sigma}$ of Section~\ref{sssec:alg-leftmod}.  Since $z(\Gamma_{\op{Q}}(V)) = \op{Q} \circ z(V)$, $z$ restricts to define a functor, $z : \coalg{Q} \rightarrow \Comod{Q}{}$.

\subsection{A twisting structure for $\cat{dgProj}^\Sigma$}\label{subsec:twist-chcx}

We now show that a reasonable subcategory of symmetric sequences of chain complexes is twistable, and then use the operadic bar construction with coefficients to construct an explicit twisting structure.

\subsubsection{Twistability of $\cat{dgProj}^\Sigma$}\label{sssec:twistability}

In this section, we prove that a certain subcategory of symmetric sequences is twistable.  Let $\cat{Proj}$ be the full subcategory of $\cat{M}$ spanned by finitely generated projective $R$-modules.  Then $\cat{gProj}$ and $\cat{dgProj}$ are the categories of graded and differential graded projective $R$-modules of finite type, respectively.  Let $\cat{dgProj}^{\Sigma}$ denote the category of symmetric sequences in $\cat{dgProj}$.

Let $\op{P}$ be an operad.  Let $\op{X}$ and $\op{Y}$ be left and right $\op{P}$-modules, respectively. Recall that composition over $\op{P}$ is defined by a reflexive coequalizer,
\[
    \xymatrix{
        \op{X} \circ \op{P} \circ \op{Y}
        \ar@<1ex>[r] \ar@<-1ex>[r] &
        \op{X} \circ \op{Y} \ar[r] & \op{X} \underset{\op{P}}{\circ} \op{Y},
    }
\]
where the mutual section is provided by the unit morphism $\op{J} \rightarrow \op{P}$.  Furthermore, $-\otimes-$ commutes with reflexive coequalizers simultaneously in both factors.  Since the composition product is built from colimits, we obtain the following result.

\begin{prop}\label{prop:circ-over-assoc}\cite[2.3.12,2.3.13]{rezk:96}
Let $\op{O}$, $\op{P}$, and $\op{Q}$ be operads.
\begin{enumerate}
    \item\cite[2.3.12]{rezk:96} Let $\op{X}$ be an $(\op{O},\op{P})$-bimodule, and let $\op{Y}$ be a $(\op{P},\op{Q})$-bimodule.  Then $\op{X} \underset{\op{P}}{\circ} \op{Y}$ has a natural $(\op{O},\op{Q})$-bimodule structure.
    \item\cite[2.3.13]{rezk:96} Let $\op{X}$ be a right $\op{O}$-module, let $\op{Y}$ be an $(\op{O},\op{P})$-bimodule, and let $\op{Z}$ be a right $\op{P}$-module. Then there is a unique isomorphism
        \[
            (\op{X} \underset{\op{O}}{\circ} \op{Y}) \underset{\op{P}}{\circ} \op{Z}
            \cong
            \op{X} \underset{\op{O}}{\circ} (\op{Y} \underset{\op{P}}{\circ} \op{Z})
        \]
        commuting with the natural maps $(\op{X} \circ \op{Y}) \circ \op{Z} \rightarrow (\op{X} \underset{\op{O}}{\circ} \op{Y}) \underset{\op{P}}{\circ} \op{Z}$ and $\op{X} \circ ( \op{Y} \circ \op{Z}) \rightarrow \op{X} \underset{\op{O}}{\circ} (\op{Y} \underset{\op{P}}{\circ} \op{Z})$.
\end{enumerate}
\end{prop}

Dually, a \emph{coreflexive equalizer} is a diagram,
\[
    \xymatrix{
        X \ar[r]^{k} & Y \ar@<1ex>[r]^{f} \ar@<-1ex>[r]_{g}
        & Z, & Z \ar[r]^{r} & Y,
    }
\]
in which $rf = rg = 1_{Y}$ and $X \xrightarrow{k} Y$ is final among morphisms equalizing $f$ and $g$.  As in the dual case, if $F : \cat{C} \times \cat{C} \rightarrow \cat{C}$ is a bifunctor such that $F(X,-)$ and $F(-,X)$ preserve coreflexive equalizers, then $F$ preserves coreflexive equalizers simultaneously in each variable~\cite[Corollary 1.2.12]{johnstone:02}.

\begin{prop}\label{prop:circ-preserve}
Let $\op{X} \in \cat{dgProj}^{\Sigma}$ be connected (that is, $\op{X}(0)=0$). Then $\op{X} \circ -$ and $- \circ \op{X}$ preserve coreflexive equalizers.
\end{prop}

\begin{proof}
At the level of DG modules, if $M$ is projective, then $M \otimes -$ and $- \otimes M$ preserve coreflexive equalizers of DG modules, since these can be defined as kernels.  Therefore $- \otimes -$ preserves coreflexive equalizers in $\cat {dgProj}$. It follows that if
\[
    \xymatrix{
        X_{i} \ar[r]^{k_{i}} & Y_{i} \ar@<1ex>[r]^{f_{i}} \ar@<-1ex>[r]_{g_{i}}
        & Z_{i}, & Z_{i} \ar[r]^{r_{i}} & Y_{i},
    }
\]
is a coreflexive equalizer for $i = 1, \ldots, m$ in $\cat {dgProj}$, then
\[
    \xymatrix{
        \bigotimes_{i=1}^{m}X_{i} \ar[r]^{\otimes k_{i}}
            & \bigotimes_{i=1}^{m} Y_{i} \ar@<1ex>[r]^{\otimes f_{i}} \ar@<-1ex>[r]_{\otimes g_{i}}
            & \bigotimes_{i=1}^{m} Z_{i},
            & \bigotimes_{i=1}^{m} Z_{i} \ar[r]^{\otimes r_{i}} & \bigotimes_{i=1}^{m} Y_{i},
    }
\]
is a coreflexive equalizer.

Let $\vec{n} = (n_{1},\ldots,n_{m})$ be an $m$-partition of $n$ consisting of positive integers. Let $\Sigma_{\vec{n}} = \Sigma_{n_{1}} \times \cdots \times \Sigma_{n_{m}}$.  In general, if $G$ is a discrete group and $H < G$ is a subgroup, then $R[G]$ is free as an $R[H]$-module, with basis $G/H$.  Correspondingly, write $R[\Sigma_{n}] = R[\Sigma_{\vec{n}}] \otimes V$, where $V$ is the free $R$-module on basis $\Sigma_{n}/\Sigma_{\vec n}$.  Then if $\op{Y} \in \cat{Proj}^{\Sigma}$,
\[
    \left(\op{Y}(n_{1}) \otimes \cdots \otimes  \op{Y}(n_{m})\right)
    \otimes_{\Sigma_{\vec{n}}} R[\Sigma_{n}]
    \cong
    \op{Y}(n_{1}) \otimes \cdots \otimes  \op{Y}(n_{m}) \otimes V.
\]
Since finite direct sums of chain complexes are naturally isomorphic to finite direct products, it follows that
\begin{equation}\label{eq:odot-ce}
    \xymatrix{
        \op{X}^{\odot m}(n) \ar[r]^{k^{\odot m}}
            & \op{Y}^{\odot m}(n)
                \ar@<1ex>[r]^{f^{\odot m}}
                \ar@<-1ex>[r]_{g^{\odot m}}
            & \op{Z}^{\odot m}(n),
            & \op{Z}^{\odot m}(n)
                \ar[r]^{r^{\odot m}}
            & \op{Y}^{\odot m}(n)
    }
\end{equation}
is a coreflexive equalizer.

Now, by~\cite[Lemma 1.1.16]{fresse:00},  for all $m,n$, there exists $\op X_{m}(n) \in \cat{dgProj}$ and  a natural isomorphism $\op{X}^{\odot m}(n) \cong R[\Sigma_{m}] \otimes \op{X}_{m}(n)$, since $\op{X}$ is connected.  It follows that $h^{\odot m}:\op{X}^{\odot m}(n) \rightarrow \op{Y}^{\odot m}(n)$ is of the form $R[\Sigma_{m}] \otimes h_{m}$, where $h_{m}:\op{X}_{m}(n) \rightarrow \op{Y}_{m}(n)$, for any morphism $h:\op{X} \rightarrow \op{Y}$.  Therefore, if \eqref{eq:odot-ce} is a coreflexive equalizer, so too is
\[
    \xymatrix{
        \op{X}_{m}(n) \ar[r]^{k_{m}}
            & \op{Y}_{m}(n)
                \ar@<1ex>[r]^{f_{m}}
                \ar@<-1ex>[r]_{g_{m}}
            & \op{Z}_{m}(n),
            & \op{Z}_{m}(n) \ar[r]^{r_{m}}
            & \op{Y}_{m}(n).
        }
\]
Let $\op{W} \in \cat{Proj}^{\Sigma}$.  Then for all $m, n$ we have a natural isomorphism $\op{W}(m) \otimes_{\Sigma_{m}} \op{X}^{\odot m}(n) \cong \op{W}(m) \otimes \op{X}_{m}(n)$, and similarly for $\op{Y}$ and $\op{Z}$.  It follows that
\[
    \xymatrix{
        \op{W} \otimes_{\Sigma_{m}} \op{X}^{\odot m}(n) \ar[r]^{k^{\odot m}}
            & \op{W} \otimes_{\Sigma_{m}} \op{Y}^{\odot m}(n)
                \ar@<1ex>[r]^{f^{\odot m}}
                \ar@<-1ex>[r]_{g^{\odot m}}
            & \op{W} \otimes_{\Sigma_{m}} \op{Z}^{\odot m}(n), \\
            & \op{W} \otimes_{\Sigma_{m}} \op{Z}^{\odot m}(n)
                \ar[r]^{r^{\odot m}}
            & \op{W} \otimes_{\Sigma_{m}} \op{Y}^{\odot m}(n)
    }
\]
is a coreflexive equalizer.

Finally, we need to show that
\begin{equation}\label{eq:sum-ce}
    \xymatrix{
        \op{W} \circ \op{X}
            \ar[r]^{1 \circ k}
            & \op{W} \circ \op{Y}
                \ar@<1ex>[r]^{1 \circ f}
                \ar@<-1ex>[r]_{1 \circ g}
            & \op{W} \circ \op{Z},
            & \op{W} \circ \op{Z}
                \ar[r]^{1 \circ r}
            & \op{W} \circ \op{Y}
    }
\end{equation}
is a coreflexive equalizer.  Let $\ell : \op{U} \rightarrow \bigoplus_{m \geq 1} \op{W} \otimes_{\Sigma_{m}} \op{Y}^{\odot m}$ equalize $\sum f^{\odot m}$ and $\sum g^{\odot m}$. Then $\ell(u) = \sum \ell_{m}(u)$, with $\ell_{m}(u) \in \op{W} \otimes_{\Sigma_{m}} \op{Y}^{\odot m}$, for all $ u \in \op{U}$.  It is immediate that $\ell_{m}$ equalizes $f^{\odot m}$ and $g^{\odot m}$.  Therefore we get morphisms $j_{m} : \op{U} \rightarrow \op{W}(m) \otimes_{\Sigma_{m}} \op{X}^{\odot m}$ that satisfy $k^{\odot m}j_{m} = \ell_{m}$.  Since for each $w \in \op{W}$, only finitely many $\ell_{m}(w)$ are nonzero, we may add the $j_{m}$'s to obtain a morphism $j : \op{U} \rightarrow \op{W} \circ \op{X}$ such that $kj=\ell$.  Therefore \eqref{eq:sum-ce} is a coreflexive equalizer, and so $\op{W} \circ -$ preserves coreflexive equalizers.

The proof that $- \circ \op{Y}$ preserves coreflexive equalizers follows the same lines, but is easier.
\end{proof}

Let $\op{Q}$ be a cooperad; let $\op{X}$ be a right $\op{Q}$-comodule and let $\op{Y}$ be a left $\op{Q}$-comodule. The equalizer that defines $\op{X} \underset{\op{Q}}{\square} \op{Y}$ is coreflexive; the common left inverse $\op{X} \circ \op{Q} \circ \op{Y} \rightarrow \op{X} \circ \op{Y}$ is provided by the counit, $\op{Q} \rightarrow \op{J}$.  From Proposition~\ref{prop:circ-preserve}, it follows that $(\op{X} \underset{\op{Q}}{\square} \op{Y}) \circ \op{Z} \cong \op{X} \underset{\op{Q}}{\square} ( \op{Y} \circ \op{Z} )$ for any $\op{Z} \in \cat{Proj}^{\Sigma}$, and likewise on the other side.  We therefore obtain the following results.

\begin{prop}\label{prop:square-bicomod}
Let $\op{O}$, $\op{P}$ and $\op{Q}$ be cooperads in $\cat{dgProj}^\Sigma$.  Let $\op{X}$ be an $(\op{O},\op{P})$-bicomodule and $\op{Y}$  a $(\op{P},\op{Q})$-bicomodule, both in $\cat{dgProj}^\Sigma$.  Then $\op{X} \underset{\op{P}}{\square} \op{Y}$ has a natural $(\op{O},\op{Q})$-bicomodule structure.
\end{prop}

\begin{prop}\label{prop:square-assoc}
Let $\op{P}$ and $\op{Q}$ be cooperads in $\cat{dgProj}^\Sigma$.  Let $\op{X}$ be a right $\op{P}$-comodule, let $\op{Y}$ be a $(\op{P},\op{Q})$-bicomodule, and let $\op{Z}$ be a left $\op{Q}$-comodule, all in $\cat{dgProj}^\Sigma$.  Then there is a unique isomorphism
\[
    (\op{X} \underset{\op{P}}{\square} \op{Y}) \underset{\op{Q}}{\square} \op{Z}
    \cong \op{X} \underset{\op{P}}{\square} (\op{Y} \underset{\op{Q}}{\square} \op{Z})
\]
that commutes with the natural maps $(\op{X} \underset{\op{P}}{\square} \op{Y}) \underset{\op{Q}}{\square} \op{Z} \rightarrow (\op{X} \circ \op{Y}) \circ \op{Z}$ and $\op{X} \underset{\op{P}}{\square} (\op{Y} \underset{\op{Q}}{\square} \op{Z}) \rightarrow \op{X} \circ ( \op{Y} \circ \op{Z})$.
\end{prop}

Our next proposition of the section shows that relative (co)composition products behave well with respect to bi(co)module structures.

\begin{prop}\label{prop:mix}
Let $\op{O}$ and $\op{Q}$ be either operads or cooperads, or one of each, both in $\cat{dgProj}^\Sigma$.
\begin{enumerate}
    \item Let $\op{P}$ be an operad. Then
    \[
        - \underset{\op{P}}{\circ} - : {}_{\op{O}}\cat{Mix}_{\op{P}} \times {}_{\op{P}}\cat{Mix}_{\op{Q}} \rightarrow {}_{\op{O}}\cat{Mix}_{\op{Q}}.
    \]
    \item Let $\op{P}$ be a cooperad. Then
    \[
        - \underset{\op{P}}{\square} - :{}_{\op{O}}\cat{Mix}_{\op{P}} \times {}_{\op{P}}\cat{Mix}_{\op{Q}} \rightarrow {}_{\op{O}}\cat{Mix}_{\op{Q}}.
    \]
\end{enumerate}
\end{prop}

\begin{proof}
The proof consists of a sequence of straightforward verifications, using natural isomorphisms
\[
    \op{X} \circ (\op{Y} \underset{\op{P}}{\circ} \op{Z}) \cong (\op{X} \circ \op{Y}) \underset{\op{P}}{\circ} \op{Z}
\]
\[
    \op{X} \underset{\op{P}}{\circ} (\op{Y} \circ \op{Z}) \cong (\op{X} \underset{\op{P}}{\circ} \op{Y}) \circ \op{Z}
\]
\[
    \op{X} \circ (\op{Y} \underset{\op{P}}{\square} \op{Z}) \cong (\op{X} \circ \op{Y}) \underset{\op{P}}{\square} \op{Z}
\]
\[
    \op{X} \underset{\op{P}}{\square} (\op{Y} \circ \op{Z}) \cong (\op{X} \underset{\op{P}}{\square} \op{Y}) \circ\op{Z}.
\]
\end{proof}

Our final proposition shows that relative cocomposition behaves well with respect to relative composition.

\begin{prop}\label{prop:cocomp-assoc}
Let $\op{P}$ be an operad and let $\op{Q}$ be a cooperad.
\begin{enumerate}
    \item Let $\op{X}$ be a right $\op{P}$-module, let $\op{Y}$ be a $(\op{P},\op{Q})$-mixed module, and let $\op{Z}$ be a left $\op{Q}$-comodule.  Then there is a unique isomorphism
        \[
            \left( \op{X} \underset{\op{P}}{\circ} \op{Y} \right) \underset{\op{Q}}{\square} \op{Z}
            \cong
            \op{X} \underset{\op{P}}{\circ} \left( \op{Y} \underset{\op{Q}}{\square} \op{Z} \right)
        \]
        compatible with all relevant natural maps.
    \item Let $\op{X}$ be a right $\op{Q}$-comodule, let $\op{Y}$ be a $(\op{Q},\op{P})$-mixed module, and let $\op{Z}$ be a left $\op{P}$-module.  Then there is a unique isomorphism
        \[
            \left( \op{X} \underset{\op{Q}}{\square} \op{Y} \right) \underset{\op{P}}{\circ} \op{Z}
            \cong
            \op{X} \underset{\op{Q}}{\square} \left( \op{Y} \underset{\op{P}}{\circ} \op{Z} \right)
        \]
        compatible with all relevant natural maps.
\end{enumerate}
\end{prop}

\begin{proof}
Let $\cat{C}$ be a category.  Consider the following property.
\begin{equation}\label{property}
\text{Filtered colimits commute with finite limits in }\cat{C}.
\end{equation}
Property~\eqref{property} holds in $\cat{Set}$, by direct verification.  The forgetful functor $\cat{Ab} \rightarrow \cat{Set}$ preserves and reflects filtered colimits and all limits. Let $\cat{M}$ be the category of $R$-modules.  The forgetful functor $\cat{M} \rightarrow \cat{Ab}$ preserves and reflects all limits and colimits.  So~\eqref{property} holds in $\cat{M}$.

Let $\cat{gM}$ be the category of graded $R$-modules.  Let $E_{i}:\cat{gM} \rightarrow \cat{M}$ be defined on objects by $E_{i}(X) = X_{i}$.  Limits and colimits in $\cat{gM}$ are defined degree-wise, so if $F : \cat{D} \rightarrow \cat{gM}$ is a diagram, then
\[
    L = \operatorname{(co)lim}F \quad \Leftrightarrow \quad E_{i}L = \operatorname{(co)lim}E_{i}F \quad \forall i.
\]
It follows that~\eqref{property} holds in $\cat{gM}$.

Let $W : \cat{dgM} \rightarrow \cat{gM}$ be the forgetful functor, and let $F : \cat{D} \rightarrow \cat{dgM}$ be a diagram.  Let $L = \operatorname{(co)lim}WF$.  We consider the differential to be a natural transformation $\partial_{i}:E_{i}W \Rightarrow E_{i-1}W$.  Thus we get a morphism $\partial_{i} E_{i}L \rightarrow E_{i-1}L$ that is a differential by naturality.  It follows that $(L,\partial) = \operatorname{(co)lim}F$, and so $W$ preserves and reflects limits and colimits.  Hence~\eqref{property} holds in $\cat{dgM}$.

Let $U : \cat{dgM}^{G} \rightarrow \cat{dgM}$ be the forgetful functor, where $G$ is a finite group and $\cat{dgM}^{G}$ is the category of right $G$-modules.  Each $g \in G$ determines a natural transformation from $U$ to $U$.  Therefore, if $F:\cat{D} \rightarrow \cat{dgM}^{G}$ is a diagram and $L = \operatorname{(co)lim}UF$, then we can equip $L$ with a right $G$-action compatible with $UF(D) \rightarrow L$ for all $D \in \cat{D}$.  Therefore $L = \operatorname{(co)lim}F$, and so ~\eqref{property} holds in $\cat{dgM}^{G}$.

For $i \geq 1$, let $A_{i}:\cat{dgM}^{\Sigma} \rightarrow \cat{dgM}^{\Sigma_{i}}$ be defined by $A_{i}(\op{X}) = \op{X}(i)$.  Since limits and colimits in $\cat{dgM}^{\Sigma}$ are defined arity-wise, we have that
\[
    L = \operatorname{(co)lim}F \quad \Leftrightarrow \quad A_{i}L = \operatorname{(co)lim}A_{i}F \quad \forall i
\]
whenever $F : \cat{D} \rightarrow \cat{dgM}^{\Sigma}$ is a diagram.  It follows that filtered colimits commute with finite limits in $\cat{dgM}^{\Sigma}$.

We now observe that a coequalizer (such as composition over an operad) is a filtered colimit, and an equalizer (such as cocomposition over a cooperad) is a finite limit, to complete the proof.
\end{proof}

To summarize, we have proved the following result.

\begin{thm}\label{thm:twistable} The category $\cat {dgProj}^\Sigma$ is twistable.
\end{thm}

\subsubsection{Operadic bar construction}\label{sssec:opbar}

We begin by recalling the bar construction of an operad; for details, see~\cite{ginzburg-kapranov:94}.

Let $\op{P}$ be an operad.   Let $I$ and $J$ be finite indexing sets, and let $i \in I$.  We define the $i$th \emph{partial composition},
\[
    \circ_{i}: \op{P}(I) \otimes \op{P}(J) \rightarrow \op{P}((I-\{i\})\amalg J),
\]
by substituting the operation $q \in \op{P}(J)$ in the $i$th entry of the operation $p \in \op{P}(I)$.  Likewise, we set $I \circ_{i} J = (I-\{i\}) \cup J$.

A \emph{tree} is a connected directed graph $\tau$ without loops, in which each vertex has at least one incoming edge and exactly one outgoing edge.  We allow edges without a source and without a target; these are called the \emph{inputs} and \emph{output} of the $\tau$, respectively.  If $\tau$ has an edge that is both an input and an output, then by connectedness, $\tau$ is the degenerate graph $1$, with no vertices and one edge.  If an edge is neither an input nor an output, then it is \emph{internal}.  We denote the set of vertices of $\tau$ by $V(\tau)$ and the internal edges of $\tau$ by $E(\tau)$.  Let $\In(\tau)$ and $\Out(\tau)$ be the set of input edges and the unique output edge of $\tau$, respectively. For $v \in V(\tau)$, let $\In(v)$ be the set of edges entering $v$ and let $\Out(v)$ be the unique output edge of $v$.

For a set $I$, an $I$-\emph{labelled tree} is a tree $\tau$ equipped with a bijection $I \xrightarrow{\cong} \In(\tau)$.  If $I = \{ 1, \ldots, n\}$, then we refer to an $I$-labelled tree simply as an $n$-tree.

An isomorphism of trees is a bijection of vertex sets that preserves the edge relations and the labelling. The category of labelled $n$-trees and isomorphisms is denoted $T(n)$; we note that $T(n) = \amalg_{r \geq 0} T_{(r)}(n)$ where $T_{(r)}(n)$ is the full subcategory of $n$-trees with $r$ vertices.  The trivial tree, $1$, is the degenerate $1$-tree with no vertices.  Note that $\Sigma_{n}$ acts on $T(n)$ by permuting the labels of the entries, and that this action preserves each $T_{(r)}(n)$.  Thus $T = (T(n))$ forms a symmetric sequence in the category of small categories.  Grafting of trees, that is, identifying output edges of $\sigma_{1},\ldots,\sigma_{n}$ with the input edges of the $n$-tree $\tau$ to obtain the tree $\tau(\sigma_{1},\ldots,\sigma_{n})$,  defines an operad structure on $T$.

Let $\op{X}$ be a symmetric sequence.  If $\tau$ is a tree, then we set
\[
    \tau(\op{X}) = \bigotimes_{v \in V(\tau)} \op{X}(In(v)).
\]
We think of an element of $\tau(\op{X})$ as being the tree $\tau$, with each vertex $v$ labeled with an element of $\op{X}(In(v))$.  We set
\[
    FO(\op{X})(n) = \colim_{\tau \in T(n)} \tau(\op{X}).
\]
The \emph{free operad on} $\op{X}$ is the symmetric sequence $FO(\op{X}) = \big(FO(\op{X})(n)\big)$. The grafting of trees determines a composition operation in $FO(\op{X})$.

We can also define a \emph{pruning} operation, dual to grafting, on the symmetric sequence $FO(\op{X})$ to obtain the \emph{free cooperad} $FC(\op{X})$.  Let $\tau$ be a tree.  A \emph{full subtree} $\sigma$ of $\tau$ is a subtree such that $\In_{\sigma}(v) = \In_{\tau}(v)$ for all $v \in V(\sigma)$.  Let $\sigma \subseteq \tau$ be a full subtree containing $\Out(\tau)$; we allow $\sigma=1$ and $\sigma=\tau$.  For each $e_{i} \in \In(\sigma)\cap \In(\tau)$, set $\sigma_{i}=1$.  If $e_{i} \in \In(\sigma)\cap E(\tau)$, then let $\sigma_{i}$ be the largest subtree of $\tau$ such that $e_{i} = \Out(\sigma_{i})$.  Then
\[
    \Delta_{\sigma}(\tau) = \sigma \otimes \sigma_{1} \otimes \cdots \otimes \sigma_{n}.
\]
This pruning operation, summed over all full subtrees $\sigma \subseteq \tau$ containing $\Out(\tau)$, determines the composition diagonal in $FC(\op{X})$.  For example, consider the tree $\tau$:
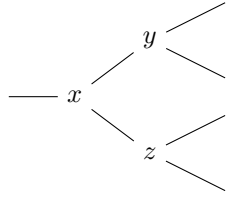
\begin{center}
\begin{tikzpicture}
    \node {} [grow'=right][level distance=10mm,level 3/.style={sibling distance=10mm}]
        child {node {$x$}
            child { node {$y$}
               child child }
            child { node {$z$} child child }
            };
\end{tikzpicture}
\end{center}
labelled with elements of $\op{X}(2)$ at each node.  The tree $\tau$ has three nontrivial subtrees containing $\Out(\tau)$, and so the ``reduced'' composition diagonal of is as pictured below.
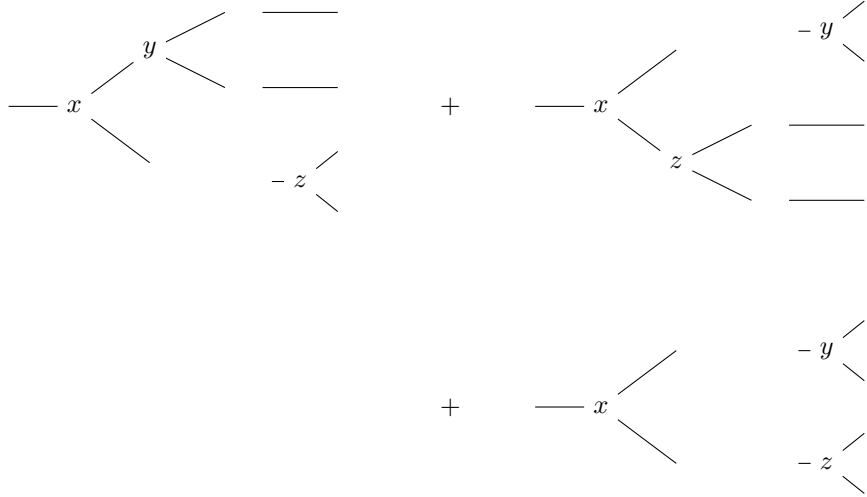
\vspace{.25cm}

\begin{tikzpicture}
    \node at (-7,0) {} [grow'=right,
        level distance=10mm,
        level 3/.style={sibling distance=10mm}]
        child {node {$x$}
            child { node {$y$}
                child child }
                child
                };
    \draw (-3.5,1.25)--(-2.5,1.25);
    \draw (-3.5,0.25)--(-2.5,0.25);
    \node at (-3.5,-1) {} [grow'=right,
        level distance=5mm,
    level 2/.style={sibling distance=8mm}]

        child {node {$z$}
            child child};
    \node at (-1,0) {$+$};
    \node at (0,0) {} [grow'=right,
        level distance=10mm,
        level 3/.style={sibling distance=10mm}]
        child {node {$x$}
            child
            child { node {$z$}
                child child }
            };
    \node at (3.5,1) {} [grow'=right,
        level distance=5mm,
        level 2/.style={sibling distance=8mm}]
        child {node {$y$}
            child child };
    \draw (3.5,-0.25)--(4.5,-0.25);
    \draw (3.5,-1.25)--(4.5,-1.25);
    \node at (-1,-4) {+};
    \node at (0,-4) {} [grow'=right,
        level distance=10mm]
        child { node {$x$}
            child child };
    \node at (3.5,-3.25) {} [grow'=right,
        level distance=5mm,sibling distance=8mm]
            child { node {$y$}
            child child };
    \node at (3.5,-4.75) {} [grow'=right,
        level distance=5mm,sibling distance=8mm]
            child {node {$z$}
            child child };
\end{tikzpicture}

\vspace{0.25in}

Let $\op{P}$ be a connected, augmented operad with augmentation ideal $\tilde{\op{P}}$.  The \emph{bar construction} $\opbar{\op{P}}$, forgetting differentials, is the cofree cooperad $FC(s\tilde{\op{P}})$.  The differential is a perturbation of the internal differential (from the differential $\op{P}$) by a bar differential.  The bar differential of $\tau$ has one term for each internal edge $sp \xleftarrow{e} sq$, where $p,q \in \tilde{\op{P}}$.  We collapse the edge $e$ to a vertex that we label by $-s(p \circ_{e} q)$.

If $\op{M}$ is a right $\op{P}$-module and $\op{N}$ is a left $\op{P}$-module, then the \emph{bar construction with coefficients}, $\opbar(\op{M},\op{P},\op{N})$, is the symmetric sequence $\op{M} \circ \opbar(\op{P}) \circ \op{N}$, with differential perturbed by two terms, $d_{L}$ and $d_{R}$, that come from the actions of $\op{P}$ on $\op{M}$ and $\op{N}$, respectively, that are defined as follows.
We consider an element of $\opbar(\op{M},\op{P},\op{N})$ to be a tree $\tau$ with a partition $V(\tau) = V_L \amalg V_{P} \amalg V_{R}$, where $V_{L}$ consists of the root vertex, and $V_{R}$ consists of all the targets of $\In(\tau)$.  The root vertex $v$ is labeled with an element of $\op{M}(\In(v))$.  The vertices $w \in V_{R}$ are labelled with elements of $\op{N}(\In(w))$. The vertices $u \in V_{P}$ are labelled with elements of $s\tilde{\op{P}}(\In(u))$.

The perturbation $d_{L}$ has one term for each edge $e \in \In(v)$, where $v$ is the unique element of $V_{L}$.  The edge $e$ is collapsed, and the root vertex is relabelled with $x \circ_{u} p$, where $x \in \op{M}$ is the label of $v$ and $sp$ is the label of the source $u$ of $e$.

The perturbation $d_{R}$ has one term for each vertex $w$ that is the target of an edge starting in $V_{R}$.  Let $sp$ be the label of $w$.  Suppose $\In({w}) = \{ e_{1}, \ldots, e_{n}\}$, where  $sp \xleftarrow{e_{i}} y_{i}$ for $y_{i} \in \op{N}$.  Then the contribution to $d_{R}$ for $e_{i}$ comes from collapsing the subtree with root $w$ to a vertex labeled with $-p(y_{1}, \ldots, y_{n})$.

For example, the following picture shows the differential of $x \circ_{1} sp \circ_{1} y$, $x \in \op{M}(2)$, $sp \in s\tilde{\op{P}}(2)$, $y \in \op{N}(2)$.

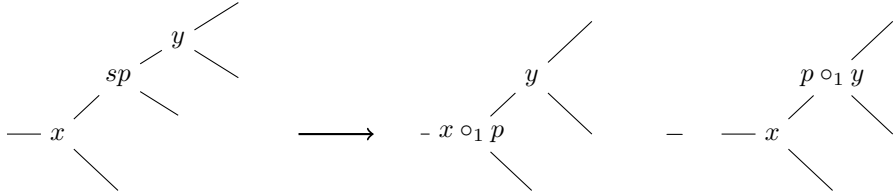
\vspace{0.25in}

\begin{center}
\begin{tikzpicture}
    \node at (-11,0) {} [grow'=right,
        level distance=8mm,
        level 3/.style={sibling distance=10mm}]
        child {node {$x$}
            child { node {$sp$}
                child { node {$y$} child child }
                child }
                child
                };
    \draw[thick,->] (-7,0)--(-6,0);
    \node at (-5.5,0) {} [grow'=right,
        level distance=8mm]
        child {node {$x \circ_{1} p$}
            child { node {$y$} child child }
            child };

    \node at (-2,0) {$-$};
    \node at (-1.5,0) {} [grow'=right,
        level distance=8mm]
        child { node {$x$}
            child { node {$p \circ_{1} y$}
                child child }
            child };

\end{tikzpicture}

\end{center}

\vspace{0.25in}

Note that twistability implies that if $\op M$ is an $(\op P', \op P)$-bimodule and $\op N$ is an $(\op P, \op P'')$-bimodule, then $\opbar (\op M, \op P, \op N)$ is an $(\op P', \op P'')$-bimodule.

Let $\op{P}$ be a connected, augmented operad.  Set
$$\widetilde{\mathbf E}\op{P} = \opbar(\op{P},\op{P},\op{J})\in {}_{\op P}\cat {Mix}_{\opbar \op P}$$ and
$$\mathbf E\op{P} = \opbar(\op{J},\op{P},\op{P})\in {}_{\opbar\op P}\cat {Mix}_{\op P}.$$
We define a natural morphism of $\opbar \op P$-bicomodules
$$\delta_{\op{P}} : \opbar(\op{P}) \rightarrow \mathbf E\op{P} \circ_{\op{P}} \widetilde{\mathbf {E}}\op{P}.$$  It suffices to construct a natural map $\hat{\delta}_{\op{P}} : \opbar(\op{P}) \rightarrow \opbar(\op{P}) \circ \op{P} \circ \opbar(\op{P})$ such that $\im(\hat{\delta}_{\op{P}}) \subseteq \ker(d_{R} \circ 1 + 1 \circ d_{L})$, where $d_{R}$ is the component of the differential in $\mathbf E\op{P}$ that comes from the left action of $\op{P}$ on itself, and $d_{L}$ is the component of the differential in $\widetilde{\mathbf E}\op{P}$ that comes from the right action of $\op{P}$ on itself. As a morphism of symmetric sequences, forgetting differentials, $\hat{\delta}_{\op{P}}$ is the composite
\[
    \opbar\op{P} \xrightarrow{\Delta} \opbar\op{P} \circ \opbar\op{P} \xrightarrow{\cong} \opbar\op{P} \circ \op{J} \circ \opbar\op{P} \xrightarrow{1 \circ \eta \circ 1} \opbar\op{P} \circ \op{P} \circ \opbar\op{P}.
\]
Concretely, if $\Delta(\tau) = \sum \sigma \otimes \{ \sigma_{1} \otimes \cdots \otimes \sigma_{m}\}$, where $\tau$ is a tree with vertices labeled with elements of $s\tilde{\op{P}}$, then
\[
    \hat{\delta}_{\op{P}}(\tau) = \sum \sigma \otimes 1^{\otimes m} \otimes \{ \sigma_{1} \otimes \cdots \otimes \sigma_{m} \}.
\]

\begin{prop}
$\im(\hat{\delta}_{\op{P}}) \subseteq \ker(d_{R} \circ 1 + 1 \circ d_{L})$.
\end{prop}

\begin{proof}
We show that there is a one-to-one correspondence between the terms of $(1 \circ d_{L})\hat{\delta}_{\op{P}}$ and the terms of $(d_{R}\circ 1)\hat{\delta}_{\op{P}}$, with opposite signs.

The image of $\hat{\delta}$ is spanned by elements of the form
\[
    \sigma \otimes \{ 1^{\otimes n} \} \otimes \{ \sigma_{1} \otimes \cdots \otimes \sigma_{n} \},
\]
where $\sigma, \sigma_{1},\ldots,\sigma_{n} \in \opbar\op{P}$.

From the definitions,
\begin{multline*}
    1 \circ d_{L} (\sigma \otimes 1^{\otimes n} \otimes \{ \sigma_{1} \otimes \cdots \sigma_{n} \})\\
    = \sum_{i=1}^{n} \sigma \otimes \{ 1^{\otimes(i-1)} \otimes p_{i} \otimes 1^{\otimes(n-i)} \} \otimes \{ \sigma_{1} \otimes \cdots \otimes (\rho_{1} \otimes \cdots \otimes \rho_{m_{i}}) \otimes \cdots \otimes \sigma_{n} \}
\end{multline*}
where $\sigma_{i}$ has root labeled by $p_{i}$ of arity $m_{i}$, that roots the subtrees $\rho_{1},\ldots,\rho_{m_{i}}$.

Let $\sigma'$ be the subtree generated by $\sigma$ and the additional vertex $p_{i}$ from the $i$th term above.  Since $\sigma$ is an $n$-tree and $p_{i}$ has arity $m_{i}$, $\sigma'$ has arity $n + m_{i}-1$.  We find, in $\hat{\delta}_{\op{P}}(\tau)$, the term
\[
    \Phi = \sigma' \otimes 1^{\otimes (n+m_{i}-1)} \otimes \{ \sigma_{1} \otimes \cdots \otimes (\rho_{1} \otimes \cdots \otimes \rho_{m_{i}}) \otimes \cdots \sigma_{n} \}.
\]
Recall that $d_{R}(\sigma')$ has one component for each vertex which is the target only of ingoing edges; $p_{i}$ is one such vertex.  The contribution to the differential is $-\sigma \otimes 1^{\otimes(i-1)} \otimes p_{i}(1,\ldots,1) \otimes 1^{\otimes(n-i)}$; thus $(d_{R}\circ 1)(\Phi)$ contains the term
\[
    - \sigma \otimes \{ 1^{\otimes(i-1)} \otimes p_{i} \otimes 1^{\otimes(n-i)} \} \otimes \{ \sigma_{1} \otimes \cdots \otimes (\rho_{1} \otimes \cdots \otimes \rho_{m_{i}}) \otimes \cdots \sigma_{n} \}
\]
that cancels with the $i$th term in $1 \circ d_{L} (\sigma \otimes 1^{\otimes n} \otimes \{ \sigma_{1} \otimes \cdots \sigma_{n} \})$.
\end{proof}

To construct the natural morphism of $\op P$-bimodules
$$\mu_{\op{P}} : \widetilde{\mathbf E}\op{P} \underset {\opbar{\op{P}}}\square \mathbf E\op{P} \rightarrow \op P,$$
it suffices to construct a map $\hat{\mu}_{\op{P}} : \op{P} \circ \opbar\op{P} \circ \op{P} \rightarrow \op{P}$ in such a way that $\im(d_{L} \circ 1 + 1 \circ d_{R}) \subseteq \ker{\hat{\mu}_{\op{P}}}$.  We define $\hat{\mu}_{\op{P}}$ as the composite
\[
    \op{P} \circ \opbar\op{P} \circ \op{P} \xrightarrow{1 \circ \epsilon \circ 1} \op{P} \circ \op{J} \circ \op{P} \xrightarrow{\cong} \op{P} \circ \op{P} \xrightarrow{\gamma} \op{P}.
\]

\begin{prop}
$\im(d_{L} \circ 1 + 1 \circ d_{R}) \subseteq \ker{\hat{\mu}_{\op{P}}}$.
\end{prop}

\begin{proof}
Since the differentials reduce weight by one, and the augmentation $\epsilon:\opbar\op{P} \rightarrow \op{J}$ kills everything of nonzero weight, we only need to concern ourselves with elements of $\op{P} \circ \opbar\op{P} \circ \op{P}$ of the form
\[
    \Phi = p \otimes \{ 1^{\otimes(i-1)} \otimes \tau_{i} \otimes 1^{\otimes(n-i)}\} \otimes \{q_{1} \otimes \cdots \otimes \mathbf{q}_{i} \otimes \cdots \otimes q_{n}\}
\]
where $\tau_{i}$ is represented by an $m$-tree with one vertex labeled by $sr \in s\tilde{P}(m)$ and $\mathbf{q}_{i} = q_{i1} \otimes \cdots \otimes q_{im}$. From the definitions,
\[
    \hat{\mu}_{\op{P}}d_{L}\Phi = (p \circ_{i} r) (q_{1}, \ldots, \mathbf{q}_{i}, \ldots, q_{n})
\]
while
\[
    \hat{\mu}_{\op{P}}d_{R}\Phi = -p(q_{1},\ldots,r(\mathbf{q}_{i}),\ldots,q_{n}).
\]
By associativity of the composition product $\gamma$, $\hat{\mu}_{\op{P}}d_{L}\Phi+\hat{\mu}_{\op{P}}d_{R}\Phi = 0$ as desired.
\end{proof}

\begin{thm}
The category $\cat{dgProj}^\Sigma$ admits a right twisting structure $(\mathbf B, \mathbf E, \widetilde {\mathbf E}, \delta, \mu)$.
\end{thm}

\begin{proof}
{Having identified}
$$\widetilde{\mathbf E}\op{P} \underset {\opbar{\op{P}}}\square \mathbf E\op{P} \cong \op{P} \circ \opbar{\op{P}} \circ \op{P}\quad\text{and}\quad \mathbf E\op{P} \circ_{\op{P}} \widetilde{\mathbf E}\op{P} \cong \opbar{\op{P}} \circ \op{P} \circ \opbar{\op{P}}$$
as above, the verification that the two diagrams of Definition~\ref{defn:twisting} commute is a straightforward diagram chase using  that the diagonal   $\Delta$ of $\opbar{\op{P}}$ is counital and the composition product  $\gamma$  of $\op{P}$ is unital.
\end{proof}

{Viewing chain complexes as symmetric sequences of chain complexes concentrated in arity $0$, we obtain the following immediate consequence of the theorem above.}

 \begin{cor}  The category $\cat{dgProj}$ admits a right twisting structure.
 \end{cor}

 {The twisting structure on $\cat {dgM}^\Sigma$ induces important adjunctions on the level of (co)algebras over (co)operads.} Recall the definition of twisted products with respect to classifying morphisms (Definition \ref{defn:classifying-twisting}).

 \begin{defn}\label{defn:bar-cobar} Let $\op Q$ be a cooperad, and let $\op P$ be an operad, both in $\cat {dgProj}^\Sigma$.  Let $g:\op Q \to \mathbf B P$ be a classifying morphism.  The \emph{$g$-cobar construction}
 \[
    \Omega_{g}:{}_{\op Q}\cat{Comod} \rightarrow {}_{\op{P}}\cat{Mod}
 \]
 and the \emph{$g$-bar construction}
\[
    B_{g}:  {}_{\op{P}}\cat{Mod} \rightarrow {}_{\op Q}\cat{Comod}
\]
are given by $\Omega_{g}\op{M} = \op{P} \circ_{g} \op{M}$ and $B_{g}\op{N} = \op{Q} \circ_{g} \op{N}$.
 \end{defn}

\begin{rmk}
Since $z : \cat{dgProj} \rightarrow \cat{dgProj}^{\Sigma}$ is fully faithful when restricted to both $\alg{P}$ and $\coalg{Q}$, the functors of Definition~\ref{defn:bar-cobar} restrict to define
\[
    \Omega_{g} : \coalg{Q} \rightarrow \alg{P} \quad \text{and} \quad B_{g} : \alg{P} \rightarrow \coalg{Q}.
\]
\end{rmk}

 The next result follows immediately from Corollary \ref{cor:adjunct-oneside}, but is important enough to be formulated as a separate statement.

 \begin{prop}\label{prop:gbar-gcobar}
 Let $\op Q$ be a cooperad, and let $\op P$ be an operad, both in $\cat {dgProj}^\Sigma$.  For any classifying morphism $g:\op Q \to \mathbf B P$, the $g$-cobar construction $\Om _{g}$ is left adjoint to the $g$-bar construction $B _{g}$.
 \end{prop}

 We prove in the next section that when $g$ is the canonical classifying morphism of a quadratic operad, then $(\Om_{g},B_{g})$ is the usual cobar/bar adjunction.

\section{Categories with morphisms up to strong homotopy}\label{sec:shp}

In this section, we consider classifying morphisms $g:\op{Q} \rightarrow \opbar\op{P}$ with the property that the counit of the associated standard construction, $\epsilon:K(g) \rightarrow \op{P}$, is a quasi-isomorphism.  This is the case with quadratic Koszul operads and their resolutions; we show that the corresponding Kleisli categories are isomorphic to the classic ``strong homotopy'' categories.  In general, the following argument, translated from Markl~\cite{markl:04} where it is presented in the language of coloured operads, shows that the morphism sets of the Kleisli category are homotopy-invariant, and {therefore} model categories of algebras and morphisms up to strong homotopy.

\begin{notn}
If $\op{Q}$ is a coaugmented cooperad in $\cat{dgProj}^\Sigma$ and $g: \op{Q} \rightarrow \opbar\op{P}$ is a classifying morphism, then the coaugmentation in $\op{Q}$, along with the unit in $\op{P}$, define a coaugmentation $\eta_{g}:\op{P} \rightarrow K(g)$.  Let $F:K(g) \circ_{\op P} \op{M} \rightarrow \op{N}$ be a morphism of left $\op{P}$-modules.  The \emph{underlying morphism $F_{0}:\op{M} \rightarrow\op{N}$ associated to $F$} is the composite,
\[
    \op{M} \cong \op{P} \circ_{\op{P}} \op{M} \xrightarrow{\eta_{g}\circ 1} K(g) \circ_{\op{P}} \op{M} \xrightarrow{F} \op{N}.
\]
\end{notn}

\begin{prop}\cite[Proposition 35]{markl:04}\label{prop:htyp-inv}
Let $\op{Q}$ be a coaugmented cooperad and $\op{P}$ an augmented operad, both in $\cat{dgProj}^\Sigma$. Let $g:\op{Q} \rightarrow \opbar\op{P}$ be a classifying morphism such that $\epsilon:K(g) \rightarrow \op{P}$ is a surjective quasi-isomorphism. Let $F:K(g)\circ_{\op{P}}\op{M} \rightarrow \op{N}$ be a morphism of left $\op{P}$-modules.

 If $f:\op{M} \rightarrow \op{N}$ is a morphism of symmetric sequences homotopic to $F_{0}$, then $f$ is the underlying morphism of a morphism of left $\op{P}$-modules, $f:K(g)\circ_{\op{P}}\op{M} \rightarrow \op{N}$.
\end{prop}

\begin{proof}
Let $\op{I}$ be the symmetric sequence concentrated in arity $1$, with $\op{I}(1) = R\{e_{0},e_{1},s\}$ and $\partial s = e_{1} - e_{0}$. Then $F_{0}$ and $f$ are chain homotopic if and only if there exists a morphism of symmetric sequences $\Phi : \op{I} \circ \op{X} \rightarrow \op{Y}$ such that $\Phi(e_{0} \otimes -) = f$ and $\Phi(e_{1} \otimes -)= F_{0}$.

 The symmetric sequence $\op{I}$ is a coaugmented cooperad.  Indeed, the diagonal is defined by $\psi(e_{i}) = e_{i} \otimes e_{i}$ for $i = 0,1$, and $\psi(s)=s \otimes e_{1} + e_{0} \otimes s$.  The {counit} is defined by $\epsilon(e_{i})=1$ for $i=0,1$ and $\epsilon(s) = 0$. The coaugmentation is defined by $\eta(1)=e_{0}$.

 Since $\op{I}$ is concentrated in arity $1$, $\op{I} \circ \op{Q}$ is also a coaugmented cooperad.  The classifying morphism $g$ extends to a classifying morphism $g':\op{I} \circ \op{Q} \rightarrow B\op{P}$, via $g'(e_{i} \otimes -) = g$ for $i=0,1$, and $g'(s \otimes -) = 0$.  Define the coaugmentation by $\eta'(1) = e_{0} \otimes \eta(1)$  If $\epsilon : K(g) \rightarrow \op{P}$ is a surjective quasi-isomorphism, then so too is $\epsilon':K(g')\rightarrow \op{P}$.

Let $\Hom_{\Sigma}(\op{M},\op{N})$ be the internal morphism symmetric sequence in $\cat{dgM}$; that is, $\Hom_{\Sigma}(\op{M},\op{N})(n) = \prod_{m}\Hom_{R[\Sigma_{m}]}(\op{M}^{\odot n}(m),\op{N}(m))$.  It suffices to construct a $\op{P}$-bimodule morphism $\widetilde F^{\sharp} : K(g) \rightarrow \Hom_{\Sigma}(\op{M},\op{N})$ such that $\widetilde F^{\sharp}\eta = f^{\sharp}$, where $f^{\sharp}:\op{J} \rightarrow \Hom_{\Sigma}(\op{M},\op{N})$ is the right adjoint of $f$.  We will then have that $\widetilde F$, the left adjoint of $\widetilde F^{\sharp}$, has $f$ as its underlying morphism.

Define $H:\op{I} \circ \op{Q} \rightarrow \Hom_{\Sigma}(\op{M},\op{N})$ by $H(e_{0} \otimes \eta(r))=f^{\sharp}(r)$ for $r \in \op{J}$, $H(e_{0} \otimes q) = 0$ otherwise; $H(e_{1} \otimes q) = F^{\sharp}(q)$ for all $q \in \op{Q}$; and $H(s \otimes q) = 0$ for all $q \in \op{Q}$.  Then $H$ extends uniquely to a $\op{P}$-bimodule morphism $K(g') \rightarrow \Hom_{\Sigma}(\op{M},\op{N})$, also called $H$.

We note that $K(g)$ is cellular, in the sense that it is almost free as a $\op{P}$-bimodule, and has an increasing filtration such that the differential strictly reduces filtration degree.  Furthermore, since we are assuming that $\op{P}$ and $\op{Q}$ are projective as symmetric sequences, so too is $K(g)$.  Therefore the solid square
\[
    \xymatrix{
        \op{J} \ar[r]^{\eta'} \ar[d]_{\eta} & K(g') \ar@{->>}[d]^{\epsilon'}_{\sim} \ar[r]^-{H} & \Hom_{\Sigma}(\op{M},\op{N}) \\
        K(g) \ar[r]_{\epsilon} \ar@{.>}[ur]^{\sigma} & \op{P}
    }
\]
has the lifting $\sigma$ as indicated by the dotted arrow.  We set $\widetilde F^\sharp = H\sigma$ to complete the proof.
\end{proof}

\subsection{Strongly homotopy morphisms of $\op P$-algebras}

{We now} use the right twisting structure for symmetric sequences of chain complexes to provide two operadic descriptions of strongly homotopy morphisms of algebras and of coalgebras over a quadratic operad.  We begin by recalling the basic theory of weight-graded operads.

\subsubsection{Weight-graded operads}\label{ssec:weight}

Weight-graded operads, first considered by Fresse in~\cite{fresse:04}, allow the construction of quadratic dual cooperads and Koszul resolutions for non-quadratic operads.  A \emph{weight grading} on a {chain complex} $X$ is a grading $X \cong \bigoplus_{s \geq 0} X_{(s)}$, where the differential preserves each $X_{(s)}$. If $X$ and $Y$ are weight-graded, then $X \otimes Y$ is weight-graded, with $(X \otimes Y)_{(n)} = \bigoplus_{s+t=n} X_{(s)} \otimes Y_{(t)}$, so weight-graded chain complexes form a symmetric monoidal category in which we can define symmetric sequences and operads.

If $\op{P}$ is any connected operad, then we can give it the \emph{canonical weight grading}, defined by
\[
    \op{P}_{(r)}(n) = \left\{
        \begin{array}{cl}
            \op{P}(n)  & r = n-1 \\
            0 & \text{otherwise.}
        \end{array}
        \right.
\]
If $\op{M}$ is a weight-graded symmetric sequence, then  {its} grading {induces weight gradings on}  the free operad and cooperad on $\op{M}$.  Essentially the same argument that shows that $T^{d}(V)_{(d)} = T^{d}(V_{1})$ and $T^{>d}(V)_{(d)} = 0$ for the tensor algebra on a connected {chain complex} $V$, establishes that $F^{c}_{(r)}(\op{M})_{(r)} = F^{c}_{(r)}(\op{M}_{(1)})$ and $F^{c}_{(>r)}(\op{M})_{(r)}=0$.  (Recall that $F^{c}_{(r)}(\op{M})$ is the sub symmetric sequence generated by trees with $r$ vertices.)

The bar construction $\opbar\op{P}$  of a weight-graded operad is therefore naturally weight-graded.  Set
\[
    \op{P}^{\bot}_{(s)} = H_{s}(\opbar_{s}(\op{P})_{(s)}),
\]
where the homological grading is with respect to the bar wordlength, that is, the number of vertices in a representative tree.  Since $F^{c}_{(>s)}(s\tilde{\op{P}})_{(s)}=0$, we have that
\[
    \op{P}^{\bot}_{(s)} = \ker(d:\opbar_{s}(\op{P})_{(s)} \rightarrow \opbar_{s-1}(\op{P})_{(s)}),
\]
and so $\op{P}^{\bot}$ is a sub weight-graded symmetric sequence of $\opbar\op{P}$.  Fresse proves that in fact, $\op{P}^{\bot}$ is a sub weight-graded {(in fact, quadratic)} cooperad of $B\op{P}$.  The inclusion
$$\kappa_{\op P} : \op{P}^{\bot} \rightarrow \opbar{\op{P}}$$
 is the \emph{canonical classifying morphism} of $\op{P}$.

\begin{rmk}
We take this opportunity to observe that $\op{P}^{\bot}(2) \cong s\tilde{\op{P}}(2)$.  If $\op{P}$ is quadratic and $\op{P}^{!}$ is its quadratic dual~\cite{ginzburg-kapranov:94}, then $\op{P}^{\bot} \cong \op{S} \otimes \op{P}^{!\dual}$.
\end{rmk}

If $\op{M}$ is a left $\op{P}$-module, then the $\op{P}$-\emph{bar construction} of $\op{M}$ is the cofree left $\op{P}^{\bot}$-module, $B_{\op{P}}(\op{M}) = \op{P}^{\bot} \circ \op{M}$, with internal differential from $\op{P}^{\bot}$ and $\op{M}$, perturbed by the unique coderivation determined by the composite
\[
    \op{P}^{\bot}(2) \otimes \op{M}(k) \otimes \op{M}(\ell) \xrightarrow{s^{-1} \otimes 1 \otimes 1} \op{P}(2) \otimes \op{M}(k) \otimes \op{M}(\ell) \rightarrow \op{M}(k+\ell).
\]
Observe that this makes sense for a general weight-graded operad, since $\op P^\bot$ is always quadratic.  For a full treatment of bar and cobar constructions, see for example Fox and Markl,~\cite{fox-markl:97}, or Getzler and Jones~\cite{getzler-jones:94}.

{If} $\op{Q}$ {is} a weight-graded cooperad{, then} the {\emph{operadic cobar construction on $\op Q$}} [ref],  $\mathbf{\Omega}\op{Q}$, is a weight-graded operad.  We define
\[
    \op{Q}^{\bot}_{(s)} = \coker(d:\mathbf{\Omega}^{s-1}\op{Q}_{(s)} \rightarrow \mathbf{\Omega}^{s}\op{Q}_{(s)}).
\]
Fresse~\cite{fresse:04} shows that $\op{Q}^{\bot}$ is a weight-graded {(in fact, quadratic)} quotient operad of $\mathbf{\Omega}\op{Q}$.

Let $\op{M}$ be a left $\op{Q}$-comodule.  The $\op{Q}$-\emph{cobar construction} of $\op{M}$ is the free $\op{Q}^{\bot}$-module, $\Omega_{\op{Q}}(\op{M}) = \op{Q}^{\bot} \circ \op{M}$, whose internal differential is perturbed by the unique derivation determined by
\[
    C \rightarrow \op{Q}(2) \otimes C \otimes C \rightarrow \op{Q}^{\bot}(2) \otimes C \otimes C.
\]
Again,  $\op Q^\perp$ is quadratic because it is a quotient of a quadratic operad, so the above equation makes sense.

Recall~\cite[5.2.8]{fresse:04},~\cite[Definition 2.23]{getzler-jones:94} that a weight-graded operad is called \emph{Koszul} if the canonical classifying morphism $\kappa_{\op P}:\op{P}^{\bot} \rightarrow \opbar\op{P}$ is a quasi-isomorphism.  Recall moreover the bar and cobar constructions of Definition \ref{defn:bar-cobar}.

\begin{prop}\label{prop:bar-cobar}  Let $\op P$ be a weight-graded Koszul operad in $\cat{dgProj}^\Sigma$ with canonical classifying morphism $\kappa_{\op P}:\op{P}^{\bot} \rightarrow \opbar \op{P}$. Then
\[
    \Omega_{\op{P}^{\bot}}=\Omega_{\kappa_{\op P}}\quad\text{and}\quad B_{\op{P}}=B_{\kappa_{\op P}}.
\]
\end{prop}

\begin{proof}
For brevity, we set $\op{Q} = \op{P}^{\bot}$. The isomorphism $\op{Q} \cong \op{Q} \underset{\opbar{\op{P}}}\square\opbar{\op{P}}$ is the corestriction of the  composite
\[
    \op{Q} \xrightarrow{\psi} \op{Q} \circ \op{Q} \xrightarrow{1 \circ \kappa} \op{Q} \circ \opbar{\op{P}}.
\]
Let $\op{M}$ be a left $\op{P}$-module. We observe that
$$\mathbf E\op{P} \circ_{\op{P}} \op M = \opbar(\op{J};\op{P};\op{P})\circ_{\op{P}}\op M \cong \opbar(\op{J};\op{P};\op{M})$$ by~\cite[4.1.2]{fresse:04}. Thus is suffices to show that the map
\[
    \op{Q} \circ \op{M} \rightarrow \op{Q} \circ \opbar(\op{J};\op{P};\op{M})
\]
commutes with differentials.  The map {obviously} commutes with the internal differential in $\op{Q}$, if any, and with the differential in $\opbar{\op{P}}$.  Since $\im\kappa_{\op P} \subset F^{c}(s\tilde{P}(2))$, and $\op{Q}$ is quadratic, it suffices to show that
\[
    \xymatrix{
        \op{Q}(2) \otimes \op{M}(n_1) \otimes \op{M}(n_2)
            \ar[r]^{\kappa \otimes 1 \otimes 1}
            \ar[d]^{d}
        & \opbar{\op{P}}(2) \otimes \op{M}(n_1) \otimes \op{M}(n_2)
            \ar[d]^{d}
        \\
        \op{Q}(1) \otimes \op{M}(n_{1}+n_{2}) \ar[r]_{\kappa \otimes 1}
        & \opbar{\op{P}}(1) \otimes \op{M}(n_{1}+n_{2})
    }
\]
commutes.  From the definitions, $\op{Q}(2) = \opbar{\op{P}}(2) = s\tilde{P}(2)$, and both differentials come from desuspending then applying the action of $\op{P}$ on $\op{M}$, so the diagram commutes.

The proof for cobar constructions is similar, and exploits the fact that $(\op{P}^{\bot})^{\bot} = \op{P}$ for Koszul operads $\op{P}$.
\end{proof}

Propositions~\ref{prop:bar-cobar} and~\ref{prop:gbar-gcobar} yield another proof of the following well known result.

\begin{prop}\label{prop:bar-cobar-adjoint}~\cite[Theorem 2.25]{getzler-jones:94}
Let $\op{P}$ be a Koszul operad in $\cat{dgProj}^\Sigma$. Then the $\op{P}^{\bot}$-cobar and $\op{P}$-bar
constructions are adjoint functors,
\[
    \Omega_{\op{P}^{\bot}} : \coalg{P^{\bot}} \rightleftarrows \alg{P} :
    B_{\op{P}}.
\]
\end{prop}

\subsubsection{The Kleisli category description of strong homotopy}

In this section, we observe that the standard construction, $K(\op{P})$, for the canonical classifying morphism, $\kappa_{\op P}:\op{P}^{\perp}\rightarrow \opbar\op{P}$ is the two-sided Koszul resolution of $\op{P}$, for a weight-graded operad $\op{P}$.  We then show that the Kleisli category associated to the comonad determined by $K(\op{P})$ is isomorphic to the category of $\op{P}$-algebras and strongly homotopy $\op{P}$-algebra morphisms.

\begin{defn}  Let $\op P$ be a weight-graded Koszul operad. The objects of the category $\alg{\op P}_{sh}$ are all $\op P$-algebras, while the morphisms are \emph{strongly homotopy morphisms} of $\op P$-algebras, i.e.,
$$\alg{\op P}_{sh}(A,E)=\coalg{P^{\bot}}(B_{\op P}A,B_{\op P}E).$$
Dually, the objects of the category $\coalg{P^{\bot}}_{sh}$ are all $\op{P}^{\bot}$-coalgebras, while the morphisms are \emph{strongly homotopy morphisms} of $\op{P}^{\bot}$-coalgebras, i.e.,
$$\coalg{P^{\bot}}_{sh}(A,E)=\alg{\op P}(\Omega_{\op P^\bot}A,\Omega_{\op P^\bot}E).$$
\end{defn}

Recall the canonical classifying morphism $\kappa_{\op P}:\op{P}^{\bot} \rightarrow \opbar\op{P}$ constructed in Section~\ref{ssec:weight} and the standard construction $K(g)$ for any classifying morphism $g$, from Definition~\ref{defn:std-constr}.  {It is easy to see from the definitions that} the underlying $\op{P}$-bimodule of the standard construction $K(\kappa_{\op P})$ is the two-sided Koszul resolution of $\op{P}$; see Fresse~\cite{fresse:04}.  For this reason, we write $K(\op{P})$ instead of $K(\kappa_{\op P})$ {and are} motivated to make the following definition.

\begin{defn}  The \emph{two-sided dual Koszul resolution} of $\op{P}^{\bot}$ is the dual standard construction applied to the classifying morphism $\kappa_{\op P}$:
 \[
    T(\op P)=\op{P}^{\bot} \circ_{\kappa_{\op P}} \op{P} \circ_{\kappa_{\op P}} \op{P}^{\bot},
\]
 which is a $\op{P}^{\bot}$-bicomodule.
 \end{defn}

 \begin{prop}\label{prop:quad-co-ring} If $\op P$ is a Koszul operad in $\cat{dgProj}^\Sigma$, then $K(\op P)$ is a $\op P$-co-ring and $T(\op P)$ is a $\op{P}^{\bot}$-ring.
 \end{prop}

 \begin{proof} The proposition is an instance of Proposition ~\ref{prop:std-coring}, for the classifying morphism $\kappa_{\op P}:\op{P}^{\bot} \rightarrow \opbar\op{P}$.
 \end{proof}

{ \begin{notn}  Let
 \[
    K_{\op P}=K(\op P)\underset {\op P}\circ -:{}_{\op P}\cat{Mod}\to {}_{\op P}\cat {Mod}.
 \]
 Observe that $K_{\op P}$ is (the underlying functor of) a comonad, since $K(\op P)$ is a co-ring. Moreover there is an induced comonad on $\alg{\op P}$, denoted $\mathsf K_{\op P}$, with underlying endofunctor $K(\op P)\underset {\op P}\circ z(-)$.

 Similarly, there is a monad on ${}_{\op{P}^{\bot}}\cat {Comod}$ with underlying endofunctor
 \[
    T_{\op P}=T(\op P)\underset {\op{P}^{\bot}}\square-:{}_{\op{P}^{\bot}}\cat {Comod}\to {}_{\op{P}^{\bot}}\cat {Comod},
 \]
which in turn restricts and corestricts to a monad on $\coalg{\op P^\bot}$, denoted $\mathsf T_{\op P}$, with underlying endofunctor $T(\op P)\underset {\op{P}^{\bot}}\square z(-)$.
 \end{notn}}

 Recall the definition of the (co)Kleisli category determined by a (co)monad (Notation \ref{notn:intro}).
{
 \begin{thm}\label{thm:kleisli}
If $\op P$ is a Koszul operad in $\cat{dgProj}^\Sigma$, then there are isomorphisms of categories
 $$\Theta_{\op P}:\alg{\op P}_{sh} \xrightarrow\cong  {}_{\mathsf K_{\op P}}\alg{\op P}$$
 and
 $$\Upsilon_{\op P}:\coalg{\op{P}^{\bot}}_{sh} \xrightarrow\cong \coalg{\op P^\bot}_{\mathsf T_{\op P}},$$
 which are the identity on objects.
 \end{thm}}

\begin{proof}  The definition of $\Theta_{\op P}$ on morphisms follows from the sequence of natural isomorphisms below, where $A$ and $B$ are $\op P$-algebras.
\begin{align*}
    {}_{\op P}&\cat {Mod}\big(K(\op P)\underset {\op P}\circ z(A), z(E)\big) \\
        &\cong{}_{\op P}\cat {Mod}\big(\op{P} \circ_{\kappa_{\op P}} \op{P}^{\bot} \circ_{\kappa_{\op P}} z(A),z(E)\big) \\
        &\cong{}_{\op P^{\bot}}\cat{Comod}\big( \op P^\perp \circ_{\kappa_{\op P}} z(A), \op P^\perp \circ_{\kappa_{\op P}} z(E)\big)\\
        &\cong{}_{\op{P}^{\bot}}\cat{Comod}\big( z(B_{\op P}A), z(B_{\op P}E)\big) \\
        &\cong \coalg{\op{P}^{\bot}}(B_{\op P}A,B_{\op P}E).
\end{align*}

{Similarly, the definition of the functor $\Upsilon_{\op P}$ on morphisms follows} from the sequence of isomorphisms below, where $C$ and $D$ are $\op{P}^{\bot}$-coalgebras.
\begin{align*}
    {}_{\op{P}^{\bot}}&\cat {Comod}\big( z(C), T(\op P)\underset {\op{P}^{\bot}}\square z(D)\big)\\
        &\cong {}_{\op P^\perp}\cat {Comod}\big( z(C), \op P^\perp \circ_{\kappa_{\op P}} \op{P} \circ_{\kappa_{\op P}} z(D)\big)\\
        &\cong {}_{\op P}\cat {Mod}\big( \op{P} \circ_{\kappa_{\op P}}  z(C), \op{P} \circ_{\kappa_{\op P}} z(D)\big)\\
        &\cong{}_{\op P}\cat {Mod}\big( z(\Omega_{\op{P}^{\bot}}C), z(\Omega_{\op{P}^{\bot}}D)\big)\\
        &\cong \alg{\op P}(\Om_{\op{P}^{\bot}}C, \Om_{\op{P}^{\bot}}D).
\end{align*}

Note that in both of the sequences of isomorphisms above, Corollary \ref{cor:adjunct-oneside} and Proposition \ref{prop:bar-cobar} play a key role.
\end{proof}

 \subsection{Strongly homotopy morphisms of $\op P_{\infty}$-algebras, the Koszul case}

 Let $\op{P}$ be an operad.  A $\op{P}_{\infty}$-algebra is defined to be an {algebra} over a cofibrant replacement $\op{P}'$ of $\op{P}$.  If $\op{P}$ is Koszul, then we may choose our cofibrant replacement to be $\mathbf{\Omega}\op{P}^{\bot}$, the operadic cobar construction on the quadratic dual cooperad $\op{P}^{\bot}$~\cite[4.2.14]{ginzburg-kapranov:94}. The unit of the {operadic} bar-cobar adjunction defines a classifying morphism,
 $$\eta_{\op P^\bot} : \op{P}^{\bot} \rightarrow \opbar\mathbf{\Omega}\op{P}^{\bot}.$$
 {Observe that the associated standard construction, $K(\eta_{\op P ^\bot})$, is a $\Omega \op P^\bot$-co-ring, by Proposition \ref{prop:std-coring}.}

 Let $A$ be a $\mathbf{\Omega}\op{P}^{\bot}$-algebra.  By analogy with the homotopy-associative case~\cite{stasheff:63}, we define the \emph{bar-tilde construction} $\widetilde{B}_{\op{P}}(A)$ to be the cofree $\op P^\bot$-coalgebra $\Gamma_{\op{P}^{\bot}}(A)$, with differential perturbed by the composition
 \[
    \op{P}^{\bot}(n) \otimes A^{\otimes n} \rightarrow s^{-1}\op{P}^{\bot}(n) \otimes A^{\otimes n} \hookrightarrow \mathbf{\Omega}(\op{P}^{\bot})(n) \otimes A^{\otimes n} \rightarrow A.
 \]
 Here we are using the fact that the operad $\mathbf{\Omega}\op{P}^{\bot}$ is generated by the symmetric sequence $s^{-1}\op{P}^{\bot}$.

\begin{prop}
Let $\op P$ be a Koszul operad in $\cat{dgProj}^\Sigma$. If $A$ is a $\mathbf{\Omega}\op{P}^{\bot}$-algebra, then $B_{\eta_{\op P^\bot}}(A)= \widetilde{B}_{\op{P}}(A)$.
\end{prop}

\begin{proof}
Definitions \ref{defn:classifying-twisting} and  \ref{defn:bar-cobar} imply that
$$B_{\eta_{\op P^\bot}}(A) = \op P^\bot \circ_{\eta_{\op P^\bot}} z(A)=\op P^\bot \underset {\opbar\mathbf{\Omega}\op{P}^{\bot}}\square \mathbf E\mathbf{\Omega}\op{P}^{\bot}\underset {\mathbf{\Omega}\op{P}^{\bot}}\circ z(A)\cong \widetilde{B}_{\op{P}}(A).$$
\end{proof}

Let $\alg{P}_{\infty}$ denote the category of $\op{P}_{\infty}$-algebras and morphisms up to strong homotopy, i.e., the full subcategory  of $\coalg{\op{P}^{\bot}}$ spanned by the almost-cofree coalgebras.  (A DG coalgebra is \emph{almost cofree} if the underlying coalgebra of graded modules is cofree.)

\begin{thm}\label{thm:p-infty} Let $\op P$ be a Koszul operad in $\cat{dgProj}^\Sigma$.
The $\eta_{\op P^\bot}$-bar construction {induces} an isomorphism of categories,
{$$B_{\eta_{\op P^\bot}}:{}_{\mathsf K(\eta_{\op P^\bot})}\big(\alg{\mathbf{\Omega} \op P^\bot}\big) \to \alg{P}_{\infty},$$
where $\mathsf K(\eta_{\op P^\bot})$ denotes the comonad on $\alg{\mathbf{\Omega} \op P^\bot}$ with underlying endofunctor $K(\eta_{\op P^\bot})\underset {\mathbf{\Omega}\op P^\bot}\circ z(-)$.}
\end{thm}

\begin{proof} If $A$ and $E$ are $\Omega\op P^\bot$-algebras, then
\begin{align*}
{}_{\mathsf K(\eta_{\op P^\bot})}\big(\alg{\mathbf{\Omega} \op P^\bot}\big)(A,E)&={}_{\mathbf{\Omega} \op P^\bot}\cat {Mod}\big(K(\eta_{\op P^\bot})\circ_{\eta_{\op P^\bot}} z(A), z(E)\big)\\
&={}_{\mathbf{\Omega} \op P^\bot}\cat {Mod}\big(\mathbf{\Omega} \op P^\bot\circ_{\eta_{\op P^\bot}} \op P^\bot \circ_{\eta_{\op P^\bot}} z(A), z(E)\big)\\
&\overset {(\star)}\cong{}_{\op P^\bot}\cat {Comod}\big( \op P^\bot \circ_{\eta_{\op P^\bot}} z(A),  \op P^\bot\circ_{\eta_{\op P^\bot}}z(E)\big)\\
&=\coalg{\op P^\bot}\big(B_{\eta_{\op P^\bot}}(A),B_{\eta_{\op P^\bot}}(E)\big),
\end{align*}
where the isomorphism $(\star)$ is an instance of Corollary \ref{cor:adjunct-oneside}.
It follows that $B_{\eta_{\op P^\bot}}$ is full and faithful.

To see that $B_{\eta_{\op P^\bot}}$ is {essentially} surjective on almost-{cofree} coalgebras, {observe first that if $(C,d)$ is almost-cofree, then there is a graded module $X$ such that $C\cong \Gamma_{\op P^\bot}(X)$.  Moreover,  the} differential $d$ is determined by a map of degree $-1$,
{$d:\Gamma_{\op P^\bot}(X) \rightarrow X$.  The composite
$$\Gamma_{s^{-1}\op P^\bot}(X)\cong s^{-1}\Gamma_{\op P^\bot}(X)\xrightarrow {ds} X$$
extends in the obvious way to an action of $\Omega \op{P}^{\bot}$ on $X$, endowing $X$ with the structure of a $\mathbf{\Omega}\op P^\bot$-algebra.}
\end{proof}

 \subsection{The {nonKoszul} case}

{Let $\op P$ be any operad. The standard construction on the identity map} $id:\opbar \op P\to \opbar\op P$ is simply the two-sided bar construction,
 $$K(\mathrm{id})=\opbar(\op{P},\op{P},\op{P}).$$
{Moreover, the} bar construction on a left $\op{P}$-module $\op{M}$ is the bar construction with coefficients,
 $$B_{\mathrm{id}}\op{M} = \opbar(\op{P},\op{P},\op{M}).$$
The category of $\op{P}$-algebras and strongly homotopy $\op{P}$-morphisms is thus the full subcategory of ${}_{\opbar\op{P}}\cat{Comod}$ spanned by objects of the form $\opbar\big(\op{P},\op{P},z(A)\big)$, which is isomorphic to the Kleisli category for the comonad with underlying functor $\opbar(\op{P},\op{P},-)$.  Note that this {description} reduces to the classical case, of left $A$-modules and strongly homotopy $A$-linear maps.

{To work with} $\op{P}_{\infty}$-algebras, we consider algebras over the resolution $\mathbf{\Omega} \opbar\op{P}$ of $\op{P}$.  We take as our classifying morphism the unit of the adjunction, $\eta:\opbar\op{P} \rightarrow \opbar\mathbf{\Omega} \opbar\op{P}$.  In this case, the standard construction is
$$K(\eta) = \opbar(\mathbf{\Omega} \opbar\op{P};\op{P};\mathbf{\Omega} \opbar \op{P}).$$

 In both cases, the augmentation is  a surjective quasi-isomorphism, and so the resulting Kleisli categories have homotopy-invariant morphism sets, by Proposition~\ref{prop:htyp-inv}.

 \section{Strong homotopy with parameters}\label{sec:delooping}

In this section we study a parametrized version of strong homotopy for morphisms of (co)associative (co)algebras.  We begin by explaining how to parametrize strong homotopy, in terms of the usual cobar/bar adjunction.  We then give a Kleisli-type operadic description of parametrized strong homotopy, which we {apply} to prove a useful existence theorem.

\subsection{Introducing parameters}

Let $F:\cat C\adjunct {}{} \cat D:U$ be a pair of adjoint functors.   Let $\cat {Monad}_{\cat C}$ and $\cat {Comonad}_{\cat D}$ denote the categories of monads on $\cat C$ and of comonads on $\cat D$, respectively. It is well known and easy to show that $UF:\cat C\to \cat C$ is a monad and $FU:\cat D\to \cat D$ is a comonad. More generally, there are functors
$$\cat {Monad}_{\cat D}\to \cat {Monad}_{\cat D}:T\mapsto UTF$$
and
$$\cat {Comonad}_{\cat C}\to \cat{Comonad}_{\cat C}:K\mapsto FKU.$$
We view these functors as providing parametrized families of monads and comonads arising from the $(F,U)$-adjunction.

We now apply this point of view to the usual cobar/bar adjunction for  the associative operad $\op{A}$ (Example~\ref{ex:assoc}), which is a slight variation on the $(\Om_{\op A^\bot}, B_{\op A})$-adjunction of Propositions \ref{prop:bar-cobar}.   Recall from Example \ref{ex:susp-op} the symmetric sequence $\op S$, which gives rise to the \emph{operadic suspension functor}
$$\susp=\op S\otimes-:\cat {dgM}^\Sigma \to \cat {dgM}^\Sigma,$$ which is clearly invertible.  It is easy to check the following useful properties of operadic suspension.

\begin{lem}\label{lem:susp-comon} Operadic suspension is a strongly monoidal functor, i.e., for all symmetric sequences $\op X$ and $\op Y$,
$$\susp(\op X\circ \op Y)\cong \susp \op X\circ \susp \op Y.$$
In particular, operadic suspension  induces endofunctors
$$\susp: \cat {Op}_{\cat {dgM}} \to \cat {Op}_{\cat {dgM}} \quad\text{and}\quad\susp: \cat{CoOp}_{\cat {dgM}} \to \cat {Comon}_{\cat {dgM}}.$$
Moreover, if $\op P$ is an operad and $\op Q$ is a cooperad, then operadic suspension induces a functor
$$\susp:\cat {Mod}_{\op P}\to \cat {Mod}_{\susp\op P}\quad\text{and}\quad \susp: \cat {Comod}_{\op Q}\to \cat {Comod}_{\susp\op Q}.$$
\end{lem}

Recall that $\op A^\bot=\susp \op A^\sharp$.

\begin{defn}\label{defn:usual-barcobar} The \emph{classical cobar construction} $\Om : \coalg {\op A^\sharp} \to \alg A$ is the composite functor
$$\coalg {\op A^\sharp} \xrightarrow\susp \coalg {\op A^\bot} \xrightarrow {\Om_{\op A^\bot}} \alg A,$$
and the \emph{classical bar construction} $B:\alg A \to \coalg {\op A^\sharp}$ is the composite functor
$$\alg A \xrightarrow {B_{\op A}} \coalg {\op A^\bot} \xrightarrow{\susp^{-1}} \coalg {\op A^\sharp}.$$
\end{defn}

It is an easy exercise to show that the classical cobar and bar constructions are indeed the usual, well-known reduced cobar and bar constructions.

The parameters we consider are constructed from symmetric sequences in the following way.
Recall the Schur functor $T:\cat{dgM}^{\Sigma} \rightarrow \End(\cat{dgM})$ from Section~\ref{sssec:alg} {and the related functor $\Gamma :\cat{dgM}^{\Sigma} \rightarrow \End(\cat{dgM}_{+})$ from Section~\ref{sssec:coalg} .}   If $\op P$ is a operad in the category $\alg A$ of associative chain algebras, i.e., an operad in $\cat{dgM}$ with a compatible level monoid structure, then $T_{\op P}$ is the {endofunctor underlying a monad $\mathsf T_{\op P}$ on $\alg{A}$.}  Similarly, if $\op{Q}$ is a cooperad in the category $\coalg {\op A^\sharp}$ of coassociative chain coalgebras, then $\Gamma_{\op Q}$ is the {endofunctor underlying  a comonad $\mathsf \Gamma_{\op Q}$ on $\coalg {\op A^\sharp}$.}
Consequently, there are functors
$$\cat{CoOp}_{\coalg {\op A^\sharp}}\to \cat {Comonad}_{\alg A}:\op Q\mapsto \Omega \Gamma_{\op Q} B$$
and
$$\cat {Op}_{\alg A}\to \cat {Monad}_{\coalg{\op A^\sharp}}:\op P \mapsto B T_{\op P} \Omega,$$
giving us families of (co)monads parametrized by (co)operads.

\begin{defn}\label{defn:parameter} Let $\op Q$ be a cooperad in $\coalg {\op A^\sharp}$, and let $A,B \in \alg A$.  A \emph{$\op Q$-parametrized strongly homotopy morphism} from $A$ to $B$ is a morphism $$\Omega \Gamma_{\op Q} B (A)\to B$$ in $\alg A$.  These are the morphisms in  the coKleisli category
{$$\alg A_{sh,\op Q}:={}_{\Omega \Gamma_{\op Q} B} \alg A.$$}

Dually, let $\op P$ be an operad in $\alg {A}$, and let $C,D\in \coalg {\op A^\sharp}$. A \emph{$\op P$-parametrized strongly homotopy morphism} from $C$ to $D$ is a morphism  $$C\to B T_{\op P} \Om (D) $$ in $\coalg {\op A^\sharp}$.  These are the morphisms in  the Kleisli category
{$$\coalg {\op A^\sharp}_{sh,\op P}:= \coalg {\op A^\sharp}_{B T_{\op P} \Om }.$$}
\end{defn}

 \begin{ex}  We obtain the usual strongly homotopy morphisms of algebras (respectively, of coalgebras) if we set $\op Q$ (respectively, $\op P$) equal to $\op J$.
 \end{ex}

 The parametrized categories defined above admit natural monoidal structures.  The \emph{level tensor product} of symmetric sequences $\op{X}$ and $\op{Y}$ is defined by $(\op{X} \otimes \op{Y})(n) = \op{X}(n) \otimes \op{Y}(n)$.  The constant symmetric sequence $\op{C}$, {where} $\op{C}(n) = R$ with the trivial $\Sigma_{n}$-action for all $n$, is the neutral object for this product.  A (co)monoid with respect to the level tensor product is called, unsurprisingly, a \emph{level (co)monoid}.

\begin{prop} \begin{enumerate}
\item If $\op Q$ is a cooperad in $\coalg {\op A^\sharp}$, then $\alg A_{sh,\op Q}$ has a natural monoidal structure.
\item If $\op P$ is an operad in $\alg {A}$, then $\coalg {\op A^\sharp}_{sh,\op P}$ has a natural monoidal structure.
\end{enumerate}
\end{prop}

\begin{proof} Recall that both $\Om$ and $B$ are monoidal and op-monoidal, i.e., there are natural transformations of functors in $\cat {dgM}$
$$\Om (-)\otimes \Om (-) \Rightarrow \Om (-\otimes -)\Rightarrow \Om (-)\otimes \Om (-)$$
and
$$B (-)\otimes B (-) \Rightarrow B (-\otimes -)\Rightarrow B (-)\otimes B (-)$$
that are appropriately associative and unital.

To prove (1), note that if $(\op X,\Delta)$ is any level comonoid, then $\widehat {\op X}$ is op-monoidal, i.e., there is an appropriately associative and unital natural transformation
$$\Gamma_{\op X}(-\otimes -)\Rightarrow \Gamma_{\op X}(-)\otimes \Gamma_{\op X}(-),$$
given by summing up the natural maps
{\small{{$$\op X(n)\underset {\Sigma_{n}}\otimes (A\otimes B)^{\otimes n}\xrightarrow {\Delta \otimes Id} \big(\op X(n)\otimes \op X(n)\big)\underset {\Sigma_{n}}\otimes (A\otimes B)^{\otimes n}\cong( \op X(n)\underset {\Sigma_{n}}\otimes A^{\otimes n}) \otimes (\op X(n)\underset {\Sigma_{n}}\otimes B^{\otimes n}).$$}}}
It follows that $\Om \Gamma_{\op{Q}}B$ is op-monoidal, since $\op Q$ is a level comonoid.

There is therefore a monoidal structure on $\alg A_{sh,\op Q}$, which is the usual monoidal product on objects. If $f:\Om \Gamma_{Q}B A\to E$ and $f':\Om \Gamma_{Q}B A'\to E'$ represent elements of $\alg A_{sh,\op Q}(A,E)$ and $\alg A_{sh,\op Q}(A',E')$, then the monoidal product of $f$ and $f'$ is equal to the composite
$$\Om \Gamma_{Q}B (A\otimes A')\to \Om \Gamma_{Q}B (A)\otimes \Om \Gamma_{Q}B (A') \xrightarrow {f \otimes f'}E\otimes E'.$$

The proof of (2) is strictly dual and left to the reader.
\end{proof}

\subsection{Diffraction and parametrized strong homotopy}

In this section we provide an operadic description of the categories $\alg A_{sh,\op Q}$ and $\coalg {\op A^\sharp}_{sh,\op P}$, analogous to Theorem \ref{thm:kleisli}.   This description is given in terms of the \emph{diffraction} of $\op Q$ (respectively,  $\op P$), denoted $\Phi(\op Q)$ (respectively, $\Psi(\op P)$), which is an $\op A$-co-ring (respectively, a $\op A^\sharp$-ring).

Our choice of terminology is motivated by the isomorphisms  established below in Theorem \ref{thm:kleisli-x}, which imply the existence of the following bijective correspondence. Let $A$
and $E$ be associative chain algebras.  If $\op Q$ is a cooperad in the category of chain coalgebras, then the set of morphisms of $\op A^\sharp $-coalgebras from $\Gamma_{\op Q} (B A)$ to $B E$ is in bijective correspondence with the set of morphisms of left $\op A$-modules from $\Phi (\op Q)\underset {\op A}\circ z(A)$ to $z(E)$.  In other words, a ``$\op Q$-parametrized''  map on the
bar constructions can be ``diffracted'' into its component pieces on the underlying algebras.  Adding
up the component pieces, we obtain a $\diffract(\op Q)$-parametrization of a morphism between the
underlying algebras.

One advantage to working with $\diffract(\op Q)$ is that  it lends itself well to
existence proofs of $\op Q$-parametrized strongly homotopy morphisms of algebras by acyclic models methods.
We formulate one such existence result and its dual in Theorems \ref{thm:x-exist-bar} and \ref {thm:x-exist-cobar}.

We  need below a few fundamental results about $\op A$-bimodules and $\op A^\sharp$-comodules. We use the following notation throughout the rest of this section.

\begin{notn}\label{notn:alpha}
For all $n\geq 1$, let
$$\delta_{n}=1\cdot Id_{\{1,...,n\}}\in \op A(n)=R[\Sigma _{n}],$$ and let $\delta_{n}^\sharp$ denote the dual element of $\op A^\sharp(n)$.
Let
$$\alpha_{n}=s^{1-n}\delta_{n}\in \susp^{-1}\op A(n)\quad \text{and}\quad \alpha^\sharp_{n}=s^{n-1}\delta_{n}^\sharp\in \op A^\perp(n). $$  For all $\vec n=(n_{1},...,n_{m})\in \N ^m
$, let
$$\delta _{\vec n}=\delta _{n_{1}}\otimes \cdots \otimes \delta _{n_{m}}\in \op A [\vec n],$$
$$\delta^\sharp _{\vec n}=\delta^\sharp _{n_{1}}\otimes \cdots \otimes \delta^\sharp _{n_{m}}\in \op A^\sharp [\vec n],$$
$$\alpha_{\vec n}=\alpha_{n_{1}}\otimes \cdots \otimes a_{n_{m}}\in \susp^{-1}\op A[\vec n]$$
and
$$\alpha^\sharp_{\vec n}=\alpha^\sharp_{n_{1}}\otimes \cdots \otimes \alpha^\sharp_{n_{m}}\in \op A^\perp[\vec n].$$
\end{notn}

\begin{rmk}\label{rmk:levelmon}  Let $\cat {Mon}_{\otimes}$ denote the category of level monoids in $\cat {dgM}^\Sigma$.  The functor $-\circ \op A: \cat {dgM}^\Sigma \to \cat {Mod}_{{\op A}}$ restricts and corestricts to a functor
$$-\circ \op A: \cat {Mon}_{\otimes} \to \bimod AA.$$
Because $\op X\circ \op A$ is a free right $\op A$-module, a
 left $\op A$-action
$$\lambda _{\op X}:\op A \circ \op X\circ \op A \to \op X\circ \op A$$
that is a morphism of right $\op A$-modules is determined by a morphism of symmetric sequences
$$\lambda_{\op X}^\flat:\op A\circ \op X\to \op X\circ \op A.$$
Moreover, since $\op A$ is generated by $\delta_{2}$, it suffices to specify
$$\lambda _{\op X}^\flat\big(\delta_{2}\otimes (x\otimes x')\big)\in \op X\circ \op A$$
for all $x\in \op X(m)$, $x'\in \op X(m')$, and $m,m'\geq0$
to define a left $\op A$-action.
We choose to define $\lambda_{\op X}$ by setting
$$\lambda _{\op X}^\flat\big(\delta_{2}\otimes (x\otimes x')\big)=\begin{cases} xx'\otimes \delta_{2}^{\otimes m}&: m=m'\\ 0&:m\not= m',\end{cases}$$
for all $x\in \op X(m)$, $x' \in \op X'(m')$ and $m,m'\geq 0$.

Similarly,  the functor $-\circ \op A^\sharp: \cat {dgM}^\Sigma \to \cat {Comod}_{A^\sharp}$ restricts and corestricts to a functor
$$-\circ \op A^\sharp: \cat {Comon}_{\otimes} \to \bicomod {\op A^\sharp}{\op A^\sharp},$$
where the left $\op A^\sharp$-coaction on $\op X\circ \op A^\sharp$ is expressed in terms of the level  comultiplication on $\op X$.
\end{rmk}

We can now  define the diffracting functor in terms of operadic suspension and twisting structures and then study its properties.

Recall the right twisting structure $( \mathbf B, \mathbf E, \widetilde {\mathbf E}, \delta, \mu)$ on the category of symmetric sequences of chain complexes (section \ref{subsec:twist-chcx}) and the definition of the adjoint functors
$$(g,g')_{*}:{}_{C}\cat {Comod}_{C'}\to \bimod A{A'}:(g,g')^*$$
induced by classifying morphisms $g:C\to B A$ and $g':C'\to B A'$ (Theorem \ref{thm:key-adjunction}), for any right twisting structure $(B , E, \widetilde E, \delta, \mu)$.  Let
$$\kappa_{\op A}:\op A^\perp \to \mathbf B\op A$$
denote the canonical classifying morphism from Section~\ref{ssec:weight}.

\begin{defn} The \emph{monoid diffracting functor}
$$\Psi:\cat{Mon}_{\otimes}\to \bicomod {\op A^\sharp}{\op A^\sharp}$$ is equal to the composite
$$\cat{Mon}_{\otimes}\xrightarrow{-\circ {\op A}}\bimod AA\xrightarrow {(\kappa_{\op A},\kappa_{\op A})^*}\bicomod{\op A^\perp}{\op A^\perp}\xrightarrow {\susp^{-1}}\bicomod{\op A^\sharp}{\op A^\sharp}.$$
The \emph{comonoid diffracting functor}
$$\Phi:\cat{Comon}_{\otimes}\to \bimod AA$$ is equal to the composite
$$\cat{Comon}_{\otimes}\xrightarrow{-\circ {\op A}^{\sharp}}{}_{\op A^\sharp}\cat {Comod}_{{\op A}^\sharp}\xrightarrow {\susp}{}_{\op A^\perp}\cat {Comod}_{{\op A}^\perp}\xrightarrow {(\kappa_{\op A}, \kappa_{\op A})_{*}} \bimod AA.$$
\end{defn}

\begin{rmk} For every level monoid $(\op X, \mu )$,
$$\Psi(\op X)=\big(\op A^\sharp\circ (\susp^{-1}\op X \circ \susp^{-1}\op A)\circ \op A^\sharp, d_{\Psi}\big),$$
where $d_{\Psi}$ is specified  as follows.  If $x'\in \op X(m')$ and $x''\in \op X(m'')$, then for all $n',n''\geq 1$ and $\vec n'\in I_{n',m'}$, $\vec n''\in I_{n'',m''}$, the perturbation of the internal differential that gives rise to $d_{\Psi}\Big(\delta_{2}^\sharp\otimes \big((s^{1-m'}x' \otimes \alpha_{\vec n'}) \otimes (s^{1-m''}x'' \otimes \alpha_{\vec n''})\big)\Big)$ is
$$\begin{cases}\pm s^{1-m}(x'x'')\otimes \alpha_{\vec n'+\vec n''}&: m'=m''\\0&: m'\not=m'',\end{cases}$$
{where the sign is given by the Koszul rule.}

On the other hand, for all $x\in \op X(m)$ and $\vec n\in I_{n,m}$, the perturbation of the internal differential that determines $d_{\Psi}\big( (s^{1-m}x\otimes \alpha_{\vec n})\otimes (\delta^\sharp_{1})^{\otimes k-1}\otimes \delta^\sharp _{2}\otimes (\delta _{1}^\sharp)^{n-k-1}\big)$ is
$\pm s^{1-m}x\otimes \alpha _{\del^k\vec n}$,
where $\del_{k}(n_{1},...,n_{m})=(n_{1},...,n_{k}+1,...,n_{m})$, and the sign is determined by the Koszul rule.
It suffices to specify these values of the differential, since the underlying bicomodule is cofree and since $\op A^\sharp$ is cogenerated by $\delta_{2}^\sharp$.
\end{rmk}

\begin{rmk}\label{rmk:phi-diffl} For every level comonoid $(\op X,\Delta)$,
$$\Phi (\op X)=\big(\op A\circ (\susp \op X\circ \op A^\perp)\circ \op A, d_{\Phi}),$$
where $d_{\Phi}$ is specified as follows.  Let $x\in \op X(m)$, $\vec n\in I_{m,n}$ and $\alpha^\sharp_{\vec n}\in \op A^{\perp}[\vec n]$ (cf., Conventions \ref{notn:alpha}).   Write $\Delta (x)=x_{i}\otimes x^{i}$ (using the Einstein summation convention).

Then
\begin{align*}
d_{\Phi}(s^{m-1}x\otimes \alpha^\sharp_{\vec n})=&\pm s^{m-1}dx \otimes \alpha^\sharp_{\vec n}\\
&+\delta_{2}\otimes \sum _{\vec n'+\vec n''=\vec n}\pm (s^{m-1}x_{i}\otimes \alpha^\sharp _{\vec n'})\otimes (s^{m-1}x^{i}\otimes \alpha^\sharp_{\vec n''})\\
&+s^{m-1}x \otimes \sum _{0<k<n}\pm \alpha^\sharp_{\vec {n}^k}\otimes (\delta_{1}^{\otimes k-1}\otimes \delta _{2}\otimes \delta _{1}^{{\otimes} n-k-1}),
\end{align*}
where $\vec {n}^k=(n_{1},....,n_{j}-1,...,n_{m})$ if $\sum _{i=1}^{j-1} n_{i}<k \leq \sum _{i=1}^{j} n_{i}${, and the signs are again determined by the Koszul rule.}
\end{rmk}

In order to construct Kleisli categories associated to $\Psi(\op X)$ for a level monoid $\op X$, it must be an $\op A^\sharp$-ring, which is a consequence of the next proposition.  Let $\cat{Ring}_{\op A}$ denote the category of $\op A$-rings.

\begin{prop}\label{prop:a-ring} The functor $-\circ \op A:\cat{Mon}_{\otimes}\to \bimod AA$ restricts and corestricts to a functor
$$-\circ \op A: \cat {Op}_{\alg A}\to \cat {Ring}_{\op A}.$$
\end{prop}

\begin{proof} Let $\op P$ be an operad in the category of associative chain algebras, with multiplication map $\gamma:\op P\circ \op P\to \op P$ and unit map $\eta_{\op P}:\op J\to \op P$.  In particular, $\gamma$ is a morphism of level monoids and therefore induces a morphism of $\op A$-bimodules
$$(\op P\circ \op A)\underset {\op A}\circ (\op P\circ \op A)\cong \op P\circ \op P\circ \op A \xrightarrow {\gamma\circ \op A} \op P\circ \op A,$$
which is the desired multiplication map  on $\op P\circ \op A$.  The unit  of $\op P\circ \op A$ is just $\eta_{\op P}\circ \op A$.
\end{proof}

Proposition \ref{prop:a-ring} and Proposition \ref{prop:param-ring}(2) together imply the following result.

\begin{cor}  If $\op P$ is an operad in the category of associative chain algebras, then $\Psi(\op P)$ is naturally an $\op A^\sharp$-ring.
\end{cor}

As usual, a dual version of Proposition \ref{prop:a-ring} holds as well; we leave its strictly dual proof to the reader. Let $\cat{CoRing}_{\op A^\sharp}$ denote the category of $\op A^\sharp$-co-rings.

\begin{prop}\label{prop:asharp-coring} The functor $-\circ \op A^\sharp:\cat{Comon}_{\otimes}\to \bicomod {\op A^\sharp}{\op A^\sharp}$ restricts and corestricts to a functor
$$-\circ \op A^\sharp: \cat {CoOp}_{\coalg {\op A^\sharp}}\to \cat {CoRing}_{\op A^\sharp}.$$
\end{prop}

Proposition \ref{prop:asharp-coring} and Proposition \ref{prop:param-ring}(1) together imply the following result.

\begin{cor} If $\op Q$ is a cooperad in the category of coassociative chain coalgebras, then $\Phi(\op Q)$ is naturally an $\op A$-co-ring.
\end{cor}

\begin{rmk}  Any simplicial or topological cooperad gives rise to a cooperad in the category of coassociative chain coalgebras, upon application of the normalized chains functor.  Moreover, since $\op A^\sharp (n)$ is the dual Hopf algebra to $R[\Sigma_{n}]$ for all $n$, $\op A^\sharp$ is itself a cooperad in the category of coassociative chain coalgebras.

Dually, if we apply the normalized {cochains functor to a cooperad in the category of simplicial sets or topological spaces}, then we obtain an operad in the category of associative chain algebras, of which the associative operad itself is another important example.
\end{rmk}

 \begin{rmk} From the definitions above,  one can deduce a formula for the comultiplication $\psi_{\op Q}:\Phi (\op Q)\to \Phi (\op Q) \acirc \Phi (\op Q)$, where $\op Q$ is any cooperad in the category of coassociative chain coalgebras.  Observe that since the $\op A$-bimodule underlying $\Phi(\op Q)$ is free, $\psi_{\op Q}$ is determined by its image on $\susp \op Q\circ \op A^\perp$.

 Let $x\in \op Q(m)$, and let $\vec n \in I_{m,n}$.  If $m$ does not divide $n$, then $\psi_{\op Q}(s^{m-1}x\otimes \alpha^\sharp _{\vec n})=0$.

 If there exists a natural number $\ell$ such that $n=m\ell$ and $\vec n\in I_{n,m}$, then let
 $$K_{\vec n, \ell}=I_{\ell, n_{1}}\times \cdots \times I_{\ell, n_{m}}.$$
 and for all $(\vec n_{1},...,\vec n_{m})\in K_{\vec n, \ell}$,  with $\vec n_{i}=(n_{i,1},...,n_{i,\ell})$,
 let
 $$\vec n^j=(n_{1,j},...,n_{m,j}),$$
 for all $1\leq j\leq \ell$.  In terms of this notation we have for $x\in \op X(m)$ and $\vec n\in I_{n,m}$ that
 {\small \begin{align*}
 \psi_{\op Q}(s^{m-1}x\otimes &\alpha^\sharp _{\vec n})\\
 =&\sum _{i, K_{\vec n, \ell},\vec\ell\in I_{m,\ell}} \pm s^{m-1}x_{i,0}\otimes \alpha^\sharp_{\vec \ell}\otimes \big((s^{m-1}x_{i,1}\otimes\alpha^\sharp _{\vec n^1})\otimes \cdots \otimes (s^{m-1}x_{i,\ell}\otimes \alpha^\sharp_{\vec n^\ell})\big),
 \end{align*}}
 where the signs are determined by the Koszul rule and
 $$\Delta ^{(\ell+1)}(x)=\sum _{i}x_{i,0}\otimes \cdots \otimes x_{i,\ell}.$$

 We leave it as a (rather technical) exercise for the reader to obtain analogous formulas for the multiplication $\Psi(\op P)\underset{\op A^\sharp}\square \Psi(\op P) \to \Psi(\op P)$.
 \end{rmk}

 \subsection{Kleisli categories and parametrized strong homotopy}

In this section we give an operadic, Kleisli category description of parametrized strong homotopy of (co)associative (co)algebras{, in the spirit of Theorem \ref{thm:kleisli}}.

 \begin{defn}  Let $\op Q$ be a cooperad in the category of coassociative chain coalgebras.  Let {$\mathsf \Phi_{\op Q}$ denote the comonad on ${}_{\op A}\cat {Mod}$ with underlying endofunctor
 $$\Phi(\op Q)\underset {\op A}\circ-:{}_{\op A}\cat {Mod} \to {}_{\op A}\cat {Mod}.$$
 The induced monad on $\alg A$ is also denoted $\mathsf \Phi (\op Q)$}

 Let $\op P$ be a operad in the category of associative chain algebras.  Let {$\mathsf\Psi_{\op P}$  denote the monad on ${}_{\op A^\sharp}\cat {Comod}$ with underlying endofunctor
 $$\Psi(\op P)\underset {\op A^\sharp}\square-:{}_{\op A^\sharp}\cat {Comod} \to {}_{\op A^\sharp}\cat {Comod}.$$
 The induced comonad on $\coalg {\op A^\sharp}$ is also denoted $\mathsf \Psi(\op P)$.}
 \end{defn}

 \begin{thm}\label{thm:kleisli-x} Let $\op Q$ be a cooperad in the category of coassociative chain coalgebras, and let $\op P$ be an operad in the category of associative chain algebras. {There are isomorphisms of categories
 $$\Theta_{\op A,\op Q}:\alg{\op A}_{sh,\op Q} \xrightarrow\cong {}_{\mathsf \Phi(\op Q)}\alg A$$
 and
 $$\Upsilon_{\op A,\op P}:\coalg{\op A^\sharp}_{sh,\op P} \xrightarrow\cong\coalg {\op A^\sharp}_{\mathsf\Psi_{\op P}},$$
 which are the identity on objects.}
 \end{thm}

 \begin{proof} The proof of this theorem strongly resembles that of Theorem \ref{thm:kleisli}. The definition of $\Theta_{\op A,\op Q}$ on morphisms follows from the sequence of natural isomorphisms below, where $A$ and $E$ are $\op A$-algebras.
\begin{align*}
{}_{\mathsf \Phi(\op Q)}&\alg A(A,E)\\
&={}_{\op A}\cat {Mod}\big(\Phi (\op Q)\underset {\op A}\circ z(A), z(E)\big)\\
&={}_{\op A}\cat {Mod}\big(\op A\circ_{\kappa_{\op A}}\susp (\op Q\circ \op A^\sharp)\circ _{\kappa_{\op A}}\op A\underset {\op A}\circ z(A), z(B)\big)\\
&\cong{}_{\op A^\perp}\cat {Comod}\big(\susp \op Q\circ \op A^\perp \circ_{\kappa_{\op A}} z(A), \op A^\perp\circ _{\kappa_{\op A}}z(E)\big)\\
&\cong{}_{\op A^\perp}\cat {Comod}\big(\susp \op Q\circ z(B_{\op A}A), z(B_{\op A}E)\big)\\
&\cong{}_{\op A^\sharp}\cat {Comod}\big(\op Q\circ z(B A), z(B E)\big) \\
&\cong \coalg{\op A^\sharp}(\Gamma_{\op Q}B A,B E)\\
&\cong \alg A\big(\Om \Gamma_{\op Q}B (A), E\big )\\
&=\alg{\op A}_{sh,\op Q}(A,E).
\end{align*}

The proof of the dual case is similar and left to the reader.
 \end{proof}

 As a consequence of Theorem \ref{thm:kleisli-x}, it is relatively easy to prove the following existence theorems{, which is expressed in terms of acyclic models. We recall the foundations of
this method before stating the theorems.

Let $\cat D$ be a  category, and let $\mathfrak M$ be a set of objects in $\cat D$.    A functor $X: \cat D\to
\cat {dgM}$ is \emph{free} with respect to $\mathfrak M$ if there is a set $\{x_{\mathfrak m}\in X(\mathfrak
m)\mid \mathfrak m\in \mathfrak M\}$ such that $\{X(f)(x_{\mathfrak m})\mid f\in \cat D(\mathfrak m,D),
\mathfrak m\in \mathfrak M\}$ is a $R$-basis of $X(D)$ for all objects $D$ in $\cat D$. The
functor $X$ is \emph{acyclic}  with respect to $\mathfrak M$ if $X(\mathfrak m)$ is acyclic for all $
\mathfrak m\in \mathfrak M$.   More generally, if $\cat C$ is a category with a forgetful functor $U$ to $
\cat {dgM}$ and $X:\cat D\to \cat C$ is a functor, we say that $X$ is \emph {free}, respectively \emph
{acyclic},  with respect to $\mathfrak M$ if $UX$ is.}

 \begin{thm} \label{thm:x-exist-bar} Let $\op Q$ be a cooperad in $\coalg {\op A^\sharp}$ with $\op Q(0)=0$ and $\op Q(1)=R$. Let $X,Y: \cat D \to \alg {\op A}$ be functors, where $\cat D$ is a category admitting a set of models $\mathfrak M$ with respect to which $X$ is free and such that $Y(\mathfrak m)$ is acyclic for every $\mathfrak m$ that is a coproduct of elements of $\mathfrak M$.  Let $U: \alg{\op P} \to \cat {dgM}$ be the forgetful functor.

 If  $\tau: UX\Rightarrow UY$ is a natural transformation, then there is a natural transformation $\hat \tau: \Om \widehat {\op Q}B\circ X\Rightarrow Y$ of functors from $\cat D$ into $\alg A$ extending $\tau$, i.e., for all $D\in \ob \cat D$, the following composite is equal to $\tau_{D}$.
$$ X(D)\hookrightarrow \Om \widehat {\op Q}B X(D)\xrightarrow {\hat \tau _{D}}Y(D)$$
In other words, for each $D\in \ob \cat D$, the natural chain map $\tau_{D}:X(D)\to Y(D)$ admits a natural,
$\op Q$-parametrized, strongly homotopy multiplicative structure. \end{thm}

 There is also a coalgebra version.

  \begin{thm}\label{thm:x-exist-cobar}  Let $\op P$ be an operad in $\alg A$ with $\op P(0)=0$ and $\op P(1)=R$. Let $X,Y: \cat D \to \coalg {\op A^\sharp}$ be functors, where $\cat D$ is a category admitting a set of models $\mathfrak M$ with respect to which $X$ is free and $Y$ is acyclic.  Let $U: \coalg{\op A^\sharp} \to \cat {dgM}$ be the forgetful functor.

If  $\tau: UX\Rightarrow UY$ is a natural transformation, then there is a natural transformation $\hat \tau: X \Rightarrow B \widehat {\op P}\Om\circ Y$ of functors from $\cat D$ into $\coalg{\op A^\sharp}$ lifting $\tau$, i.e., for all $D\in \ob \cat D$, the following composite is equal to $\tau_{D}$.
$$X(D)\xrightarrow {\hat \tau _{D}}B \widehat {\op P}\Om Y(D)\xrightarrow{\text{proj.}} Y(D)$$
In other words, for each $D\in \ob \cat D$, the natural chain map $\tau_{D}:X(D)\to Y(D)$ admits a natural,
$\op P$-parametrized, strongly homotopy comultiplicative structure.
 \end{thm}

 The proofs of these theorems are discussed in Appendix \ref{app2}.

\appendix

\section{The Alexander-Whitney co-ring, by P.-E. Parent}\label{app1}

The categories $\cat {DASH}=\alg{\op A}_{sh}$ and $\cat {DCSH}=\coalg{\op A^\sharp}_{sh}$ are particularly important in topology. For example, Gugenheim and Munkholm  showed that the Alexander-Whitney equivalence
\[
    C_{*}(K \times L) \rightarrow C_{*}(K)\otimes C_{*}(L)
\]
of normalized chains on reduced simplicial sets is naturally the linear part of  a
morphism in $\cat {DCSH}$, which implies, as shown in \cite {hpst:canonicalAH}, that the same is true of the usual comultiplication on $C_{*}K$. We devote this section to a careful analysis of the operadic description of these categories.

 Let $\op{F} =K(\op A)=(\op A\circ \op A^\perp\circ \op A, \del _{\op F})=\Phi (\op J)$. Then
\[
    \op{F} = \{ \op{F}(m) \mid m \in \N \}
\]
with generators
\[
    \{ f_{m} = s^{m-1}u_{0} \in \op{F}(m)_{m-1} \mid m \in \N \}
\]
satisfying
\[
    \partial f_{m} = \sum_{i=1}^{m-1} \delta \otimes (f_{i} \otimes
    f_{m-i}) + \sum_{i=0}^{m-2} f_{m-1} \otimes (1^{\otimes i} \otimes
    \delta \otimes 1^{\otimes (m-2-i)}).
\]

By Proposition~\ref{prop:quad-co-ring}, ${\op{F}}$ is an $\op{A}$-co-ring, which we call the \emph{Alexander-Whitney co-ring}. The formula for the
composition comultiplication $\psi_{\op F}$
is particularly simple.  For $n \geq 1$, we have
\[
    \psi_{\op F}( f_{n}) = \sum_{m \leq n} \sum_{\vec{n} \in I_{m,n}}
    f_{m}
    \otimes ( f_{n_{1}} \otimes \cdots \otimes f_{n_{m}} )
\]
where $n_{i} \geq 1$ for all $i$.  In fact $\op F$ is a counital $\op{A}$-co-ring, with counit $$\varepsilon:
\op F\to \op A$$ specified by $\varepsilon(f_{n})=0$ for all $n>1$ and $\varepsilon (f_{1})=1$. Since $\op F$ is the Koszul resolution of $\op A$, the counit is a levelwise quasi-isomorphism.

Moreover,
$(\op{F},\partial_{\op{F}},\Delta_{\op{F}})$ is a level comonoid in the category
$\op{A}$-bimodules. Explicitly,
\begin{eqnarray*}
    \Delta_{\op{F}}(f_{m})
    & = & \sum_{k=1}^{m}\sum_{\vec{\imath} \in I_{k,m}}
    \left(f_{k} \otimes \delta^{(i_{i})}\otimes \cdots \otimes
    \delta^{(i_{k})}\right)
    \otimes
    \left( \delta^{(k)} \otimes f_{i_{1}} \otimes \cdots
    \otimes f_{i_{k}}\right).
\end{eqnarray*}
It is easy to check that $\Delta _{\op F}$ is coassociative.

It follows from Theorem \ref{thm:kleisli} or Theorem \ref{thm:kleisli-x} that there is a natural isomorphism
$$\cat {DASH}(A,A')\cong {}_{\op A}\cat {Mod}\big(\op F\underset {\op A}\circ z(A), z(A')\big)$$
for all associative chain algebras $A$ and $A'$.

On the other hand we can also characterize morphisms in $\cat {DCSH}$ in terms of $\op F$, as described below.  This alternate characterization has already proved useful in, e.g.,  \cite{hpst:canonicalAH}, \cite{hps:cohoch}, \cite{hess-levi:07}, \cite{hess-rognes}, \cite {naito} and \cite{boyle}.

Let $\mathsf F$ denote the comonad on ${}_{\op A}\cat {Mod}_{\op A}$ with underlying endofunctor  $-\underset {\op A}\circ \op F$.

An $\op{A}$-\emph{coalgebra} is a chain complex $C$, equipped with structure morphisms
\[
    \psi_{n} : C \otimes \op{A}(n) \rightarrow C^{\otimes n} \quad (n \geq 1)
\]
that are associative, equivariant and unital with respect to $\op{J} \rightarrow \op{A}$.  This definition does coincide with that of $\op{A}^{\sharp}$-coalgebras, since $\op{A}$ is projective.  We denote by $\coalg{A}$ the category of $\op{A}$-coalgebras and morphisms that commute strictly with the structure morphisms.

The tensor functor
$$\op T:\cat {dgM} \to \cat {dgM}^\Sigma,$$
where $\op T(X)(n)=X^{\otimes n}$, restricts and corestricts to a full and faithful functor
$$\op T: \coalg {\op A} \to {}_{\op A}\cat {Mod}_{\op A}.$$
It follows that $\mathsf F$ induces a comonad, also denoted $\mathsf F$, on $\coalg {\op A}$.

\begin{thm}\label{thm:dcsh}
There is an isomorphism of categories
$$\cat{DCSH} \cong {}_{\mathsf F}\coalg{\op A}.$$
\end{thm}

We will make use of the following lemma.

\begin{lem}\label{lem:right-adjoints}
Let $f : W \otimes Y \rightarrow Z$ and $g : X \rightarrow Y$ be
morphisms in $\cat{dgM}$.  The right adjoint to the composite
\[
    W \otimes X \xrightarrow{1 \otimes g} W \otimes Y
    \xrightarrow{f} Z
\]
is the composite
\[
    W \xrightarrow{\hat{f}} \Hom(Y,Z) \xrightarrow{\Hom(g,Z)}
    \Hom(X,Z),
\]
where $\hat{f}$ is the right adjoint to $f$.
\end{lem}

\begin{proof}
Recall that for all $w \in X$, $\hat{f}_{w} \in \Hom(Y,Z)$ is
defined by $\hat{f}_{w}(y) = f(w \otimes y)$.  We set $g^{\dual} =
\Hom(g,Z)$, so that $g^{\dual}(\phi) = \phi g$ for all $\phi \in
\Hom(Y,Z)$.

For all $x \in X$, we have
\begin{eqnarray*}
    (g^{\dual}\hat{f}_{w})(x)
        & = & g^{\dual}(\hat{f}_{w})(x) \\
        & = & \hat{f}_{w}(g(x)) \\
        & = & f(w \otimes g(x)) \\
        & = & f(1 \otimes g)(w \otimes x) \\
        & = & (f(1\otimes g))^{\wedge}_{w}(x).
\end{eqnarray*}
Therefore, the right adjoint of $f(1 \otimes g)$ is
$g^{\dual}\hat{f}$.
\end{proof}

\begin{proof}[Proof of Theorem~\ref{thm:dcsh}.]
First, we show that there is a functor
\[
    \ind : {}_{\mathsf F}\coalg{\op A} \rightarrow \cat{DCSH}
\]
such that $\ind(\op{T}(C)) = C$ on objects.

A morphism  $\theta$ in ${}_{\mathsf F}\coalg{A}(C,D) = {}_{\op{A}}\cat{Mod}_{\op{A}}(\op{T}(C) \circ_{\op{A}}\op{F},\op{T}(D))$ is determined by a sequence of equivariant morphisms,
\[
    C \otimes \op{A}^{\bot}(n) \rightarrow D^{\otimes n}.
\]
Desuspend $n$ times to obtain the equivariant morphism
\[
    (s^{-1}C) \otimes \op{A}^{\dual}(n) \rightarrow
    (s^{-1}D)^{\otimes n}.
\]
Now take the right adjoint, then pass to orbits.  The end result is a
morphism
\[
    s^{-1}C \rightarrow \op{A}(n) \otimes_{\Sigma_{n}}
    (s^{-1}D)^{\otimes n}
\]
that extends uniquely to a morphism of $\op{A}$-algebras,
\[
    \ind(\theta) : T_{\op{A}}(s^{-1}C) \rightarrow
    T_{\op{A}}(s^{-1}D).
\]
Since the internal differentials from $C$ and $D$ automatically commute with $\ind(\theta)$, we may suppose that they are zero. In particular, the differentials in $\Omega(C)$ and $\Omega(D)$ are then entirely quadratic: $d = d_{2}$.

Let $F = \ind(\theta)$; then $Fd$ is defined by the composition
\begin{multline}\label{eq:Fd}
    s^{-1}C \xrightarrow{d} \op{A}(2) \otimes (s^{-1}C)^{\otimes
    2}  \xrightarrow{1 \otimes F} \op{A}(2) \otimes
    \op{A}[2,n] \otimes (s^{-1}D)^{\otimes n} \\
    \xrightarrow{\gamma \otimes 1} \op{A}(n) \otimes (s^{-1}D)^{\otimes n}.
\end{multline}
Once we identify $\Hom(X,-)$ with $X^{\dual} \otimes -$,
Lemma~\ref{lem:right-adjoints} exhibits \eqref{eq:Fd} as the right
adjoint to the following composition:
\begin{multline}\label{eq:Fd-left}
    s^{-1}C \otimes \op{A}^{\dual}(n) \xrightarrow{1 \otimes
    \gamma^{\dual}} s^{-1}C \otimes \op{A}^{\dual}(2) \otimes
    \op{A}^{\dual}[2,n] \\ \xrightarrow{\tilde{d} \otimes 1}
    (s^{-1}C)^{\otimes 2}
    \otimes \op{A}^{\dual}[2,n] \xrightarrow{\tilde{F}}
    (s^{-1}D)^{\otimes n}
\end{multline}
where $W = s^{-1}C$, $X = \op{A}^{\dual}(n)$, $Y =
\op{A}^{\dual}(2) \otimes \op{A}^{\dual}[2,n]$, $Z =
(s^{-1}D)^{\otimes n}$, $g = \gamma$, $\hat{f} = (1 \otimes
F)d$, and $\tilde{d}$ and $\tilde{F}$ are the left adjoint to the
cobar differential and $F$, respectively. If we suspend
\eqref{eq:Fd-left} $n$ times, we obtain $(-1)^{n}$ \emph{times} the
composition,
\begin{multline}\label{eq:Fd-left-susp}
    C \otimes \op{A}^{\bot}(n) \xrightarrow{1 \otimes
    \gamma^{\bot}} C \otimes \op{A}^{\bot}(2) \otimes
    \op{A}^{\bot}[2,n] \xrightarrow{1 \otimes s^{-1} \otimes 1}
    C \otimes \op{A}(2) \otimes \op{A}^\bot [2,n] \\
    \xrightarrow{\rho \otimes 1}
    C^{\otimes 2}
    \otimes \op{A}^{\bot}[2,n] \xrightarrow{\theta}
    D^{\otimes n}
\end{multline}
where $\rho : C \otimes \op{A}(2) \rightarrow C^{\otimes 2}$ is the
structure map for $C$.  The sign of $(-1)^{n}$ is introduced because
$(1 \otimes s^{-1} \otimes 1)$ has degree $-1$ and $s^{n}$ has degree $n$. We note that \eqref{eq:Fd-left-susp} defines $\theta d_{L}$.

Now, since $\theta$ commutes with differentials, and since we are
assuming that internal differentials vanish, $\theta(d_{L} - d_{R})
= 0$.  Therefore \eqref{eq:Fd-left-susp} coincides with$(-1)^{n}\theta
d_{R}$, namely, with $(-1)^{n}$ times  the composition
\begin{multline}\label{eq:dF-left-susp}
    C \otimes \op{A}^{\bot}(n) \xrightarrow{1 \otimes
    \gamma^{\bot}} C \otimes \op{A}^{\bot}(n-1) \otimes
    \op{A}^{\bot}[n-1,n] \xrightarrow{1 \otimes 1 \otimes s^{-1}}
    C \otimes \op{A}^{\bot}(n-1) \otimes \op{A}[n-1,n] \\
    \xrightarrow{\theta \otimes 1} D^{\otimes n-1} \otimes
    \op{A}[n-1,n] \xrightarrow{\rho} D^{\otimes n}.
\end{multline}
This time, when we desuspend \eqref{eq:dF-left-susp} $n$ times, we
obtain $(-1)^{n}$ times:
\begin{multline}\label{eq:dF-left}
    s^{-1}C \otimes \op{A}^{\dual}(n) \xrightarrow{1 \otimes
    \gamma^{\dual}} s^{-1}C \otimes \op{A}^{\dual}(n-1) \otimes
    \op{A}^{\dual}[n-1,n]  \\
    \xrightarrow{\tilde{F} \otimes 1} (s^{-1}D)^{\otimes n-1} \otimes
    \op{A}^{\dual}[n-1,n] \xrightarrow{\tilde{d}} (s^{-1}D)^{\otimes n}.
\end{multline}
The sign is introduced for essentially the same reason: $1 \otimes 1 \otimes s^{-1}$ has degree $-1$ while $s^{n}$ has degree $n$.

By Lemma~\ref{lem:right-adjoints}, \eqref{eq:dF-left} is left
adjoint to the composite,
\begin{multline}\label{eq:dF}
    s^{-1}C \xrightarrow{F} \op{A}(n-1) \otimes
    (s^{-1}D)^{\otimes n-1} \\
        \xrightarrow{1 \otimes d} \op{A}(n-1) \otimes
        \op{A}[n-1,n] \otimes (s^{-1}D)^{\otimes n}
        \xrightarrow{\gamma \otimes 1} \op{A}(n) \otimes
        (s^{-1}D)^{\otimes n},
\end{multline}
which defines $dF$.
Therefore \eqref{eq:dF} and \eqref{eq:Fd} coincide, and so
$\ind(\theta)$ is a chain map.

Let $\theta : \op{T}(C) \circ_{\op{A}} \op{F} \rightarrow
\op{T}(D)$ and $\varphi : \op{T}(D) \circ_{\op{A}} \op{F}
\rightarrow \op{T}(E)$ be morphisms of $(\op{A},\op{A})$-bimodules.
The morphism $\ind(\varphi) \circ \ind(\theta)$ is the sum of
composites,
\[
    s^{-1}C \rightarrow \op{A}(m) \otimes (s^{-1}D)^{\otimes m}
    \rightarrow \op{A}(m) \otimes \op{A}[m,n] \otimes
    (s^{-1}E)^{\otimes n} \rightarrow \op{A}(n) \otimes
    (s^{-1}E)^{\otimes n}.
\]
The left adjoint of the above composite coincides with the
composite,
\[
    s^{-1}C \otimes \op{A}^{\dual}(n) \rightarrow s^{-1}C \otimes
    \op{A}^{\dual}(m) \otimes \op{A}^{\dual}[m,n] \rightarrow
    (s^{-1}D)^{\otimes m} \otimes \op{A}^{\dual}[m,n] \rightarrow
    (s^{-1}E)^{\otimes n}.
\]
Suspending, we obtain the component in arity $n$ of $\varphi \circ
\theta$.  It follows that $\ind(\varphi \circ \theta) =
\ind(\varphi) \circ \ind(\theta)$.

Next, we show that  $\ind$ has an inverse functor
\[
    \lin :\cat{DCSH} \to {}_{\mathsf F}\coalg A
    \]
that coincides with $\op{T}$ on objects.

Since $\op{A}$ and $\op{A}^{\dual}$ are projective, the natural morphism $\pi :
(\op{A}^{\dual}(n) \otimes V^{\otimes n})^{\Sigma_{n}} \rightarrow
(\op{A}^{\dual}(n) \otimes V^{\otimes n})_{\Sigma_{n}}$ is
invertible for all $n$.

Let $F : \Omega(A) \rightarrow \Omega(C)$ be a
morphism in $\cat{DCSH}$.   We take the adjoint
of the composite
\begin{multline*}
        s^{-1}A \xrightarrow F
        \op{A}(n) \otimes_{\Sigma_{n}} (s^{-1}C)^{\otimes n} \\
        \xrightarrow{\nu}
        \Hom(\op{A}^{\dual}(n),(s^{-1}C)^{\otimes n})_{\Sigma_{n}}
         \xrightarrow{\pi^{-1}}
        \Hom(\op{A}^{\dual}(n),C^{\otimes n})^{\Sigma_{n}}
\end{multline*}
to obtain a $\Sigma_{n}$-equivariant morphism
\[
    s^{-1}A \otimes \op{A}^{\dual}(n) \rightarrow
    C^{\otimes n}.
\]
Suspending $n$ times, we obtain a $\Sigma_{n}$-equivariant morphism
\[
    \varphi_{n} : A \otimes \op{A}^{\bot}(n) \rightarrow C^{\otimes n}.
\]
The sequence $(\varphi_{n})$ defines a morphism of symmetric sequences,
\[
    \op{T}(A) \circ \op{A}^{\bot} \rightarrow \op{T}(C)
\]
that extends to define a morphism of right $\op{A}$-modules,
\[
    \lin(F) : \op{T}(A) \circ_{\op{A}} \op{F} \rightarrow
    \op{T}(C).
\]
The morphism $\lin(F)$ commutes with differentials and respects
composition by symmetric arguments to those for $\ind$.

Clearly $\ind\lin$ and $\lin\ind$ are the identity on objects.
Furthermore, the algorithm that yields $\lin(F)$ is the reverse of
the algorithm for $\ind(\theta)$.  Therefore $\ind\lin = 1$ and
$\lin\ind = 1$.
\end{proof}

Specializing Theorems \ref{thm:x-exist-bar} and \ref{thm:x-exist-cobar} to the case $\op P=\op J=\op Q$, we obtain the following existence results, which have already been applied to great effect in \cite{hps:james}, \cite{hess-levi:07}, \cite{hps:cohoch}, and \cite{hess-rognes}.

\begin{thm}\label{thm:dash-exist}
Let $X,Y: \cat D \to \alg {\op A}$ be functors, where $\cat D$ is a category admitting a set of models $\mathfrak M$ with respect to which $X$ is free and such that $Y(\mathfrak m)$ is acyclic for all objects $\mathfrak m$ that are coproducts of objects in $\mathfrak M$.  Let $U: \alg{\op A} \to \cat {dgM}$ be the forgetful functor.

If $\tau: UX\Rightarrow UY$ is a natural transformation, then
 there is a natural transformation $\hat \tau:  \Omega B\circ X\Rightarrow  Y$ of functors from $\cat D$ into $\alg A$ extending $\tau$, i.e., for all $D\in \ob \cat D$, the following composite is equal to $\tau_{D}$.
$$ X(D)\hookrightarrow \Om B X(D)\xrightarrow {\hat \tau _{D}}Y(D)$$
In other words, for each $D\in \ob \cat D$, the natural chain map $\tau_{D}:X(D)\to Y(D)$ admits a natural
DASH-structure.
\end{thm}

\begin{thm}\label{thm:dcsh-exist}
Let $X,Y: \cat D \to \coalg {\op A^\sharp}$ be functors, where $\cat D$ is a category admitting a set of models $\mathfrak M$ with respect to which $X$ is free and $Y$ is acyclic.  Let $U: \coalg{\op A^\sharp} \to \cat {dgM}$ be the forgetful functor.

 If $\tau: UX\Rightarrow UY$ is a natural transformation, then there is a natural transformation $\hat \tau:  X\Rightarrow B\Om\circ Y$ of functors from $\cat D$ into $\alg{\op P}$ lifting $\tau$, i.e., for all $D\in \ob \cat D$, the following composite is equal to $\tau_{D}$.
$$ X(D)\xrightarrow {\hat \tau _{D}}B\Om Y(D)\xrightarrow{\text{proj.}} Y(D)$$
In other words, for each $D\in \ob \cat D$, the natural chain map $\tau_{D}:X(D)\to Y(D)$ admits a natural
DCSH-structure.
\end{thm}

\begin{rmk} Since the bimodule ${\op{F}}$ is a free
${\op{A}}$-bimodule resolution of ${\op{A}}$, we
may use it to do homological algebra. If $\op{M}$ and $\op{N}$ are
right and left ${\op{A}}$-modules, respectively, then
\[
    \mathrm{Tor}^{{\op{A}}}(\op{M},\op{N})
        := H( \op{M} \underset{{\op{A}}}{\circ}
        {\op{F}} \underset{{\op{A}}}{\circ}
        \op{N}).
\]
In particular, we may read off the isomorphisms
\[
    \mathrm{Tor}^{{\op{A}}}(\op{J},\op{J}) \cong \Ssigma
    \quad \text{and} \quad
    \mathrm{Tor}^{{\op{A}}}(\op{J},{\op{A}}) \cong
    {\op{J}}.
\]
\end{rmk}

\begin{ex}\label{ex:counter}
We are now in a position to construct over $R=\Zmod{2}$ an example
of a chain coalgebra $M$ such that
\begin{itemize}
\item  its cohomology
algebra is realizable, i.e., there is a topological space $X$ such that $H^*(X;\mathbb F_{2})\cong H^*M$
as graded algebras, but
\item  $M$ is not
quasi-isomorphic to the chains on any space.
\end{itemize}
This example,
along with~\cite[Example 3.8]{ndombol-thomas:04}, show that the
concepts of ``\emph{shc} algebras'' and ``algebras with cup-$i$
products'' are independent of one another.

Let $M = \Zmod{2}\{ 1_{0}, u_{2}, x_{3}, y_{3}, z_{3}, v_{4},
w_{6} \}$, where subscript indicates degree.  The only non-zero
differential in $M$ is $\partial(v)=x+y$.  All elements other than
$v$ and $w$ are primitive, while $B{\psi}(v) = u \otimes u$ and
$B{\psi}(w) = x \otimes z + z \otimes y$.  It is readily
verified that $(M,\partial,\psi)$ is a coassociative chain
coalgebra.

Let $W$ be the usual $\Zmod{2}[\Sigma_{2}]$-free resolution of
$\Zmod{2}$.  Specifically, $W_{i}$ is generated by an element
$e_{i}$ with $\partial e_{i} = (1+\tau)e_{i-1}$, where $\tau \in
\Sigma_{2}$ is the transposition.

\begin{prop}
There exists an equivariant morphism
\[
    g : W \otimes M \rightarrow M \otimes M
\]
such that $g ( e_{0} \otimes - ) = \psi$.
\end{prop}

\begin{proof}
We construct the morphism $g$; verification that it is a
chain map is routine and left to the reader.

It suffices to define $g$ on generators.  The only non-zero
values that $g$ takes on generators are $g(e_{1} \otimes
w) = v \otimes z + z \otimes v$, $g(e_{3} \otimes v) = v
\otimes x + y \otimes v$, and $g(e_{|a|} \otimes a) = a
\otimes a$ for $a \in \{ u, v, w, x, y, z \}$.
\end{proof}

By~\cite{may:70}, $g$ defines cup-$i$ products in the $\Zmod{2}$-dual
$M^{\sharp}$, and so $H^{*}(M,\Zmod{2})$ comes equipped with an
action of the mod $2$ Steenrod algebra.  In fact, we have the
following proposition.

\begin{prop}The algebra $H^{*}(M;\Zmod{2})$ admits the structure of an
unstable algebra over the mod $2$ Steenrod algebra, where the only
non-trivial operation is the $\mathrm{Sq}^0$. Moreover, this
algebra is isomorphic to
\[
    H^*(S^2\vee(S^3\times S^3);\Zmod{2})
\]
as unstable algebras over the mod $2$ Steenrod algebra.
\end{prop}

\begin{proof}
An easy exercise in $\Zmod{2}$-linear algebra.
\end{proof}

\begin{prop}\label{prop:example}
The chain coalgebra $M$ is not realizable, i.e., $M$ is not of the
same homotopy type as $C_{*}(X;\Zmod{2})$ for any space $X$.
\end{prop}

\begin{proof}
For the duration of the proof, we suppress the coefficients from
the notation. By~\cite{gugenheim-munkholm:74}, if $X$ is a space,
then the diagonal on $C_{*}(X)$ is strongly
homotopy-comultiplicative, that is, there is a morphism of
symmetric sequences, $\Delta : \op{T}(C_{*}(X)) \circ \op{F}
\rightarrow \op{T}(C_{*}(X) \otimes C_{*}(X))$, such that
$\Delta_{1}$ is the diagonal.  If $C_{*}(X)$ and $M$ are connected
by a sequence of chain coalgebra quasi-isomorphisms, then we may
construct a morphism $\Psi : \op{T}(M) \circ \op{F} \rightarrow
\op{T}(M \otimes M)$ such that $\Psi_{1} = \psi$.  The homotopy
class of such a $\Psi$, compatible with $\Delta$, is unique. We
show that no such $\Psi$ exists.

We attempt to define $\Psi$ on generators $a \otimes f_{k}$, for
$a \in M$ and $k \geq 1$.  Necessarily, $\Psi(a \otimes f_{1}) =
\psi(a)$.  We may define
\[
    \Psi(w \otimes f_{2}) = (1 \otimes v) \otimes (z \otimes 1)
    + (1 \otimes z) \otimes (v \otimes 1)
\]
and $\Psi(a \otimes f_{2}) = 0$ for $a \neq w$.  Any other choice
of morphism $\Psi':M \otimes \op{F}(2) \rightarrow (M \otimes
M)^{\otimes 2}$ is necessarily homotopic to $\Psi$.

Now we try to define $\Psi$ on $M \otimes \op{F}(3)$.  In order to
find a value for $\Psi(w \otimes f_{3})$, we must find an element
that bounds
\[
    (1 \otimes z) \otimes (u \otimes 1) \otimes (u \otimes 1)
    + (1 \otimes u) \otimes (1 \otimes u) \otimes (z \otimes 1),
\]
but no such element exists. Therefore the diagonal on $M$ does not
extend to an $\op{F}$-parametrized morphism.
\end{proof}

\end{ex}

\section{Proof of Theorems \ref{thm:x-exist-bar} and \ref{thm:x-exist-cobar}}\label{app2}

We prove here Theorem \ref{thm:x-exist-bar}.  The proof of Theorem \ref{thm:x-exist-cobar} is essentially dual and therefore left to the reader.  The one, slightly subtle difference in the dual case is that we no longer need  the functor $Y$ to be acyclic on coproducts of models.

According to Theorem \ref{thm:kleisli-x}, we need to prove the existence of a natural transformation
$$\hat\tau:\Phi(\op Q)\underset {\op A}\circ z(X) \Rightarrow z(Y)$$
of functors into the category of left $\op A$-modules that extends $\tau$.  Since the symmetric sequence of graded modules underlying $\Phi(\op Q)\underset {\op A}\circ z(X)$ is  a free $\op A$-module on $\susp \op Q \circ \op A^\bot \circ z(X)$, it is sufficient to define a morphism of symmetric sequences
 $$\hat \tau:\susp \op Q \circ \op A^\bot \circ z(X) \to z(Y)$$
 that we then extend to a morphism of $\op A$-modules
 $$\hat \tau:\op A\circ \susp \op Q \circ \op A^\bot \circ z(X) \to z(Y).$$
If $\hat \tau$ commutes with the differentials, then we can conclude.

Let $\overline\Delta$ denote the reduced, levelwise comultiplication on $\op Q$. It follows from the formula in Remark \ref{rmk:phi-diffl} for $d_{\Phi}$ that the differential $\tilde d_{\Phi}$ in $\Phi (\op Q)\underset {\op A}\circ z(X)$ is specified as follows.  For all $D\in \ob \cat D$, $1\leq m\leq n$, $q\in \op Q(m)$, $\vec n\in I_{n,m}$, and $x_{i}\in X(D)$ for $1\leq i\leq n$,
{\small\begin{align*}
&\tilde d_{\Phi}\big(s^{m-1} q \otimes \alpha_{\vec n}^\sharp \otimes (x_{1}\otimes \cdots \otimes x_{n})\big)\\
&=\pm s^{m-1} dq \otimes \alpha_{\vec n}^\sharp \otimes (x_{1}\otimes \cdots \otimes x_{n})\\
&+s^{m-1} q \otimes \alpha_{\vec n}^\sharp \otimes \sum _{1\leq k\leq n}\pm(x_{1}\otimes \cdots\otimes dx_{k}\otimes \cdots \otimes x_{n})\\
&+\delta_{2}\otimes \sum _{\vec n'+\vec n''=\vec n}\pm \big(s^{m-1} q_{i} \otimes \alpha_{\vec n'}^\sharp \otimes (x_{1}\otimes \cdots \otimes x_{n'})\big)\otimes \big(s^{m-1} q^{i} \otimes \alpha_{\vec n''}^\sharp \otimes (x_{n'+1}\otimes \cdots \otimes x_{n})\big)\\
&+s^{m-1}q \otimes \sum _{0<k<n}\pm \alpha ^\sharp_{\vec n^k}\otimes (x_{1}\otimes \cdots \otimes x_{k}x_{k+1}\otimes \cdots \otimes x_{n}),
\end{align*}}
where the signs are determined by the Koszul rule, $d$ denotes the internal differential in $\op Q$ and in $X(D)$, and $\overline\Delta (q) = q_{i}\otimes q^{i}$ (using Einstein summation notation).

To prepare the proof of the existence of $\hat \tau$, we define a natural, increasing bilfiltration of the symmetric sequence $\susp\op Q \circ\op A^\bot\circ z\big(X(D)\big )$, for every object $D$ of $\cat D$.  Let
$$F^{s,t}(D)=\bigoplus_{k< s}(\susp\op Q\circ \op A^\bot)(k) \otimes X(D)^{\otimes k} \oplus \big((\susp\op Q\circ \op A^\bot)(s) \otimes X(D)^{\otimes s}\big)_{\leq s+t}$$
and
$$F^s(D)=\bigoplus_{k\leq s}(\susp\op Q\circ \op A^\bot)(k) \otimes X(D)^{\otimes k}.$$
Observe that
$$\tilde d_{\Phi}F^{s,t}(D) \subset \op A\circ F^{s, t-1}(D)\; \forall s,t,$$

$$F^{s}(D)=\bigcup_{t\geq 0}F^{s,t}(D)=F^{s+1, t'}(D)\;\forall t'<s,$$
and
$$\susp\op Q \circ\op A^\bot\circ z\big(X(D)\big )=\bigcup_{s\geq 0}F^s(D).$$
Moreover, since $\op Q(1)=R$,
$$F^1(D)=X(D)$$
for all $D\in \ob D$.

We now prove the existence of $\hat \tau$ by induction on $s$ and $t$. We start by setting $\hat \tau_{D}$ equal to $\tau_{D}$ on $F^1(D)$, for all $D\in\ob\cat D$.

Suppose now that for some $s,t$,  a morphism of symmetric sequences of graded modules $\hat \tau_{D}: F^{s,t} (D)\to Y(D)$ has been defined naturally for all $D\in \ob \cat D$, so that its extension to a morphism of left $\op A$-modules is a differential map and so that its restriction to $X(D)$ is exactly $\tau_{D}$.  We now show that $\hat \tau_{D}$ can be naturally extended to $F^{s,t+1}(D)$ for all $D\in \ob \cat D$.

Note that for all $q\in \op Q(k)$ and all $\vec s \in I_{s,k}$,
$$|s^{k-1}q \otimes \alpha^\sharp_{\vec s}|= |q| + s -1,$$
since $\susp (\op Q\circ \op A^\sharp)\cong \susp \op Q \circ \op A^\bot$, and $\op A^\sharp$ is concentrated in degree 0 in every arity.

For all $k\geq 1$, $r\geq 0$ and $D\in \ob \cat D$, let
$$G_{k,r}(D)=\{s^{k-1}q \otimes \alpha^\sharp_{\vec s}\otimes (z_{1}\otimes \cdots \otimes z_{s})\in F^{s,t+1}(D) \mid q\in \op Q(k)_{<r}\}.$$
Note that  $G_{k,0}(D)=\{0\}$, and a natural extension of $\hat \tau_{D}$ over $G_{k,0}(D)$ therefore exists trivially.

For each $k$, suppose that there is some $r_{k}\geq 0$ such that $\hat \tau_{D}$ has been extended from $F^{s,t}(D)$ naturally over $G_{k,r_{k}}(D)$. We now prove that $\hat \tau_{D}$ can then be extended over $G_{k,r_{k}+1}(D)$, for all $k$.

 Let $\mathfrak m_{1},..., \mathfrak m_{s}\in \mathfrak M$ be models with corresponding generators $x_{i}\in X(\mathfrak m_{i})$ for $1\leq i\leq s$ such that
$$\sum_{i=1}^s |x_{i}|=t- r_{k}-s+2.$$
Let $\mathfrak m= \mathfrak m_{1}\coprod \cdots \coprod \mathfrak m_{s}$, and let $B x_{i}$ denote the image of $x_{i}$ under the morphism $X(\mathfrak m_{i})\to X(\mathfrak m)$ induced by the natural map $\mathfrak m_{i}\to \mathfrak m$.

Let $q\in \op Q(k)_{r_{k}}$, and $\vec s\in I_{s,k}$. Since $$|s^{k-1}q \otimes \alpha^\sharp_{\vec s}\otimes (B x_{1}\otimes \cdots \otimes B x_{s})|=t+1,$$
it follows that
$$s^{k-1}q \otimes \alpha^\sharp_{\vec s}\otimes (B x_{1}\otimes \cdots \otimes B x_{s})\in F^{s, t+1}(\mathfrak m).$$
The formula for $\tilde d_{\Phi}$ and the minimality of $r_{k}$ together imply then that
$$\tilde d_{\Phi}\big(s^{k-1}q \otimes \alpha^\sharp_{\vec s}\otimes (B x_{1}\otimes \cdots \otimes B x_{s})\big)\in
F^{s,t}(\mathfrak m)+G_{k,r_{k}},$$
and therefore, by the induction hypothesis,
$$\hat \tau_{\mathfrak m}\Big(\tilde d_{\Phi}\big(s^{k-1}q \otimes \alpha^\sharp_{\vec s}\otimes (B x_{1}\otimes \cdots \otimes B y_{s})\big)\Big)$$
is a well defined element of $Y(\mathfrak m)$, which, moreover, must be a cycle.  Since $Y(\mathfrak m)$ is acyclic, there exists $y\in Y(\mathfrak m)$ such that
$$dy =\hat \tau_{\mathfrak m}\Big(\tilde d_{\Phi}\big(s^{k-1}q \otimes \alpha^\sharp_{\vec s}\otimes (B x_{1}\otimes \cdots \otimes B x_{s})\big)\Big).$$
We set
$$\hat \tau_{\mathfrak m}\big(s^{k-1}q \otimes \alpha^\sharp_{\vec s}\otimes (B x_{1}\otimes \cdots \otimes B x_{s})\big)= y.$$

Let $D\in \ob\cat D$. We can now extend $\hat \tau_{D}$ naturally over $G_{k,r_{k}}(D)$
as follows.
Let $q\in \op Q(k)_{r_{k}}$. For any $D\in \ob \cat D$ and any $z_{1},...,z_{s}\in X(D)$ such that
$$|s^{k-1}q \otimes \alpha^\sharp_{\vec s}\otimes (z_{1}\otimes \cdots \otimes z_{s})|=t+1,$$
let $\xi_{i}:\mathfrak m_{i}\to D$ denote the representing morphism for $z_{i}$, for $1\leq i\leq s$, i.e.,
$$X(\xi_{i})(x_{i})=z_{i}\;\forall 1\leq i\leq s.$$
Let $\mathfrak m =\mathfrak m_{1}\coprod \cdots \coprod \mathfrak m_{s}$, and $\xi=\xi _{1}+\cdots +\xi_{s}: \mathfrak m\to D$.  Set
$$\hat \tau_{D}\big(s^{k-1}q \otimes \alpha^\sharp_{\vec s}\otimes (z_{1}\otimes \cdots \otimes z_{s})\big)= Y(\xi)\Big(\hat \tau_{\mathfrak m}\big(s^{k-1}q \otimes \alpha^\sharp_{\vec s}\otimes (B x_{1}\otimes \cdots \otimes B x_{s})\big)\Big).$$
It is clear that, thus defined, $\hat \tau_{D}$ is natural and commutes with the differential.


\begin{thebibliography}{29}
\expandafter\ifx\csname natexlab\endcsname\relax\def\natexlab#1{#1}\fi
\providecommand{\bibinfo}[2]{#2}
\ifx\xfnm\relax \def\xfnm[#1]{\unskip,\space#1}\fi
\bibitem[{Berger and Moerdijk(2003)}]{berger-moerdijk:03}
\bibinfo{author}{C.~Berger}, \bibinfo{author}{I.~Moerdijk},
  \bibinfo{title}{Axiomatic homotopy theory for operads},
  \bibinfo{journal}{Comment. Math. Helv.} \bibinfo{volume}{78}
  (\bibinfo{year}{2003}) \bibinfo{pages}{805--831}.
\bibitem[{Berger and Moerdijk(2007)}]{berger-moerdijk:07}
\bibinfo{author}{C.~Berger}, \bibinfo{author}{I.~Moerdijk},
  \bibinfo{title}{Resolution of coloured operads and rectification of homotopy
  algebras}, in: \bibinfo{booktitle}{Categories in algebra, geometry and
  mathematical physics}, volume \bibinfo{volume}{431} of
  \textit{\bibinfo{series}{Contemp. Math.}}, \bibinfo{publisher}{Amer. Math.
  Soc.}, \bibinfo{address}{Providence, RI}, \bibinfo{year}{2007}, pp.
  \bibinfo{pages}{31--58}.
\bibitem[{Boyle(2009)}]{boyle}
\bibinfo{author}{M.~Boyle}, \bibinfo{title}{An Algebraic Model for the Homology
  of Pointed Mapping Spaces out of a Closed Surface}, Ph.D. thesis, University
  of Aberdeen, \bibinfo{address}{UK}, \bibinfo{year}{2009}.
\bibitem[{Chevalley and Eilenberg(1948)}]{chevalley-eilenberg:48}
\bibinfo{author}{C.~Chevalley}, \bibinfo{author}{S.~Eilenberg},
  \bibinfo{title}{Cohomology theory of {L}ie groups and {L}ie algebras},
  \bibinfo{journal}{Trans. Amer. Math. Soc.} \bibinfo{volume}{63}
  (\bibinfo{year}{1948}) \bibinfo{pages}{85--124}.
\bibitem[{Fox and Markl(1997)}]{fox-markl:97}
\bibinfo{author}{T.F. Fox}, \bibinfo{author}{M.~Markl},
  \bibinfo{title}{Distributive laws, bialgebras, and cohomology}, in:
  \bibinfo{booktitle}{Operads: Proceedings of Renaissance Conferences
  (Hartford, CT/Luminy, 1995)}, volume \bibinfo{volume}{202} of
  \textit{\bibinfo{series}{Contemp. Math.}}, \bibinfo{publisher}{Amer. Math.
  Soc.}, \bibinfo{address}{Providence, RI}, \bibinfo{year}{1997}, pp.
  \bibinfo{pages}{167--205}.
\bibitem[{Fresse(2000)}]{fresse:00}
\bibinfo{author}{B.~Fresse}, \bibinfo{title}{On the homotopy of simplicial
  algebras over an operad}, \bibinfo{journal}{Trans. Amer. Math. Soc.}
  \bibinfo{volume}{352} (\bibinfo{year}{2000}) \bibinfo{pages}{4113--4141}.
\bibitem[{Fresse(2004)}]{fresse:04}
\bibinfo{author}{B.~Fresse}, \bibinfo{title}{Koszul duality of operads and
  homology of partition posets}, in: \bibinfo{booktitle}{Homotopy theory:
  relations with algebraic geometry, group cohomology, and algebraic
  $K$-theory}, volume \bibinfo{volume}{346} of
  \textit{\bibinfo{series}{Contemp. Math.}}, \bibinfo{publisher}{Amer. Math.
  Soc.}, \bibinfo{address}{Providence, RI}, \bibinfo{year}{2004}, pp.
  \bibinfo{pages}{115--215}.
\bibitem[{Getzler and Jones(1994)}]{getzler-jones:94}
\bibinfo{author}{E.~Getzler}, \bibinfo{author}{J.~Jones},
  \bibinfo{title}{Operads, homotopy algebra, and iterated integrals for double
  loop spaces}, \bibinfo{year}{1994}. \bibinfo{note}{{\tt
  arxiv:hep-th/9403055}}.
\bibitem[{Ginzburg and Kapranov(1994)}]{ginzburg-kapranov:94}
\bibinfo{author}{V.~Ginzburg}, \bibinfo{author}{M.~Kapranov},
  \bibinfo{title}{Koszul duality for operads}, \bibinfo{journal}{Duke Math. J.}
  \bibinfo{volume}{76} (\bibinfo{year}{1994}) \bibinfo{pages}{203--272}.
\bibitem[{Gugenheim and Munkholm(1974)}]{gugenheim-munkholm:74}
\bibinfo{author}{V.K.A.M. Gugenheim}, \bibinfo{author}{H.J. Munkholm},
  \bibinfo{title}{On the extended functoriality of {T}or and {C}otor},
  \bibinfo{journal}{J. Pure Appl. Algebra} \bibinfo{volume}{4}
  (\bibinfo{year}{1974}) \bibinfo{pages}{9--29}.
\bibitem[{Hess and Lack(2010)}]{hess-lack}
\bibinfo{author}{K.~Hess}, \bibinfo{author}{S.~Lack}, \bibinfo{title}{On
  twisting structures}, \bibinfo{year}{2010}. \bibinfo{note}{In preparation}.
\bibitem[{Hess and Levi(2007)}]{hess-levi:07}
\bibinfo{author}{K.~Hess}, \bibinfo{author}{R.~Levi}, \bibinfo{title}{An
  algebraic model for the loop space homology of a homotopy fiber},
  \bibinfo{journal}{Algebr. Geom. Topol.} \bibinfo{volume}{7}
  (\bibinfo{year}{2007}) \bibinfo{pages}{1699--1765}.
\bibitem[{Hess et~al.(2005)Hess, Parent and Scott}]{hps:05}
\bibinfo{author}{K.~Hess}, \bibinfo{author}{P.E. Parent},
  \bibinfo{author}{J.~Scott}, \bibinfo{title}{Co-rings over operads
  characterize morphisms}, \bibinfo{year}{2005}. \bibinfo{note}{{\tt
  arxiv:math.AT/0505559}}.
\bibitem[{Hess et~al.(2007)Hess, Parent and Scott}]{hps:james}
\bibinfo{author}{K.~Hess}, \bibinfo{author}{P.E. Parent},
  \bibinfo{author}{J.~Scott}, \bibinfo{title}{A chain coalgebra model for the
  {J}ames map}, \bibinfo{journal}{Homology, Homotopy Appl.} \bibinfo{volume}{9}
  (\bibinfo{year}{2007}) \bibinfo{pages}{209--231}.
\bibitem[{Hess et~al.(2009)Hess, Parent and Scott}]{hps:cohoch}
\bibinfo{author}{K.~Hess}, \bibinfo{author}{P.E. Parent},
  \bibinfo{author}{J.~Scott}, \bibinfo{title}{Co{H}ochschild homology of chain
  coalgebras}, \bibinfo{journal}{J. Pure Appl. Algebra} \bibinfo{volume}{213}
  (\bibinfo{year}{2009}) \bibinfo{pages}{536--556}.
\bibitem[{Hess et~al.(2006)Hess, Parent, Scott and Tonks}]{hpst:canonicalAH}
\bibinfo{author}{K.~Hess}, \bibinfo{author}{P.E. Parent},
  \bibinfo{author}{J.~Scott}, \bibinfo{author}{A.~Tonks}, \bibinfo{title}{A
  canonical enriched {A}dams-{H}ilton model for simplicial sets},
  \bibinfo{journal}{Adv. Math.} \bibinfo{volume}{207} (\bibinfo{year}{2006})
  \bibinfo{pages}{847--875}.
\bibitem[{Hess and Rognes(2010)}]{hess-rognes}
\bibinfo{author}{K.~Hess}, \bibinfo{author}{J.~Rognes}, \bibinfo{title}{Power
  maps in algebra and topology}, \bibinfo{year}{2010}. \bibinfo{note}{In
  preparation}.
\bibitem[{Iwase and Mimura(1989)}]{iwase-mimura:89}
\bibinfo{author}{N.~Iwase}, \bibinfo{author}{M.~Mimura}, \bibinfo{title}{Higher
  homotopy associativity}, in: \bibinfo{booktitle}{Algebraic topology
  ({A}rcata, {CA}, 1986)}, volume \bibinfo{volume}{1370} of
  \textit{\bibinfo{series}{Lecture Notes in Math.}},
  \bibinfo{publisher}{Springer}, \bibinfo{address}{Berlin},
  \bibinfo{year}{1989}, pp. \bibinfo{pages}{193--220}.
\bibitem[{Johnstone(2002)}]{johnstone:02}
\bibinfo{author}{P.T. Johnstone}, \bibinfo{title}{Sketches of an elephant: a
  topos theory compendium. {V}ol. 1}, volume~\bibinfo{volume}{43} of
  \textit{\bibinfo{series}{Oxford Logic Guides}}, \bibinfo{publisher}{The
  Clarendon Press Oxford University Press}, \bibinfo{address}{New York},
  \bibinfo{year}{2002}.
\bibitem[{Joyal(1981)}]{joyal:81}
\bibinfo{author}{A.~Joyal}, \bibinfo{title}{Une th\'eorie combinatoire des
  s\'eries formelles}, \bibinfo{journal}{Adv. in Math.} \bibinfo{volume}{42}
  (\bibinfo{year}{1981}) \bibinfo{pages}{1--82}.
\bibitem[{Leinster(2000)}]{leinster:00}
\bibinfo{author}{T.~Leinster}, \bibinfo{title}{Homotopy algebras for operads},
  \bibinfo{year}{2000}. \bibinfo{note}{{\tt arxiv:math/0002180v1 [math.QA]}}.
\bibitem[{Markl(2002)}]{markl:02}
\bibinfo{author}{M.~Markl}, \bibinfo{title}{Homotopy diagrams of algebras},
  \bibinfo{journal}{Rend. Circ. Mat. Palermo (2) Suppl.}, in:
  \bibinfo{booktitle}{Proceedings of the 21st {W}inter {S}chool ``{G}eometry
  and {P}hysics'' ({S}rn\'\i, 2001)}, \bibinfo{number}{69}, pp.
  \bibinfo{pages}{161--180}.
\bibitem[{Markl(2004)}]{markl:04}
\bibinfo{author}{M.~Markl}, \bibinfo{title}{Homotopy algebras are homotopy
  algebras}, \bibinfo{journal}{Forum Math.} \bibinfo{volume}{16}
  (\bibinfo{year}{2004}) \bibinfo{pages}{129--160}.
\bibitem[{May(1970)}]{may:70}
\bibinfo{author}{J.P. May}, \bibinfo{title}{A general algebraic approach to
  {S}teenrod operations}, in: \bibinfo{booktitle}{The Steenrod Algebra and its
  Applications (Proc. Conf. to Celebrate N. E. Steenrod's Sixtieth Birthday,
  Battelle Memorial Inst., Columbus, Ohio, 1970)}, Lecture Notes in
  Mathematics, Vol. 168, \bibinfo{publisher}{Springer},
  \bibinfo{address}{Berlin}, \bibinfo{year}{1970}, pp.
  \bibinfo{pages}{153--231}.
\bibitem[{Naito(2009)}]{naito}
\bibinfo{author}{T.~Naito}, \bibinfo{title}{Coalg\`{e}bres
  d'{A}lexander-{W}hitney: un mod\`{e}le alg\'{e}brique pour les espaces
  topologiques}, Ph.D. thesis, \'{E}cole Polytechnique F\'{e}d\'{e}rale de
  Lausanne, \bibinfo{address}{Switzerland}, \bibinfo{year}{2009}.
\bibitem[{Ndombol and Thomas(2004)}]{ndombol-thomas:04}
\bibinfo{author}{B.~Ndombol}, \bibinfo{author}{J.C. Thomas},
  \bibinfo{title}{Steenrod operations in {$shc$}-algebras},
  \bibinfo{journal}{J. Pure Appl. Algebra} \bibinfo{volume}{192}
  (\bibinfo{year}{2004}) \bibinfo{pages}{239--264}.
\bibitem[{Rezk(1996)}]{rezk:96}
\bibinfo{author}{C.W. Rezk}, \bibinfo{title}{Spaces of Algebra Structures and
  Cohomology of Operads}, Ph.D. thesis, Massachusetts Institute of Technology,
  \bibinfo{year}{1996}.
\bibitem[{Smirnov(1981)}]{smirnov:81}
\bibinfo{author}{V.A. Smirnov}, \bibinfo{title}{On the cochain complex of
  topological spaces}, \bibinfo{journal}{Mat. Sb. (N.S.)}
  \bibinfo{volume}{115(157)} (\bibinfo{year}{1981}) \bibinfo{pages}{146--158,
  160}.
\bibitem[{Stasheff(1963)}]{stasheff:63}
\bibinfo{author}{J.D. Stasheff}, \bibinfo{title}{Homotopy associativity of
  {$H$}-spaces. {I}, {II}}, \bibinfo{journal}{Trans. Amer. Math. Soc. 108
  (1963), 275-292; ibid.} \bibinfo{volume}{108} (\bibinfo{year}{1963})
  \bibinfo{pages}{293--312}.

\end{thebibliography}

 \end{document}